\documentclass{article}

\usepackage{a4wide,cite,latexsym,amsfonts,amssymb,exscale,enumerate}
\usepackage[centertags,sumlimits,intlimits,namelimits,reqno]{amsmath}
\usepackage{amsthm}
\usepackage{hyperref}

\usepackage{pstricks}
\usepackage{pst-grad}
\usepackage{pst-plot}
\usepackage[tiling]{pst-fill}
\usepackage{pstcol}

\usepackage{graphics}
\usepackage[all, knot]{xy}
\xyoption{arc}
\usepackage{picins}

\usepackage{fancyheadings}
\def\emph#1{{\sl #1\/}}
\def\ie{{\sl i.e.\/}}
\def\etc{{\sl etc.\/}}
\def\cf{{\sl c.f.\/}}

\def\cal#1{\mathcal{#1}}%
\numberwithin{equation}{section}

\newcommand{\hpeprint}[1]{%
  \href{http://arXiv.org/abs/#1}{\texttt{#1}}}%
\newcommand{\hpmathsci}[1]{%
  \href{http://www.ams.org/mathscinet-getitem?mr=#1}{\texttt{MR #1}}}%

\newcommand{\cat}[1]{\ensuremath{\mbox{\bfseries {\upshape {#1}}}}}
\newcommand{\Hom}{{\rm Hom}}
\newcommand{\End}{{\rm End}}
\newcommand{\id}{{\rm id}}
\newcommand{\im}{{\rm im\ }}
\newcommand{\coim}{{\rm coim\ }}
\newcommand{\coker}{{\rm coker\ }}
\newcommand{\chr}{{\rm char\ }}
\def\ev{\mathrm{ev}}%
\def\coev{\mathrm{coev}}%
\def\tr{\mathrm{tr}}%
\def\st{\mathrm{st}}%
\newcommand{\del}{\partial}
\newcommand{\twocob}{\cat{2Cob}^{{\rm ext}}}
\newcommand{\ontop}[2]{\genfrac{}{}{0pt}{2}{\scriptstyle #1}{\scriptstyle #2}}
\newcommand{\maps}{\colon}
\def\nn{\notag}
\def\1{\mathbbm{1}}%

\usepackage{bbm}
\def\C{{\mathbbm C}}
\def\N{{\mathbbm N}}
\def\R{{\mathbbm R}}

\newcommand{\scs}{\scriptstyle}

\theoremstyle{definition}
\newtheorem{thm}{Theorem}[section]
\newtheorem{cor}[thm]{Corollary}

\newtheorem{lem}[thm]{Lemma}
\newtheorem{rem}[thm]{Remark}
\newtheorem{prop}[thm]{Proposition}
\newtheorem{defn}[thm]{Definition}
\newtheorem{myexample}[thm]{Example}

\let\phi=\varphi
\let\theta=\vartheta
\let\epsilon=\varepsilon

\def\theauthor{}
\def\empty{}
\def\theaffiliation{}
\def\preprint#1{
  \thispagestyle{plain}
  \def\theauthor{#1}
  \ifx\theauthor\empty
  \else
    \begin{flushright}{\small #1\par}\end{flushright}
  \fi
  \begin{center}}
\def\title#1{
  {\LARGE #1\par}\vskip 1em}
\def\author#1{
  \ifx\theaffiliation\empty
  \else
    \par
  \fi
  \def\theauthor{#1}\def\theaffiliation{}}
\def\email#1{
  \vskip 1em{\large\theauthor\footnote{\small email: {\tt #1}}\par}\vskip .5em}
\def\affiliation#1{
  \ifx\theaffiliation\empty
    \def\theaffiliation{second}
  \else
    \par and\par
  \fi
  {\small\sl #1}}
\def\date#1{
  \vskip 1em{(#1)\par}\end{center}\vskip 2em}

\SelectTips{cm}{}

\psset{linewidth=0.3pt,dimen=middle}
\psset{xunit=.70cm,yunit=0.70cm}

\newcommand{\multl}{
  \pscustom[fillcolor=lightgray, fillstyle=solid]{
        \psbezier(1.5,2.5)(1.5,1.1)(.4,1.6)(.5,0)
        \psline(-0.5,0)
        \psbezier(-0.5,0)(-.4,1.6)(-1.5,1.1)(-1.5,2.5)
        \psline(-.5,2.5)
        \psbezier(-.5,2.5)(-.6,1.5)(0.6,1.5)(.5,2.5)
        \psline(1.5,2.5)
    }
}

\newcommand{\comultl}{
  \pscustom[fillcolor=lightgray, fillstyle=solid]{
        \psbezier(1.5,0)(1.5,1.4)(.4,.9)(.5,2.5)
        \psline(-0.5,2.5)
        \psbezier(-0.5,2.5)(-.4,.9)(-1.5,1.4)(-1.5,0)
        \psline(-.5,0)
        \psbezier(-.5,0)(-.6,1)(0.6,1)(.5,0)
        \psline(1.5,0)
    }
}

\newcommand{\ctl}{
  \begin{psclip}{
    \pscustom{
        \psline(-.58,2)(-.58,0)
        \psline(-.58,0)(.42,0)
        \psline(.42,0)(.42,2)
        \psellipse(-.08,2)(.5,.2)
    }
  }
    \pspolygon[fillcolor=lightgray,fillstyle=gradient,
    gradbegin=lightgray,gradend=gray,gradmidpoint=1,gradangle=110](-.58,0)(-.58,2.4)(.42,2.4)(.42,0)(-.58,0)
 \end{psclip}
 \pscustom[fillcolor=lightgray,fillstyle=gradient,
        gradbegin=white, gradend=gray,gradmidpoint=0,gradangle=88]{
    \psline(-.58,2)(-.58,0)
    \psbezier(-.58,0)(-.48,.5)(-.48,.7)(-.08,1)
    \psbezier(-.08,1)(.32,.7)(.32,.5)(.42,0)
    \psellipse(-.08,2)(.5,.2)
 }
 \psellipse[fillcolor=lightgray,fillstyle=gradient,
        gradbegin=lightgray, gradend=gray,gradmidpoint=1,gradangle=110](-.08,2)(.5,.2)
}

\newcommand{\ltc}{
    \pspolygon[fillcolor=lightgray,fillstyle=gradient,
    gradbegin=lightgray,gradend=gray,gradmidpoint=1,gradangle=60](.58,2)(.58,.4)(-.42,.4)(-.42,2)(.58,2)
 \pscustom[fillcolor=lightgray,fillstyle=gradient,
        gradbegin=white, gradend=gray,gradmidpoint=0,gradangle=88]{
    \psline(.58,0)(.58,2)
    \psbezier(.58,2)(.48,1.5)(.48,1.3)(.08,1)
    \psbezier(.08,1)(-.32,1.3)(-.32,1.5)(-.42,2)
    \psline(-.42,0)
    \psbezier(-.42,0)(-.32,-.25)(.48,-.25)(.58,0)
 }
 \begin{psclip}{
 \pspolygon[linestyle=none](.58,0)(.58,.3)(-.42,.3)(-.42,0)(.58,0)
 }
 \psellipse[linestyle=dotted](.08,0)(.5,0.2)
 \end{psclip}
}

 \newcommand{\birthl}{
 \pscustom[fillcolor=lightgray, fillstyle=solid]{
        \psbezier(-.5,0)(-.5,.9)(0.5,.9)(.5,0)
        \psline(-.5,0)
    }
 }

  \newcommand{\deathl}{
 \pscustom[fillcolor=lightgray, fillstyle=solid]{
        \psbezier(-.5,0)(-.5,-.9)(0.5,-.9)(.5,0)
        \psline(-.5,0)    }
 }

\newcommand{\identl}{
    \pspolygon[fillcolor=lightgray,fillstyle=solid](-.5,0)(.5,0)(.5,2.5)(-.5,2.5)(-.5,0)
}

\newcommand{\medidentl}{
    \pspolygon[fillcolor=lightgray,fillstyle=solid](-.5,0)(.5,0)(.5,2)(-.5,2)(-.5,0)
}

\newcommand{\curverightl}{
  \pscustom[fillcolor=lightgray, fillstyle=solid]{
        \psbezier(1.5,2.5)(1.5,1.5)(.4,1.3)(.5,0)
        \psline(-0.5,0)
        \psbezier(-.5,0)(-.6,1.3)(.5,1.5)(.5,2.5)
        \psline(1.5,2.5)
    }
}

\newcommand{\curveleftl}{
  \pscustom[fillcolor=lightgray, fillstyle=solid]{
        \psbezier(-1.5,2.5)(-1.5,1.5)(-.4,1.3)(-.5,0)
        \psline(0.5,0)
        \psbezier(.5,0)(.6,1.3)(-.5,1.5)(-.5,2.5)
        \psline(-1.5,2.5)
    }
}

\newcommand{\multc}{
      \pscustom[fillstyle=gradient,
    gradbegin=white, gradend=gray,gradmidpoint=0,gradangle=70]{
        \psbezier(1.5,2.5)(1.5,1.1)(.4,1.6)(.5,0)
        \psbezier(.5,0)(.4,-.25)(-.4,-.25)(-.5,0)
        \psbezier(-0.5,0)(-.4,1.6)(-1.5,1.1)(-1.5,2.5)
        \psline(-.5,2.5)
        \psbezier(-.5,2.5)(-.6,1.5)(0.6,1.5)(.5,2.5)
        \psline(1.5,2.5)
    }
    \psellipse[fillcolor=lightgray,fillstyle=gradient,
        gradbegin=lightgray, gradend=gray,gradmidpoint=1,gradangle=110](-1,2.5)(.5,.2)
    \psellipse[fillcolor=lightgray,fillstyle=gradient,
        gradbegin=lightgray, gradend=gray,gradmidpoint=1,gradangle=110](1,2.5)(.5,.2)
     \begin{psclip}{
 \pspolygon[linestyle=none](.5,0)(.5,.3)(-.5,.3)(-.5,0)(.5,0)
 }
 \psellipse[linestyle=dotted](0,0)(.5,0.2)
 \end{psclip}
 }

\newcommand{\comultc}{
  \pscustom[fillstyle=gradient,
    gradbegin=white, gradend=gray,gradmidpoint=0,gradangle=110]{
        \psbezier(1.5,0)(1.5,1.4)(.4,.9)(.5,2.5)
        \psline(-0.5,2.5)
        \psbezier(-0.5,2.5)(-.4,.9)(-1.5,1.4)(-1.5,0)
        \psbezier(-1.5,0)(-1.4,-.25)(-.6,-.25)(-.5,0)
        \psbezier(-.5,0)(-.6,1)(0.6,1)(.5,0)
        \psbezier(.5,0)(.6,-.25)(1.4,-.25)(1.5,0)
    }
  \psellipse[fillcolor=lightgray,fillstyle=gradient,
        gradbegin=lightgray, gradend=gray,gradmidpoint=1,gradangle=110](0,2.5)(.5,.2)
\begin{psclip}{
 \pspolygon[linestyle=none](1.5,0)(1.5,.3)(-1.5,.3)(-1.5,0)(1.5,0)
 }
 \psellipse[linestyle=dotted](1,0)(.5,0.2)
 \psellipse[linestyle=dotted](-1,0)(.5,0.2)
 \end{psclip}
 }

\newcommand{\birthc}{
 \pscustom[fillstyle=gradient,
    gradbegin=white, gradend=gray,gradmidpoint=0,gradangle=110]{
        \psbezier(-.5,0)(-.5,.9)(0.5,.9)(.5,0)
        \psbezier(.5,0)(.4,-.25)(-.4,-.25)(-.5,0)
    }
 \begin{psclip}{
 \pspolygon[linestyle=none](.5,0)(.5,.3)(-.5,.3)(-.5,0)(.5,0)
 }
 \psellipse[linestyle=dotted](0,0)(.5,0.2)
 \end{psclip}
 }

\newcommand{\deathc}{
 \pscustom[fillstyle=gradient,
    gradbegin=white, gradend=gray,gradmidpoint=0,gradangle=70]{
        \psbezier(-.5,1)(-.5,.1)(0.5,.1)(.5,1)
        \psline(-.5,1)
 }
  \psellipse[fillcolor=lightgray,fillstyle=gradient,
        gradbegin=lightgray, gradend=gray,gradmidpoint=1,gradangle=110](0,1)(.5,.2)
 }

\newcommand{\zagc}{
   \pscustom[fillstyle=gradient,
    gradbegin=white, gradend=gray,gradmidpoint=0,gradangle=110]{
        \psbezier(1.5,0)(1.6,2)(-1.6,2)(-1.5,0)
        \psbezier(-1.5,0)(-1.4,-.25)(-.6,-.25)(-.5,0)
        \psbezier(-.5,0)(-.6,.8)(0.6,.8)(.5,0)
        \psbezier(.5,0)(.6,-.25)(1.4,-.25)(1.5,0)
    }
  \begin{psclip}{
 \pspolygon[linestyle=none](1.5,0)(1.5,.3)(-1.5,.3)(-1.5,0)(1.5,0)
 }
 \psellipse[linestyle=dotted](1,0)(.5,0.2)
 \psellipse[linestyle=dotted](-1,0)(.5,0.2)
 \end{psclip}
 }

\newcommand{\zigc}{
       \pscustom[fillstyle=gradient,
    gradbegin=white, gradend=gray,gradmidpoint=0,gradangle=70]{
        \psbezier(1.5,2)(1.6,0)(-1.6,0)(-1.5,2)
        \psline(-.5,2)
        \psbezier(-.5,2)(-.6,1.2)(0.6,1.2)(.5,2)
        \psline(1.5,2)
    }
 \psellipse[fillcolor=lightgray,fillstyle=gradient,
        gradbegin=lightgray, gradend=gray,gradmidpoint=1,gradangle=110](1,2)(.5,.2)
        \psellipse[fillcolor=lightgray,fillstyle=gradient,
        gradbegin=lightgray, gradend=gray,gradmidpoint=1,gradangle=110](-1,2)(.5,.2)
}

\newcommand{\identc}{
 \pscustom[fillcolor=lightgray,fillstyle=gradient,
        gradbegin=white, gradend=gray,gradmidpoint=0,gradangle=88]{
 \psline(.5,0)(.5,2.5)
 \psline(-.5,2.5)
 \psline(-.5,0)
 \psbezier(-.5,0)(-.4,-.25)(.4,-.25)(.5,0)
 }
\psellipse[fillcolor=lightgray,fillstyle=gradient,
        gradbegin=lightgray, gradend=gray,gradmidpoint=1,gradangle=110](0,2.5)(.5,.2)
 \begin{psclip}{
 \pspolygon[linestyle=none](.5,0)(.5,.3)(-.5,.3)(-.5,0)(.5,0)
 }
 \psellipse[linestyle=dotted](0,0)(.5,0.2)
 \end{psclip}
 }

\newcommand{\La}{}%
\newcommand{\Lb}{}%
\newcommand{\Lc}{}%
\newcommand{\Lone}[1]{%
  \renewcommand{\La}{#1}}%
\newcommand{\Ltwo}[2]{%
  \renewcommand{\La}{#1}%
  \renewcommand{\Lb}{#2}}%
\newcommand{\Lthree}[3]{%
  \renewcommand{\La}{#1}%
  \renewcommand{\Lb}{#2}%
  \renewcommand{\Lc}{#3}}%

\newcommand{\stringcap}[1]{\xybox{%
  (-6,0)*{};(6,0)*{};**\crv{(6,12)&(-6,12)}?(.2)*\dir{>}?(.8)*\dir{>};
  (8.5,1)*[r]{#1\La};
}}
\newcommand{\stringcup}[1]{\xybox{%
  (-6,0)*{};(6,0)*{};**\crv{(6,-12)&(-6,-12)}?(.2)*\dir{>}?(.8)*\dir{>};
  (8.5,-1)*[r]{#1\La};
}}
\newcommand{\stringidd}[1]{\xybox{%
  (0,0)*{}; (0,-12)*{}; **\dir{-} ?(.5)*\dir{<};
  (2.5,-2)*[r]{#1\La};
  (2.5,-10)*[r]{#1\Lb};
}}
\newcommand{\stringidu}[1]{\xybox{%
  (0,0)*{}; (0,-12)*{}; **\dir{-} ?(.5)*\dir{>};
  (2.5,-2)*[r]{#1\La};
  (2.5,-10)*[r]{#1\Lb};
}}
\newcommand{\stringiddl}[1]{\xybox{%
  (0,0)*{}; (0,-24)*{}; **\dir{-} ?(.5)*\dir{<};
  (2.5,-2)*[r]{#1\La};
  (2.5,-22)*[r]{#1\Lb};
}}
\newcommand{\stringidul}[1]{\xybox{%
  (0,0)*{}; (0,-24)*{}; **\dir{-} ?(.5)*\dir{>};
  (2.5,-2)*[r]{#1\La};
  (2.5,-22)*[r]{#1\Lb};
}}
\newcommand{\stringmapd}[1]{\xybox{%
  (0,0)*{}="t";
  (0,-12)*\xycircle(2.65,2.65){-}="f";
  (0,-24)*{}="b";
  (3,-2)*[r]={#1\La};
  (3,-22)*[r]={#1\Lb};
  (0,-12)*{#1\Lc};
  "t";"f" **\dir{-}; ?(.25)*\dir{>};
  "f";"b" **\dir{-}; ?(.75)*\dir{>};
}}

\newcommand{\stringmapds}[1]{\xybox{%
  (0,0)*{}="t";
  (0,-6)*\xycircle(2.65,2.65){-}="f";
  (0,-12)*{}="b";
  (3,-2)*[r]={#1\La};
  (3,-10)*[r]={#1\Lb};
  (0,-6)*{#1\Lc};
  "t";"f" **\dir{-}; ?(.25)*\dir{>};
  "f";"b" **\dir{-}; ?(.75)*\dir{>};
}}
\newcommand{\stringmapdb}[1]{\xybox{%
  (0,0)*{}="t";
  (0,-12)*\xycircle(4.5,4.5){-}="f";
  (0,-24)*{}="b";
  (5,0)*{};
  (3,-2)*[r]={#1\La};
  (3,-22)*[r]={#1\Lb};
  (0,-12)*{#1\Lc};
  "t";"f" **\dir{-}; ?(.25)*\dir{>};
  "f";"b" **\dir{-}; ?(.75)*\dir{>};
}}
\newcommand{\stringmapdm}[1]{\xybox{%
  (0,0)*{}="t";
  (0,-18)*\xycircle(3.5,3.5){-}="f";
  (0,-36)*{}="b";
  (5,0)*{};
  (3,-2)*[r]={#1\La};
  (3,-34)*[r]={#1\Lb};
  (0,-18)*{#1\Lc};
  "t";"f" **\dir{-}; ?(.25)*\dir{>};
  "f";"b" **\dir{-}; ?(.75)*\dir{>};
}}
\newcommand{\stringbraidd}[1]{\xybox{%
  (-6,0)*{}="t1";
  (6,0)*{}="t2";
  (-6,-12)*{}="b1";
  (6,-12)*{}="b2";
  "t1";"b2" **\crv{(-6,-6)&(6,-6)}; ?(.25)*\dir{>};
  "t2";"b1" **\crv{(6,-6)&(-6,-6)}; ?(.25)*\dir{>};
  (-8.5,-2)*[l]{#1\La};
  (8.5,-2)*[r]{#1\Lb};
}}
\newcommand{\stringbraiddu}[1]{\xybox{%
  (-6,0)*{}="t1";
  (6,0)*{}="t2";
  (-6,-12)*{}="b1";
  (6,-12)*{}="b2";
  "t1";"b2" **\crv{(-6,-6)&(6,-6)}; ?(.25)*\dir{>};
  "t2";"b1" **\crv{(6,-6)&(-6,-6)}; ?(.25)*\dir{<};
  (-8.5,-2)*[l]{#1\La};
  (8.5,-2)*[r]{#1\Lb};
}}
\newcommand{\stringmu}[1]{\xybox{%
  (0,-12)*\xycircle(2.65,2.65){-}="f";
  (-6,0)*{}="t1";
  (6,0)*{}="t2";
  (0,-24)*{}="b";
  "t1";"f" **\crv{(-6,-12)}; ?(.25)*\dir{>};
  "t2";"f" **\crv{(6,-12)}; ?(.25)*\dir{>};
  "f";"b" **\dir{-}; ?(.75)*\dir{>};
  (-3,-2)*[r]{#1\La};
  (9,-2)*[r]{#1\La};
  (3,-22)*[r]{#1\La};
  (0,-12)*{#1\Lb}
}}
\newcommand{\stringeta}[1]{\xybox{%
  (0,0)*{};
  (0,-12)*\xycircle(2.65,2.65){-}="f";
  (0,-24)*{}="b";
  "f";"b" **\dir{-}; ?(.75)*\dir{>};
  (3,-22)*[r]{#1\La};
  (0,-12)*{#1\Lb}
}}
\newcommand{\stringdelta}[1]{\xybox{%
  (0,-12)*\xycircle(2.65,2.65){-}="f";
  (-6,-24)*{}="b1";
  (6,-24)*{}="b2";
  (0,0)*{}="t";
  "f";"b1" **\crv{(-6,-12)}; ?(.75)*\dir{>};
  "f";"b2" **\crv{(6,-12)}; ?(.75)*\dir{>};
  "t";"f" **\dir{-}; ?(.25)*\dir{>};
  (-3,-20)*[r]{#1\La};
  (9,-20)*[r]{#1\La};
  (3,-2)*[r]{#1\La};
  (0,-12)*{#1\Lb}
}}
\newcommand{\stringepsilon}[1]{\xybox{%
  (0,-12)*\xycircle(2.65,2.65){-}="f";
  (0,0)*{}="t";
  "t";"f" **\dir{-}; ?(.25)*\dir{>};
  (3,-2)*[r]{#1\La};
  (0,-12)*{#1\Lb}
}}
\newcommand{\stringpair}[1]{\xybox{%
  (0,-12)*\xycircle(2.65,2.65){-}="f";
  (-6,0)*{}="t1";
  (6,0)*{}="t2";
  "t1";"f" **\crv{(-6,-12)}; ?(.25)*\dir{>};
  "t2";"f" **\crv{(6,-12)}; ?(.25)*\dir{>};
  (-3,-2)*[r]{#1\La};
  (9,-2)*[r]{#1\La};
  (0,-12)*{#1\Lb}
}}
\newcommand{\stringdual}[1]{\xybox{%
  (0,0)*{};
  (0,-12)*\xycircle(2.65,2.65){-}="f";
  (-6,-24)*{}="b1";
  (6,-24)*{}="b2";
  "f";"b1" **\crv{(-6,-12)}; ?(.75)*\dir{>};
  "f";"b2" **\crv{(6,-12)}; ?(.75)*\dir{>};
  (-3,-20)*[r]{#1\La};
  (9,-20)*[r]{#1\La};
  (0,-12)*{#1\Lb}
}}

\newcommand{\siidd}[1]{\xybox{%
  (-3,0)*{};
  (3,0)*{};
  (-1,0);(-1,-12) **\dir{-}; ?(.5)*\dir{>};
  (1,0);(1,-12) **\dir{-}; ?(.5)*\dir{<}
}}
\newcommand{\siidds}[1]{\xybox{%
  (-3,0)*{};
  (3,0)*{};
  (-1,0);(-1,-6) **\dir{-}; ?(.5)*\dir{>};
  (1,0);(1,-6) **\dir{-}; ?(.5)*\dir{<}
}}
\newcommand{\siidu}[1]{\xybox{%
  (-3,0)*{};
  (3,0)*{};
  (-1,0);(-1,-12) **\dir{-}; ?(.5)*\dir{<};
  (1,0);(1,-12) **\dir{-}; ?(.5)*\dir{>}
}}
\newcommand{\siidus}[1]{\xybox{%
  (-3,0)*{};
  (3,0)*{};
  (-1,0);(-1,-6) **\dir{-}; ?(.5)*\dir{<};
  (1,0);(1,-6) **\dir{-}; ?(.5)*\dir{>}
}}
\newcommand{\sidlu}[1]{\xybox{%
  (5,0);(-1,-12) **\crv{(5,-6)&(-1,-6)}; ?(.25)*\dir{<};
  (7,0);(1,-12) **\crv{(7,-6)&(1,-6)}; ?(.25)*\dir{>}
}}
\newcommand{\simu}[1]{\xybox{%
  (-7,0);(-1,-12) **\crv{(-7,-6)&(-1,-6)}; ?(.5)*\dir{>};
  (7,0);(1,-12) **\crv{(7,-6)&(1,-6)}; ?(.5)*\dir{<};
  (-5,0);(5,0) **\crv{(-5,-6)&(5,-6)}; ?(.25)*\dir{<}
}}
\newcommand{\sieta}[1]{\xybox{%
  (-3,0)*{};
  (3,0)*{};
  (-1,-12);(1,-12) **\crv{(-1,-2)&(1,-2)}; ?(.15)*\dir{<}
}}
\newcommand{\sidelta}[1]{\xybox{%
  (-7,-12);(-1,0) **\crv{(-7,-6)&(-1,-6)}; ?(.5)*\dir{<};
  (7,-12);(1,0) **\crv{(7,-6)&(1,-6)}; ?(.5)*\dir{>};
  (-5,-12);(5,-12) **\crv{(-5,-9)&(3,-9)&(3,-5)&(-3,-5)&(-3,-9)&(5,-9)}; ?(.1)*\dir{>}
}}
\newcommand{\siepsilon}[1]{\xybox{%
  (-3,0)*{};
  (3,0)*{};
  (-1,0);(1,0) **\crv{(-1,-5)&(1,-5)&(1,-10)&(-1,-10)&(-1,-5)&(1,-5)}; ?(.1)*\dir{>}
}}
\newcommand{\sicup}[1]{\xybox{%
  (-5,0);(5,0) **\crv{(-5,-9)&(5,-9)}; ?(.25)*\dir{<};
  (-7,0);(7,0) **\crv{(-7,-12)&(7,-12)}; ?(.25)*\dir{>}
}}
\newcommand{\sicap}[1]{\xybox{%
  (0,0)*{};
  (-7,-12);(7,-12) **\crv{(-7,0)&(7,0)}; ?(.25)*\dir{<};
  (-5,-12);(5,-12) **\crv{(-5,-3)&(5,-3)}; ?(.25)*\dir{>}
}}
\newcommand{\sipair}[1]{\xybox{%
  (-5,0);(5,0) **\crv{(-5,-5)&(5,-5)}; ?(.25)*\dir{<};
  (-7,0);(7,0) **\crv{(-7,-7)&(3,-7)&(3,-12)&(-3,-12)&(-3,-7)&(7,-7)}; ?(.1)*\dir{>}
}}
\newcommand{\sidual}[1]{\xybox{%
  (-7,-12);(7,-12) **\crv{(-7,0)&(7,0)}; ?(.25)*\dir{<};
  (-5,-12);(5,-12) **\crv{(-5,-9)&(3,-9)&(3,-4)&(-3,-4)&(-3,-9)&(5,-9)}; ?(.1)*\dir{>}
}}
\newcommand{\sibraiddu}[1]{\xybox{%
  (-7,0);(5,-12) **\crv{(-7,-7)&(5,-7)}; ?(.25)*\dir{>};
  (-5,0);(7,-12) **\crv{(-5,-5)&(7,-5)}; ?(.25)*\dir{<};
  (5,0);(-7,-12) **\crv{(5,-5)&(-7,-5)}; ?(.25)*\dir{<};
  (7,0);(-5,-12) **\crv{(7,-7)&(-5,-7)}; ?(.25)*\dir{>}
}}

\newcommand{\bbmu}[1]{\xybox{%
  (-3,0)*{};
  (3,0)*{};
  (0,-4)*{\bullet}="f";
  (-2,0)*{}="t1";
  (2,0)*{}="t2";
  (0,-8)*{}="b";
  "t1";"f" **\crv{(-2,-2)}; ?(.35)*\dir{>};
  "t2";"f" **\crv{(2,-2)}; ?(.35)*\dir{>};
  "f";"b" **\dir{-}; ?(.75)*\dir{>};
}}
\newcommand{\bbmediummu}[1]{\xybox{%
  (-5,0)*{};
  (5,0)*{};
  (0,-4)*{\bullet}="f";
  (-4,0)*{}="t1";
  (4,0)*{}="t2";
  (0,-8)*{}="b";
  "t1";"f" **\crv{(-4,-2)}; ?(.35)*\dir{>};
  "t2";"f" **\crv{(4,-2)}; ?(.35)*\dir{>};
  "f";"b" **\dir{-}; ?(.75)*\dir{>};
}}
\newcommand{\bbeta}[1]{\xybox{%
  (-1,0)*{};
  (1,0)*{};
  (0,-4)*{\bullet}="f";
  (0,-8)*{}="b";
  "f";"b" **\dir{-}; ?(.75)*\dir{>};
}}
\newcommand{\bbstickball}[1]{\xybox{%
  (-1,0)*{};
  (1,0)*{};
  (0,-5.5)*\xycircle(2.5,2.5){-}="f";
  (0,-5.5)*{\scriptstyle #1};
}}
\newcommand{\bbelement}[1]{\xybox{%
  (-1,0)*{};
  (1,0)*{};
  (0,-4)*\xycircle(2.5,2.5){-}="f";
  (0,-4)*{\scriptstyle #1};
  (0,-8)*{}="b";
  "f";"b" **\dir{-}; ?(.9)*\dir{>};
}}
\newcommand{\bbcardy}[1]{\xybox{%
  (-2,0)*{};
  (2,0)*{};
  (0,-4)*{\bullet}="g";
  (0,4)*{\bullet}="f";
  (0,8)*{}="t";
  (0,-8)*{}="b";
  "t";"f" **\dir{-}; ?(.25)*\dir{>};
  "g";"b" **\dir{-}; ?(.75)*\dir{>};
  "f";"g" **\crv{(-6,3)&(6,-3)}; ?(.4)*\dir{>};
  "f";"g" **\crv{(6,3)&(-6,-3)}; ?(.4)*\dir{>};
}}

\newcommand{\bbainv}[1]{\xybox{%
    (-4,2)*\xycircle(2.9,2.9){-}="t1";
    (-4,2)*{\scriptstyle a^{\scriptscriptstyle -#1}};
  (-4,0)*{};
  (4,0)*{};
  (0,-4)*{\bullet}="f";
  (2,0)*{}="t2";
  (0,-8)*{}="b";
  "t1";"f" **\crv{(-4,-2)}; ?(.65)*\dir{>};
  "t2";"f" **\crv{(2,-2)}; ?(.35)*\dir{>};
  "f";"b" **\dir{-}; ?(.75)*\dir{>};
    }}
    \newcommand{\bba}[1]{\xybox{%
    (-2,8)*\bbstickball{a^{#1}};
    (0,0)*\bbmu{};
    }}
\newcommand{\bbdelta}[1]{\xybox{%
  (-3,0)*{};
  (3,0)*{};
  (0,-4)*{\bullet}="f";
  (-2,-8)*{}="b1";
  (2,-8)*{}="b2";
  (0,-0)*{}="t";
  "f";"b1" **\crv{(-2,-6)}; ?(.75)*\dir{>};
  "f";"b2" **\crv{(2,-6)}; ?(.75)*\dir{>};
  "t";"f" **\dir{-}; ?(.35)*\dir{>};
}}
\newcommand{\bbmeddelta}[1]{\xybox{%
  (-5,0)*{};
  (5,0)*{};
  (0,-4)*{\bullet}="f";
  (-4,-8)*{}="b1";
  (4,-8)*{}="b2";
  (0,-0)*{}="t";
  "f";"b1" **\crv{(-4,-6)}; ?(.75)*\dir{>};
  "f";"b2" **\crv{(4,-6)}; ?(.75)*\dir{>};
  "t";"f" **\dir{-}; ?(.35)*\dir{>};
}}
\newcommand{\bbepsilon}[1]{\xybox{%
  (-1,0)*{};
  (1,0)*{};
  (0,-4)*{\bullet}="f";
  (0,-0)*{}="t";
  "t";"f" **\dir{-}; ?(.35)*\dir{>};
}}
\newcommand{\bbepsilonp}[1]{\xybox{%
  (-1,0)*{};
  (1,0)*{};
  (0,-4)*{\ast}="f";
  (0,-0)*{}="t";
  "t";"f" **\dir{-}; ?(.35)*\dir{>};
}}
\newcommand{\bbid}[1]{\xybox{%
  (-1,0)*{};
  (1,0)*{};
  (0,0);(0,-9) **\dir{-}; ?(.5)*\dir{>}
}}
\newcommand{\bbidu}[1]{\xybox{%
  (-1,0)*{};
  (1,0)*{};
  (0,0);(0,-9) **\dir{-}; ?(.5)*\dir{<}
}}
\newcommand{\bbdl}[1]{\xybox{%
  (2,0);(0,-8) **\crv{(2,-2)&(0,-6)}; ?(.5)*\dir{>}
}}
\newcommand{\bbdlu}[1]{\xybox{%
  (2,0);(0,-8) **\crv{(2,-2)&(0,-6)}; ?(.5)*\dir{<}
}}
\newcommand{\bbdr}[1]{\xybox{%
  (-2,0);(0,-8) **\crv{(-2,-2)&(0,-6)}; ?(.5)*\dir{>}
}}
\newcommand{\bbdru}[1]{\xybox{%
  (-2,0);(0,-8) **\crv{(-2,-2)&(0,-6)}; ?(.5)*\dir{<}
}}
\newcommand{\bbbraid}[1]{\xybox{%
  (-3,0)*{};
  (3,0)*{};
  (-2,0);(2,-8) **\crv{(-2,-2)&(2,-6)}; ?(.25)*\dir{>};
  (2,0);(-2,-8) **\crv{(2,-2)&(-2,-6)}; ?(.25)*\dir{>};
}}
\newcommand{\bbbraiddu}[1]{\xybox{%
  (-3,0)*{};
  (3,0)*{};
  (-2,0);(2,-8) **\crv{(-2,-2)&(2,-6)}; ?(.25)*\dir{>};
  (2,0);(-2,-8) **\crv{(2,-2)&(-2,-6)}; ?(.25)*\dir{<};
}}
\newcommand{\bbbraidr}[1]{\xybox{%
  (-5,0)*{};
  (5,0)*{};
  (-4,0);(4,-8) **\crv{(-4,-2)&(4,-6)}; ?(.25)*\dir{>};
  (0,0);(-4,-8) **\crv{(0,-2)&(-4,-6)}; ?(.25)*\dir{>};
  (4,0);(0,-8) **\crv{(4,-2)&(0,-6)}; ?(.25)*\dir{>};
}}
\newcommand{\bbbraidl}[1]{\xybox{%
  (-5,0)*{};
  (5,0)*{};
  (4,0);(-4,-8) **\crv{(4,-2)&(-4,-6)}; ?(.25)*\dir{>};
  (0,0);(4,-8) **\crv{(0,-2)&(4,-6)}; ?(.25)*\dir{>};
  (-4,0);(0,-8) **\crv{(-4,-2)&(0,-6)}; ?(.25)*\dir{>};
}}
\newcommand{\bbpair}[1]{\xybox{%
  (-3,0)*{};
  (3,0)*{};
  (0,-4)*{\bullet}="f";
  (-2,0)*{}="t1";
  (2,0)*{}="t2";
  "t1";"f" **\crv{(-2,-2)}; ?(.35)*\dir{>};
  "t2";"f" **\crv{(2,-2)}; ?(.35)*\dir{>};
}}
\newcommand{\bbpairp}[1]{\xybox{%
  (-3,0)*{};
  (3,0)*{};
  (0,-4)*{\ast}="f";
  (-2,0)*{}="t1";
  (2,0)*{}="t2";
  "t1";"f" **\crv{(-2,-2)}; ?(.35)*\dir{>};
  "t2";"f" **\crv{(2,-2)}; ?(.35)*\dir{>};
}}
\newcommand{\bbtriple}[1]{\xybox{%
  (-5,0)*{};
  (5,0)*{};
  (0,-4)*{\bullet}="f";
  (-4,0)*{}="t1";
  (0,0)*{}="t2";
  (4,0)*{}="t3";
  "t1";"f" **\crv{(-4,-2)}; ?(.35)*\dir{>};
  "t2";"f" **\dir{-}; ?(.35)*\dir{>};
  "t3";"f" **\crv{(4,-2)}; ?(.35)*\dir{>};
}}
\newcommand{\bbdual}[1]{\xybox{%
  (-3,0)*{};
  (3,0)*{};
  (0,-4)*{\bullet}="f";
  (-2,-8)*{}="b1";
  (2,-8)*{}="b2";
  "f";"b1" **\crv{(-2,-6)}; ?(.75)*\dir{>};
  "f";"b2" **\crv{(2,-6)}; ?(.75)*\dir{>};
}}
\newcommand{\bbdualp}[1]{\xybox{%
  (-3,0)*{};
  (3,0)*{};
  (0,-4)*{\ast}="f";
  (-2,-8)*{}="b1";
  (2,-8)*{}="b2";
  "f";"b1" **\crv{(-2,-6)}; ?(.75)*\dir{>};
  "f";"b2" **\crv{(2,-6)}; ?(.75)*\dir{>};
}}
\newcommand{\bbwidedual}[1]{\xybox{%
  (-5,0)*{};
  (5,0)*{};
  (0,-4)*{\bullet}="f";
  (-4,-8)*{}="b1";
  (4,-8)*{}="b2";
  "f";"b1" **\crv{(-4,-6)}; ?(.75)*\dir{>};
  "f";"b2" **\crv{(4,-6)}; ?(.75)*\dir{>};
}}
\newcommand{\bbhugedual}[1]{\xybox{%
  (-7,0)*{};
  (7,0)*{};
  (0,-4)*{\bullet}="f";
  (-6,-8)*{}="b1";
  (6,-8)*{}="b2";
  "f";"b1" **\crv{(-6,-6)}; ?(.75)*\dir{>};
  "f";"b2" **\crv{(6,-6)}; ?(.75)*\dir{>};
}}
\newcommand{\bbcup}[1]{\xybox{%
  (-3,0)*{};
  (3,0)*{};
  (-2,0);(2,0) **\crv{(-2,-4)&(2,-4)}; ?(.25)*\dir{<};
}}
\newcommand{\bbcap}[1]{\xybox{%
  (-3,0)*{};
  (3,0)*{};
  (-2,-8);(2,-8) **\crv{(-2,-4)&(2,-4)}; ?(.25)*\dir{<};
}}

\newcommand{\bblr}[1]{\xybox{%
  (-3,0)*{};
  (3,0)*{};
  (-2,0);(2,-8) **\crv{(-2,-2)&(2,-6)}; ?(.25)*\dir{>};
}}
\newcommand{\bbrl}[1]{\xybox{%
  (-3,0)*{};
  (3,0)*{};
  (2,0);(-2,-8) **\crv{(2,-2)&(-2,-6)}; ?(.25)*\dir{>};
}}
\newcommand{\bbdouble}[1]{\xybox{%
  (-5,0)*{};
  (5,0)*{};
  (0,-4)*{\bullet}="f";
  (-4,0)*{}="t1";
  (4,0)*{}="t3";
  "t1";"f" **\crv{(-4,-2)}; ?(.35)*\dir{>};
  "t3";"f" **\crv{(4,-2)}; ?(.35)*\dir{>};
}}

\newcommand{\bbrllong}[1]{\xybox{%
  (-5,0)*{};
  (5,0)*{};
  (-4,0);(4,-8) **\crv{(-4,-2)&(4,-6)}; ?(.25)*\dir{>};
 }}

\newcommand{\bblrlong}[1]{
  \xybox{%
  (-5,0)*{};
  (5,0)*{};
  (4,0);(-4,-8) **\crv{(4,-2)&(-4,-6)}; ?(.25)*\dir{>};
  }}
\newcommand{\bbwidemu}[1]{\xybox{%
  (-5,0)*{};
  (5,0)*{};
  (0,-4)*{\bullet}="f";
  (-6,0)*{}="t1";
  (0,-8)*{}="t2";
  (6,0)*{}="t3";
  "t1";"f" **\crv{(-6,-2)}; ?(.45)*\dir{>};
  "t3";"f" **\crv{(6,-2)}; ?(.45)*\dir{>};
  "f";"t2" **\dir{-}; ?(.65)*\dir{>};
}}
\newcommand{\bbmedmu}[1]{\xybox{%
  (-5,0)*{};
  (5,0)*{};
  (0,-4)*{\bullet}="f";
  (-4,0)*{}="t1";
  (0,-8)*{}="t2";
  (4,0)*{}="t3";
  "t1";"f" **\crv{(-4,-2)}; ?(.45)*\dir{>};
  "t3";"f" **\crv{(4,-2)}; ?(.45)*\dir{>};
  "f";"t2" **\dir{-}; ?(.65)*\dir{>};
}}

\newcommand{\bbproject}[1]{\xybox{%
  (-3,0)*{};
  (5,0)*{};
    (2,0)*\bbmu{};
    (0,8)*\bbmu{};
    (4,8)*\bbdl{};
    (0,16)*\bbbraid{};
    (6,16)*\bbid{};
    (-2,24)*\bbid{};
    (4,24)*\bbdual{};
}}

\newcommand{\bbsmallproject}[1]{\xybox{%
  (-4,0)*{};
  (4,0)*{};
  (0,-4)*{\bullet}="f";
  (-2,0)*{\bullet}="t1";
  (2,0)*{}="t2";
  (0,-8)*{}="b";
  "t1";"f" **\crv{(-2,-2)}; ?(.55)*\dir{>};
  "t2";"f" **\crv{(2,-2)}; ?(.55)*\dir{>};
  "f";"b" **\dir{-}; ?(.75)*\dir{>};
  (-2,8)*{}="t";
  (1,5)*{\bullet}="m";
  "t";"t1" **\dir{-}; ?(.35)*\dir{>};
  "t1";"m" **\crv{(-8,2)&(1,4.5)}; ?(.68)*\dir{<};
   "t2";"m" **\crv{(2,3)}; ?(.3)*\dir{<};
}}
\newcommand{\bbP}[1]{\xybox{%
   (0,0)*\xycircle(2.15,2.15){-}="m";
   (0,-6)*{}="b";
   (0,6)*{}="t";
   (0,0)*{p};
         "b";"m" **\dir{-} ?(.3)*\dir{<}  ;
         "t";"m" **\dir{-} ?(.3)*\dir{>}  ;
}}
\newcommand{\bblongid}[1]{\xybox{%
  (-1,0)*{};
  (1,0)*{};
  (0,6);(0,-6) **\dir{-}; ?(.5)*\dir{>}
}}

\begin{document}
%

\preprint{AEI-2005-163\\ DAMTP-2005-81}

\title{State sum construction of two-dimensional open-closed\\
       Topological Quantum Field Theories}

\author{Aaron D.\ Lauda}
\email{A.Lauda@dpmms.cam.ac.uk}
\affiliation{Department of Pure Mathematics and Mathematical Statistics,\\
  University of Cambridge, Cambridge CB3 0WB, United Kingdom}

\author{Hendryk Pfeiffer}
\email{pfeiffer@aei.mpg.de}
\affiliation{Max Planck Institute for Gravitational Physics,\\
  Am M\"uhlenberg 1, 14476 Potsdam, Germany}

\date{June 29, 2006}

%
\begin{abstract}
%

We present a state sum construction of two-dimensional extended
Topological Quantum Field Theories (TQFTs), so-called
\emph{open-closed} TQFTs, which generalizes the state sum of
Fukuma--Hosono--Kawai from triangulations of conventional
two-dimensional cobordisms to those of open-closed cobordisms,
\ie\ smooth compact oriented $2$-manifolds with corners that have
a particular global structure. This construction reveals the
topological interpretation of the associative algebra on which the
state sum is based, as the vector space that the TQFT assigns to
the unit interval. Extending the notion of a two-dimensional TQFT
from cobordisms to suitable manifolds with corners therefore makes
the relationship between the global description of the TQFT in
terms of a functor into the category of vector spaces and the
local description in terms of a state sum fully transparent. We
also illustrate the state sum construction of an open-closed TQFT
with a finite set of D-branes using the example of the groupoid
algebra of a finite groupoid.

\end{abstract}

\noindent
\begin{small}
Mathematics Subject Classification (2000): 57R56, 57M99, 81T40, 57Q20.
\end{small}

%
\section{Introduction}
%

An $n$-dimensional Topological Quantum Field Theory
(TQFT)~\cite{Atiyah} is a symmetric monoidal functor from the category
$\cat{nCob}$ of $n$-dimensional cobordisms to the category
$\cat{Vect}_k$ of vector spaces over a given field $k$.  The objects
of the category $\cat{nCob}$ are smooth compact oriented
$(n-1)$-manifolds without boundary, and the morphisms are equivalence
classes of smooth compact oriented cobordisms between these, modulo
diffeomorphisms that restrict to the identity on the boundary. An
$n$-dimensional TQFT therefore associates vector spaces with
$(n-1)$-manifolds and linear maps with $n$-dimensional cobordisms.
Disjoint unions of manifolds correspond to tensor products of vector
spaces and linear maps, and gluing cobordisms along their boundaries
corresponds to the composition of linear maps. Note that the empty
$(n-1)$-manifold plays the role of the unit object for the tensor
product and corresponds to the field $k$.

For $n=2$, the category $\cat{nCob}$ is well understood, and so
there are strong results about $2$-dimensional TQFTs. For these
classic results, we refer to~\cite{Dij1,Abrams1,Sawin95} and to
the book~\cite{Kock}. It is known, for example, that
$2$-dimensional TQFTs are characterized by commutative Frobenius
algebras. The objects of $\cat{2Cob}$ are compact $1$-manifolds
without boundary, \ie\ disjoint unions of circles $S^1$. For the
morphisms of $\cat{2Cob}$, one has a description in terms of
generators and relations. The generators are these cobordisms:
\begin{equation}
\label{eq_closedgen}
\begin{aligned}
\psset{xunit=.4cm,yunit=.4cm}
\begin{pspicture}[.2](12,3.5)
  \rput(0,1){\multc}
  \rput(4,1){\comultc}
  \rput(8,2){\birthc}
  \rput(11,1.6){\deathc}
  \rput(0,0){$\mu$}
  \rput(4,0){$\Delta$}
  \rput(8,0){$\eta$}
  \rput(11,0){$\epsilon$}
\end{pspicture}
\end{aligned}
\end{equation}
We have drawn them in such a way that their source is at the top and
their target at the bottom of the diagram. The TQFT is a functor
$Z\colon\cat{2Cob}\to\cat{Vect}_k$. If we denote by $C:=Z(S^1)$ the
vector space associated with the circle, the TQFT assigns linear maps
$\mu\colon C\otimes C\to C$, $\Delta\colon C\to C\otimes C$,
$\eta\colon k\to C$ and $\epsilon\colon C\to k$ to the morphisms
depicted in~\eqref{eq_closedgen}. The relations among the morphisms of
$\cat{2Cob}$ then imply that $(C,\mu,\eta,\Delta,\epsilon)$ forms a
commutative Frobenius algebra. Conversely, given any commutative
Frobenius algebra $C$, there is a functor
$Z\colon\cat{2Cob}\to\cat{Vect}_k$ such that $Z(S^1)=C$.  We say that
the commutative Frobenius algebra provides a
\emph{global description} of the $2$-dimensional TQFT. The
relevant algebraic structure, namely the commutative Frobenius
algebra $(C,\mu,\eta,\Delta,\epsilon)$, has an immediate
topological interpretation in terms of the vector space $C$
associated with the circle, the linear maps $\mu$, $\eta$,
$\Delta$, and $\epsilon$ associated with the
generators~\eqref{eq_closedgen}, and in terms of the relations
among the morphisms of $\cat{2Cob}$.

The state sum of Fukuma--Hosono--Kawai~\cite{FHK} forms a different
and \emph{a priori} independent way of defining a $2$-dimensional
TQFT. This construction starts with a finite-dimensional semisimple
algebra $(A,\mu,\eta)$ over a field $k$ of characteristic zero. For
every $2$-dimensional cobordism $M\colon\Sigma_1\to\Sigma_2$, one
considers a triangulation of $M$, and from the data $A$, $\mu$,
$\eta$, and from the triangulation, one computes the linear map
$Z(M)\colon Z(\Sigma_1)\to Z(\Sigma_2)$ as a \emph{state sum}. In a
state sum, roughly speaking, one colours the simplices of the
triangulated manifold $M$ with algebraic data such as the vector space
underlying $A$ or the linear maps $\mu$, $\eta$, and then one `sums
over all colourings' following certain rules. We present this
construction in detail in Section~\ref{sect_statesum} below.  In
particular, one can compute the vector space associated with the
circle, and it turns out that this is the centre
\begin{equation}
\label{eq_centre}
  Z(S^1) = Z(A)
\end{equation}
of the algebra one has started with. The first `Z'
in~\eqref{eq_centre} refers to the functor
$Z\colon\cat{2Cob}\to\cat{Vect}_k$ whereas the second `Z' means
centre. The structure of $Z(S^1)$ as a commutative Frobenius algebra
can be computed from the algebra $A$, too.

While the centre $Z(A)$ has a clear topological interpretation as
outlined above, the algebra $A$ is so far just part of a `recipe'
(the state sum construction), but it is far from obvious whether
$A$ itself plays any role in the topology of $2$-manifolds.

Given a $2$-dimensional TQFT $Z\colon\cat{2Cob}\to\cat{Vect}_k$
where $k$ is a field of characteristic zero, one can ask the
converse question, namely, whether there is a finite-dimensional
semisimple algebra $A$ over $k$ such that one can obtain the given
TQFT from the state sum of Fukuma--Hosono--Kawai. Of course, the
algebra structure of $A$ needs to be such that $Z(A)=Z(S^1)$, but
one also has to understand which Frobenius algebra structure to
choose for $A$ in order to recover the appropriate one for $Z(A)$.
In order to answer this question, a topological interpretation of
the algebra $A$ is clearly desirable.

In the present article, we extend the category $\cat{2Cob}$ from
ordinary cobordisms to the category $\cat{2Cob}^\mathrm{ext}$ of
\emph{open-closed cobordisms}. These are certain smooth
$2$-manifolds with corners that can be viewed as cobordisms between
compact $1$-manifolds \emph{with} boundary, \ie\ between disjoint
unions of circles $S^1$ and unit intervals $I=[0,1]$. We generalize
the notion of a TQFT and the state sum of Fukuma--Hosono--Kawai
accordingly, and we show that the algebra $A$ for the state sum
construction has a topological interpretation as the vector space
associated with the unit interval.

The description of the category $\cat{2Cob}^\mathrm{ext}$ in terms
of generators and relations goes back to work on boundary
conformal field theory by Cardy and Lewellen~\cite{CL,Lew},
Lazaroiu~\cite{lar2}, and to the work of Moore and Segal~\cite{MS}
and Alexeevski and Natanzon~\cite{AN}. We have shown the
sufficiency of the relations in~\cite{LP05}. In order to get some
intuition for the extended cobordism category
$\cat{2Cob}^\mathrm{ext}$, we here display the generators for its
morphisms:
\begin{equation}
\label{eq_generators}
\begin{aligned}
\psset{xunit=.4cm,yunit=.4cm}
\begin{pspicture}[.2](26,3.5)
  \rput(0,1){\multl}
  \rput(0,0){$\mu_A$}
  \rput(4,1){\comultl}
  \rput(4,0){$\Delta_A$}
  \rput(7,2){\birthl}
  \rput(7,0){$\eta_A$}
  \rput(9,2.4){\deathl}
  \rput(9,0){$\epsilon_A$}
  \rput(12,1){\multc}
  \rput(12,0){$\mu_C$}
  \rput(16,1){\comultc}
  \rput(16,0){$\Delta_C$}
  \rput(19,2){\birthc}
  \rput(19,0){$\eta_C$}
  \rput(21,1.6){\deathc}
  \rput(21,0){$\epsilon_C$}
  \rput(23,1){\ctl}
  \rput(23,0){$\imath$}
  \rput(26,1){\ltc}
  \rput(26,0){$\imath^\ast$}
\end{pspicture}
\end{aligned}
\end{equation}
An \emph{open-closed TQFT} is a symmetric monoidal functor
$Z\colon\cat{2Cob}^\mathrm{ext}\to\cat{Vect}_k$. In~\cite{LP05},
we have shown that open-closed TQFTs are characterized by what we
call \emph{knowledgeable Frobenius algebras}
$(A,C,\imath,\imath^\ast)$ where the vector space $C:=Z(S^1)$
associated with the circle has the structure of a commutative
Frobenius algebra $(C,\mu_C,\eta_C,\Delta_C,\epsilon_C)$, the
vector space $A:=Z(I)$ associated with the interval has the
structure of a symmetric Frobenius algebra
$(A,\mu_A,\eta_A,\Delta_A,\epsilon_A)$, and there are linear maps
$\imath\colon C\to A$ and $\imath^\ast\colon A\to C$ subject to
certain conditions. For the details, see Section~\ref{sect_kfrob}
below.

In the present article, we show (see Theorem~\ref{thm_statesum} below)
that for every strongly separable\footnote{It turns out that for a
field of arbitrary characteristic, the appropriate class of algebras
is that of the strongly separable ones.} algebra $A$ over any field
$k$ and for every choice of a symmetric Frobenius algebra structure
for $A$, there is a knowledgeable Frobenius algebra
$(A,Z(A),\imath,\imath^\ast)$ and a generalization to
$\cat{2Cob}^\mathrm{ext}$ of the state sum of Fukuma--Hosono--Kawai
that yields the open-closed TQFT characterized by
$(A,Z(A),\imath,\imath^\ast)$. Extending the notion of a
$2$-dimensional TQFT to suitable manifolds with corners therefore
reveals which topological role is played by the algebra $A$ that
enters the state sum construction.

Why is it important to better understand the role of the algebra $A$?
After all, $2$-dimensional TQFTs are well understood, and the state
sum of Fukuma--Hosono--Kawai is just one of several ways of finding
examples. This question has various answers depending on the view
point taken.

Open-closed cobordisms have a natural string theoretic
interpretation. Indeed, the generators in \eqref{eq_generators} can be
interpreted as the smooth manifolds with corners that underly open and
closed string worldsheets. In the state sum, $A$ therefore turns out
to be the algebra associated to the open string.

The state sum of Fukuma--Hosono--Kawai is also relevant to recent work
on boundary conformal field theory, see, for example~\cite{FRS,RFFS}
where the algebra $A$ already appears in connection with the boundary
conditions, and so the present article is immediately relevant in this
context.

Another reason for better understanding the topological significance
of the algebra $A$ is given by attempts to generalize the framework to
higher dimensions. For $n\geq 3$, the cobordism category $\cat{nCob}$
is not fully understood, \ie\ $n$-dimensional cobordisms have not been
(or even cannot be) classified, and in particular one does not have
any description of $\cat{nCob}$ in terms of generators and
relations. This makes a full understanding of $n$-dimensional TQFTs
much harder if not impossible.

On the other hand, there are some generalizations of the state sum
construction of Fukuma--Hosono--Kawai to higher dimensions, notably
the $3$-dimensional TQFT of Turaev and Viro~\cite{TV}, extended by
Barrett and Westbury~\cite{BW}, which produces a $3$-dimensional TQFT
for any given modular category or, more generally, for suitable
spherical categories~\cite{BW2}. The step from dimension $2$ to $3$,
\ie\ from the state sum of Fukuma--Hosono--Kawai to that of
Turaev--Viro, can be understood as an example of
\emph{categorification} which means replacing algebraic structures
based on sets and maps by analogues that are rather based on
categories and functors~\cite{Jeff}. The dimensional ladder of Crane
and Frenkel~\cite{CF} sketches which sort of algebraic structures one
would need in order to construct $n$-dimensional TQFTs from state
sums:
\begin{equation}
\label{eq_ladder}
\begin{aligned}
\xymatrix@!C=2.1pc@R=1pc{
  n=4&&\mbox{trialgebras}\ar@{-}[ddr]
     &&\mbox{Hopf categories}\ar@{-}[ddl]\ar@{-}[ddr]
     &&\mbox{monoidal $2$-categories}\ar@{-}[ddl]\\
  \\
  n=3&&&\mbox{Hopf algebras}\ar@{-}[ddr]
     &&\mbox{monoidal categories}\ar@{-}[ddl]\\
  \\
  n=2&&&&\mbox{associative algebras}
}
\end{aligned}
\end{equation}
In this diagram, the entry `associative algebras' refers to the state
sum of Fukuma--Hosono--Kawai whereas `monoidal categories' refers to
the Turaev--Viro state sum. For $n=2$, it should actually read
`\emph{strongly separable} associative algebras'. The appropriate
choice of adjectives for the other cases is in fact not systematically
understood. In order to settle this question and in order to extend
the diagram upwards to higher dimension, one can ask whether it is
possible to classify the algebraic structures from which one can
construct $n$-dimensional state sum TQFTs for generic $n$.

Whereas the algebraic structures of~\eqref{eq_ladder} that are
relevant to the state sum construction, are closely related to
Pachner moves~\cite{Pachner} and to the coherence conditions in
higher categories, they have no obvious relationship to the global
description of the TQFT as a functor
$Z\colon\cat{nCob}\to\cat{Vect}_k$.

By showing that the associative algebra $A$ of the
Fukuma--Hosono--Kawai state sum is precisely the vector space $A=Z(I)$
associated with the unit interval in an appropriately extended notion
of $2$-dimensional TQFT, we have revealed such a relationship for the
simplest case $n=2$ of the dimensional ladder~\eqref{eq_ladder}. This
raises the question of whether one can find topological
interpretations for the other algebraic structures featured
in~\eqref{eq_ladder}, presumably by extending the notion of TQFT from
conventional cobordisms to manifolds with corners of higher and higher
codimension. Further evidence for such a relationship is provided by
the Hopf algebra object in $3$-dimensional extended
TQFTs~\cite{CY2,Yetter,Kerler} in connection with Kuperberg's
$3$-manifold invariant which is based on certain Hopf
algebras~\cite{Kuperberg}.

In the present article, we consider Frobenius algebras not only in the
category $\cat{Vect}_k$ of vector spaces, but in any symmetric
monoidal Abelian category $\cal{C}$. This extends our results without
any additional work to Frobenius algebras in the category of graded
vector spaces or in the category of chain complexes, \etc.

The groupoid algebra $k[\cal{G}]$ of a finite groupoid $\cal{G}$
forms an example of a strongly separable algebra for suitable
fields $k$. We show that our generalized state sum for this
algebra yields an easy example of an open-closed TQFT with
D-branes.

The present article is structured as follows. In
Section~\ref{sect_frob}, we collect the key definitions and facts
about symmetric Frobenius algebras, strongly separable symmetric
Frobenius algebras, and knowledgeable Frobenius algebras. We also
introduce convenient diagrams. In Section~\ref{sect_topology}, we
recall the definition of the category $\cat{2Cob}^{\mathrm{ext}}$
of open-closed cobordisms and how to triangulate these. The state
sum construction of combinatorial open-closed TQFTs is then
presented in Section~\ref{sect_statesum}.

%
\section{Frobenius Algebras}
%
\label{sect_frob}

\subsection{Symmetric monoidal categories and string diagrams}

In this section, we review the basics of string diagrams in a
symmetric monoidal category. The symmetric monoidal categories of
interest will be Abelian symmetric monoidal categories so that the
$\Hom$ spaces are Abelian groups and the notion of kernels and
cokernels are defined.  Such categories include the categories of
vector spaces, graded vector spaces, $R$-modules for a commutative
ring $R$, and chain complexes of each of these structures. We denote a
symmetric monoidal category $\cal{C}$ as
$(\cal{C},\otimes,\1,\alpha,\lambda,\rho,\tau)$ where $\cal{C}$ is a
category and $\otimes$ provides $\cal{C}$ with a monoidal structure
with unit object $\1$ whose associator is denoted $\alpha$ and whose
left and right unit constraints are given by $\lambda$ and $\rho$. The
symmetric braiding is denoted $\tau$.

We denote the class of objects of a category $\cal{C}$ by $|\cal{C}|$
and for each object $X\in|\cal{C}|$, the identity morphism by
$\id_X\colon X\to X$.

\begin{defn}
Let $(\cal{C},\otimes,\1,\alpha,\lambda,\rho,\tau)$ be a symmetric
monoidal category.
\begin{enumerate}
\item
  An object $X$ of $\cal{C}$ is called \emph{rigid} if it has a
  \emph{left-dual} $(X^\ast,\ev_X,\coev_X)$. This is an object
  $X^\ast$ of $\cal{C}$ with morphisms $\ev_X\colon X^\ast\otimes
  X\to\1$ (\emph{evaluation}) and $\coev_X\colon\1\to X\otimes X^\ast$
  (\emph{coevaluation}) which satisfy the \emph{zig-zag identities},
\begin{eqnarray}
\label{eq_dualzigzag}
\begin{aligned}
 \rho_X \circ (\id_X\otimes\ev_X) \circ
\alpha_{X,X^\ast,X}\circ (\coev_X\otimes\id_X) \circ\lambda_X^{-1}
 &=& \id_X \\
\lambda_{X^\ast}\circ (\ev_X\otimes\id_{X^\ast})\circ
\alpha^{-1}_{X^\ast,X,X^\ast}\circ
(\id_{X^\ast}\otimes\coev_X)\circ \rho_{X^\ast}^{-1}&=&
\id_{X^\ast}
\end{aligned}
\end{eqnarray}
\item
  Let $X$ be a rigid object of $\cal{C}$ and $f\in\Hom(X,X)$. The
  \emph{categorical trace} $\tr_X(f)$ is defined by,
\begin{equation}
\label{eq_trace}
  \tr_X(f):=\ev_X\circ\tau_{X,X^\ast}\circ(f\otimes\id_{X^\ast})\circ\coev_X
  \in\Hom(\1,\1).
\end{equation}
\item
  The \emph{categorical dimension} $\dim X$ of a rigid object $X$ of
  $\cal{C}$ is defined by,
\begin{equation}
\label{eq_dim}
  \dim X := \tr_X(\id_X)\in\Hom(\1,\1).
\end{equation}
\item
  For rigid objects $X$ and $Y$ of $\cal{C}$ and $f\in\Hom(X,Y)$, the
  morphism,
\begin{eqnarray}
\label{eq_dual}
  f^\ast&:=&\lambda_{X^\ast}\circ(\ev_Y\otimes\id_{X^\ast})
    \circ((\id_{Y^\ast}\otimes f)\otimes\id_{X^\ast})\nn\\
  &&\circ\alpha^{-1}_{Y^\ast,X,X^\ast}\circ(\id_{Y^\ast}\otimes\coev_X)
    \circ\rho^{-1}_{Y^\ast}\colon Y^\ast\to X^\ast,
\end{eqnarray}
  is called the \emph{dual} of $f$.
\end{enumerate}
\end{defn}

In the following, we use \emph{string diagrams}~\cite{js2,Street5} to
visualize morphisms of a given symmetric monoidal category $\cal{C}$
and the identities between them. The diagrams are read from top to
bottom. For each object $X\in|\cal{C}|$, the identity morphism $\id_X$
is denoted by a line labeled `$X$' with an arrow pointing down. The
identity morphism $\id_{X^\ast}$ of the dual object has the arrow
pointing up. For a morphism $f\colon X\to Y$, we write a disc labeled
`$f$', called a \emph{coupon}. This disc has a white side which always
faces the reader and a black side which never does so,
\begin{equation}
  \id_X = \;\; \xy 0;/r.19pc/:
    (0,12)*\stringiddl{\Ltwo{X}{}}
  \endxy,\qquad
  \id_{X^\ast} = \;\; \xy 0;/r.19pc/:
    (0,12)*\stringidul{\Ltwo{X}{}}
  \endxy,\qquad
  f= \; \xy 0;/r.19pc/:
    (0,12)*\stringmapd{\Lthree{X}{Y}{f}}
  \endxy    \; .
\end{equation}
Composition of morphisms is depicted by vertically concatenating
the corresponding diagrams; for example, for morphisms $f\colon
X\to Y$ and $g\colon Y\to Z$,
\begin{equation}
  g\circ f= \;  \xy
    (0,12)*\stringmapds{\Lthree{X}{}{f}};
    (3,0)*{Y};
    (0,0)*\stringmapds{\Lthree{}{Z}{g}}
  \endxy \; = \; \xy
    (0,12)*\stringmapdb{\Lthree{X}{Z}{g\circ f}}
  \endxy    \;  .
\end{equation}
The tensor product of morphisms is visualized by putting diagrams
next to each other; for example, for $f\colon X_1\to Y_1$ and
$g\colon X_2\to Y_2$,
\begin{equation}
  \id_{X_1\otimes X_2} = \;\;
  \xy 0;/r.19pc/:
    (0,12)*\stringiddl{\Lone{X_1\otimes X_2}};
  \endxy =  \;\;
  \xy 0;/r.19pc/:
    (0,12)*\stringiddl{\Lone{X_1}};
    (8,12)*\stringiddl{\Lone{X_2}}
  \endxy,\qquad f\otimes g=\xy
    (0,12)*\stringmapdb{\Lthree{X_1\otimes X_2}{Y_1\otimes Y_2}{f\otimes g}};
    (12,12)*{}
  \endxy =  \;\;
  \xy 0;/r.19pc/:
    (0,12)*\stringmapd{\Lthree{X_1}{Y_1}{f}};
    (8,12)*\stringmapd{\Lthree{X_2}{Y_2}{g}}
  \endxy    \; .
\end{equation}
The symmetric braiding is denoted by,
\begin{equation}
  \tau_{X,Y}=\xy 0;/r.19pc/:
    (0,6)*\stringbraidd{\Ltwo{X}{Y}}
  \endxy.
\end{equation}
Mac Lane's coherence theorem for monoidal categories~\cite{mac1}
then ensures that one can unambiguously translate any such string
diagram into a morphism of $\cal{C}$. One therefore chooses
parentheses for all tensor products that occur in the source and
target objects of the morphism and inserts the structure
isomorphisms $\alpha$, $\lambda$, $\rho$ appropriately. The
coherence theorem implies that all possible ways of inserting the
structure isomorphisms yield equal morphisms, \ie\ that there is a
well-defined morphism of $\cal{C}$ specified by the diagram. In
addition, the coherence theorem allows us to suppress the lines
associated with the unit object $\1$.

For a rigid object $X\in\cal{C}$, evaluation and coevaluation are
represented by these diagrams:
\begin{equation}
  \ev_X=\xy 0;/r.19pc/:
    (0,6)*\stringcup{\Lone{X}}
  \endxy,\qquad
  \coev_X=\xy 0;/r.19pc/:
    (0,-6)*\stringcap{\Lone{X}}
  \endxy.
\end{equation}
The zig-zag identities of~\eqref{eq_dualzigzag} are represented in
string diagrams as:
\begin{equation}
  \xy 0;/r.19pc/:
    (0,0)*\stringcap{};
    (12,0)*\stringcup{};
    (-6,0)*\stringidd{};
    (18,12)*\stringidd{\Lone{X}};
  \endxy = \xy 0;/r.19pc/:
    (0,12)*\stringidd{\Lone{X}};
    (0,0)*\stringidd{};
    (-2,0)*{};
  \endxy,\qquad\qquad\xy 0;/r.19pc/:
    (0,0)*\stringcup{};
    (12,0)*\stringcap{};
    (-6,12)*\stringidu{\Lone{X}};
    (18,0)*\stringidu{};
  \endxy = \xy 0;/r.19pc/:
    (0,12)*\stringidu{\Lone{X}};
    (0,0)*\stringidu{};
    (-2,0)*{};
  \endxy,
\end{equation}
and the definitions of trace~\eqref{eq_trace},
dimension~\eqref{eq_dim} and dual morphism~\eqref{eq_dual} are:
\begin{equation}
  \tr_X(f):=\xy 0;/r.19pc/:
    (0,12)*\stringcap{\Lone{X}};
    (-6,12)*\stringmapds{\Lthree{}{}{f}};
    (6,12)*\stringidu{\Ltwo{}{}};
    (0,0)*\stringbraiddu{};
    (0,-12)*\stringcup{}
  \endxy,\qquad
  \dim X:=\xy 0;/r.19pc/:
    (0,6)*\stringcap{\Lone{X}};
    (0,6)*\stringbraiddu{};
    (0,-6)*\stringcup{}
  \endxy,\qquad
    \xy 0;/r.19pc/:
    (0,18)*\stringmapdm{\Lthree{Y^\ast}{X^\ast}{f^\ast}}
    \endxy
    :=
    \; \xy 0;/r.19pc/:
    (-12,18)*\stringidul{\Ltwo{Y}{}};
    (6,6)*\stringcap{};
    (0,6)*\stringmapds{\Lthree{}{}{f}};
    (-6,-6)*\stringcup{};
    (12,6)*\stringidul{\Ltwo{}{X}}
  \endxy .
\end{equation}

If $\cal{C}$ is locally small, the set $\Hom(\1,\1)$ forms a
commutative monoid with multiplication
$\xi_1\cdot\xi_2:=\lambda_\1\circ(\xi_1\otimes\xi_2)\circ\lambda_\1^{-1}$
for $\xi_1,\xi_2\in\Hom(\1,\1)$ and unit $\id_\1$. The monoid
$\Hom(\1,\1)$ acts on $\Hom(X,Y)$ for all $X,Y\in|\cal{C}|$ by
$\xi\cdot f:=\lambda_Y\circ(\xi\otimes f)\circ\lambda_X^{-1}$ where
$f\in\Hom(X,Y)$ and $\xi\in\Hom(\1,\1)$.

The coherence theorem now allows us to view the elements of
$\Hom(\1,\1)$ as scalars by which the entire diagram is
multiplied.

\subsection{Frobenius algebras}

We consider Frobenius algebras not only in the symmetric monoidal
category $\cat{Vect}_k$ of vector spaces over some fixed field $k$,
but in any generic symmetric monoidal category $\cal{C}$. The
following definitions and results specialize to the usual notions in
the case of $\cal{C}=\cat{Vect}_k$.

\begin{defn}
Let $(\cal{C},\otimes,\1,\alpha,\lambda,\rho,\tau)$ be a symmetric
monoidal category.
\begin{enumerate}
\item
  An \emph{algebra object} $(A,\mu,\eta)$ in $\cal{C}$ consists of an
  object $A$ and morphisms $\mu\colon A\otimes A\to A$ and
  $\eta\colon\1\to A$ of $\cal{C}$ such that the following equalities
  are satisfied:
\begin{equation}
\mu \circ (\id_A\otimes\mu) = \mu\circ
(\mu\otimes\id_A)\circ\alpha_{A,A,A}
\end{equation}
and
\begin{eqnarray}
\label{eq_algunit}
\mu \circ (\eta\otimes\id_A) &=& \lambda_A, \\
\mu \circ (\id_A\otimes\eta) &=& \rho_A .
\end{eqnarray}
 \item
  A \emph{coalgebra object} $(A,\Delta,\epsilon)$ in $\cal{C}$ consists
  of an object $A$ and morphisms $\Delta\colon A\to A\otimes A$ and
  $\epsilon\colon A\to\1$ of $\cal{C}$ such that:
\begin{equation}
  (\id_A\otimes\Delta)\circ\Delta
  =\alpha_{A,A,A}\circ(\Delta\otimes\id_A)\circ\Delta
\end{equation}
and
\begin{eqnarray}
 (\epsilon\otimes\id_A) \circ \Delta &=& \lambda_A^{-1}, \\
 (\id_A\otimes\epsilon)\circ \Delta &=& \rho_A^{-1}.
\end{eqnarray}
 \item
  A \emph{Frobenius algebra object} $(A,\mu,\eta,\Delta,\epsilon)$
  in $\cal{C}$ consists of an object $A$ and of morphisms $\mu$,
  $\eta$, $\Delta$, $\epsilon$ of $\cal{C}$ such that:
\begin{enumerate}
\item
  $(A,\mu,\eta)$ is an algebra object in $\cal{C}$,
\item
  $(A,\Delta,\epsilon)$ is a coalgebra object in $\cal{C}$, and
\item
  the following compatibility condition, called the \emph{Frobenius
  relation}, holds,
\begin{equation}
(\id_A\otimes\mu)\circ \alpha_{A,A,A} \circ (\Delta\otimes\id_A) =
\Delta\circ \mu = (\mu\otimes\id_A)\circ \alpha^{-1}_{A,A,A} \circ
(\id_A\otimes\Delta).
\end{equation}
\end{enumerate}
\item
  A Frobenius algebra object $(A,\mu,\eta,\Delta,\epsilon)$ in
  $\cal{C}$ is called \emph{symmetric} if
\begin{equation}
  \epsilon\circ\mu = \epsilon\circ\mu\circ\tau.
\end{equation}
  It is called \emph{commutative} if
\begin{equation}
  \mu = \mu\circ\tau.
\end{equation}
\item
  Let $\cal{C}$ be locally small. A Frobenius algebra object
  $(A,\mu,\eta,\Delta,\epsilon)$ in $\cal{C}$ is called
  \emph{special} (also see~\cite{FS}) if
\begin{equation}
\label{eq_special}
  \epsilon\circ\eta = \xi_\1 \cdot \id_{\1}
  \qquad\mbox{and}\qquad
  \mu\circ\Delta = \xi_A \cdot\id_A.
\end{equation}
for some $\xi_\1,\xi_A\in\Hom(\1,\1)$ that are invertible in the
monoid $\Hom(\1,\1)$.
\end{enumerate}
\end{defn}

The string diagrams for the operations of a Frobenius algebra
$(A,\mu,\eta,\Delta,\epsilon)$ are as follows:
\begin{equation}
  \xy
    (0,12)*\stringmu{\Ltwo{A}{\mu}}
  \endxy,\qquad\xy
    (0,12)*\stringeta{\Ltwo{A}{\eta}}
  \endxy,\qquad\xy
    (0,12)*\stringdelta{\Ltwo{A}{\Delta}}
  \endxy,\qquad\xy
    (0,12)*\stringepsilon{\Ltwo{A}{\epsilon}}
  \endxy .
\end{equation}
In order to keep the diagrams small, from now on we replace the
coupons by vertices and also drop the label '$A$' wherever it is clear
from the context:
\begin{equation}
  \mu=\xy
    (0,4)*\bbmu{}
  \endxy,\qquad\eta=\xy
    (0,4)*\bbeta{}
  \endxy,\qquad\Delta=\xy
    (0,4)*\bbdelta{}
  \endxy,\qquad\epsilon=\xy
    (0,4)*\bbepsilon{}
  \endxy
\end{equation}
It is understood that the vertices have to be replaced by discs in the
paper plane with their white side facing the reader. Furthermore, we
drop all labels $\mu$, $\eta$, $\Delta$ and $\epsilon$ where these are
evident from the context. For example, we distinguish the operation
$\Delta$ from $\mu$ by the arrows of the lines.

The axioms of an algebra and those of a coalgebra then read:
\begin{equation} \label{eq_bbalgebra}
  \xy
    (-2,8)*\bbmu{};
    (2,8)*\bbdl{};
    (0,0)*\bbmu{}
  \endxy = \xy
    (-2,8)*\bbdr{};
    (2,8)*\bbmu{};
    (0,0)*\bbmu{}
  \endxy,\qquad \xy
    (-2,8)*\bbeta{};
    (2,8)*\bbid{};
    (0,0)*\bbmu{}
  \endxy = \xy
    (0,8)*\bbid{};
    (0,0)*\bbid{}
  \endxy = \xy
    (-2,8)*\bbid{};
    (2,8)*\bbeta{};
    (0,0)*\bbmu{}
  \endxy \; ,\qquad\xy
    (0,8)*\bbdelta{};
    (-2,0)*\bbdelta{};
    (4,0)*\bbdr{}
  \endxy = \xy
    (0,8)*\bbdelta{};
    (-4,0)*\bbdl{};
    (2,0)*\bbdelta{}
  \endxy,\qquad\xy
    (0,8)*\bbdelta{};
    (-2,0)*\bbepsilon{};
    (2,0)*\bbid{}
  \endxy = \xy
    (0,8)*\bbid{};
    (0,0)*\bbid{}
  \endxy = \xy
    (0,8)*\bbdelta{};
    (-2,0)*\bbid{};
    (2,0)*\bbepsilon{}
  \endxy \; ,
\end{equation}
and the Frobenius relation, commutativity and symmetry are
depicted as follows:
\begin{equation}
  \xy
    (-2,8)*\bbdelta{};
    (4,8)*\bbid{};
    (-4,0)*\bbid{};
    (2,0)*\bbmu{}
  \endxy = \xy
    (0,8)*\bbmu{};
    (0,0)*\bbdelta{}
  \endxy = \xy
    (-4,8)*\bbid{};
    (2,8)*\bbdelta{};
    (-2,0)*\bbmu{};
    (4,0)*\bbid{}
  \endxy,\qquad\xy
    (0,8)*\bbbraid{};
    (0,0)*\bbmu{}
  \endxy = \xy
    (-2,8)*\bbid{};
    (2,8)*\bbid{};
    (0,0)*\bbmu{}
  \endxy,\qquad\xy
    (0,10)*\bbbraid{};
    (0,2)*\bbmu{};
    (0,-6)*\bbepsilon{}
  \endxy = \xy
    (-2,10)*\bbid{};
    (2,10)*\bbid{};
    (0,2)*\bbmu{};
    (0,-6)*\bbepsilon{}
  \endxy,
\end{equation}
The conditions for a the Frobenius algebra to be special are
these:
\begin{equation}
  \xy
    (0,8)*\bbeta{};
    (0,0)*\bbepsilon{}
  \endxy = \;\; \xi_\1 \;
  \qquad\mbox{and}\qquad\xy
    (0,8)*\bbdelta{};
    (0,0)*\bbmu{}
  \endxy = \xi_A \;
  \xy
    (0,8)*\bbid{};
    (0,0)*\bbid{}
  \endxy
\end{equation}

\subsection{Symmetric Frobenius algebras}

In this section, we introduce the notion of a non-degenerate
symmetric invariant pairing in order to characterize symmetric
Frobenius algebras. In the subsequent sections, we use it to
define strongly separable algebras and to classify all symmetric
Frobenius algebra structures of a strongly separable algebra.

\begin{defn}
Let $(\cal{C},\otimes,\1,\alpha,\lambda,\rho,\tau)$ be a symmetric
monoidal category and $(A,\mu,\eta)$ be an algebra object in
$\cal{C}$.
\begin{enumerate}
\item
  A \emph{pairing} on $A$ is a morphism $g\colon A\otimes A\to\1$
  of $\cal{C}$.
\item
  A pairing $g\colon A\otimes A\to\1$ is called \emph{non-degenerate}
  if there exists a morphism $g^\ast\colon\1\to A\otimes A$ of
  $\cal{C}$ (called \emph{the inverse} of $g$) such that the
  \emph{zig-zag identities} hold,
\begin{eqnarray}
\begin{aligned}
  \rho_A \circ (\id_A\otimes g) \circ \alpha_{A,A,A} \circ
  (g^\ast\otimes\id_A) \circ \lambda_A^{-1} = \id_A,\\
  \lambda_A \circ (g\otimes\id_A) \circ \alpha^{-1}_{A,A,A} \circ
(\id_A\otimes g^\ast) \circ \rho_A^{-1} = \id_A.
\end{aligned}
\end{eqnarray}
\item
  A pairing $g\colon A\otimes A\to\1$ is called \emph{symmetric}
  if $g=g\circ\tau_{A,A}$.
\item
  A pairing $g\colon A\otimes A\to\1$ is called
  \emph{invariant}\footnote{Some authors use the term
  \emph{associative} rather than \emph{invariant}, see, for
  example~\cite{Kock}.} if,
\begin{equation}
  g\circ (\id_A\otimes\mu) \circ \alpha_{A,A,A} = g\circ
 (\mu\otimes\id_A).
\end{equation}
\end{enumerate}
\end{defn}

The string diagrams for a pairing $g\colon A\otimes A\to\1$ on an
algebra object $(A,\mu,\eta)$ in some symmetric monoidal category
$\cal{C}$ are as follows:
\begin{equation}
  \xy
    (0,12)*\stringpair{\Ltwo{A}{g}}
  \endxy,\qquad\xy
    (0,12)*\stringdual{\Ltwo{A}{g^\ast}}
  \endxy.
\end{equation}
Our shorthand notation using blackboard framing then reads:
\begin{equation}
\label{eq_pairdiag}
  g=\xy
    (0,4)*\bbpair{}
  \endxy,\qquad
  g^\ast=\xy
    (0,4)*\bbdual{}
  \endxy.
\end{equation}
The conditions of non-degeneracy, symmetry and invariance are
depicted as follows:
\begin{equation}
  \xy
    (-2,8)*\bbdual{};
    (4,8)*\bbid{};
    (-4,0)*\bbid{};
    (2,0)*\bbpair{}
  \endxy = \xy
    (0,8)*\bbid{};
    (0,0)*\bbid{}
  \endxy = \xy
    (-4,8)*\bbid{};
    (2,8)*\bbdual{};
    (4,0)*\bbid{};
    (-2,0)*\bbpair{}
  \endxy \; ,\qquad\xy
    (0,8)*\bbbraid{};
    (0,0)*\bbpair{}
  \endxy = \xy
    (-2,8)*\bbid{};
    (2,8)*\bbid{};
    (0,0)*\bbpair{}
  \endxy \; ,\qquad\xy
    (-2,8)*\bbmu{};
    (2,8)*\bbdl{};
    (0,0)*\bbpair{}
  \endxy = \xy
    (-2,8)*\bbdr{};
    (2,8)*\bbmu{};
    (0,0)*\bbpair{}
  \endxy.
\end{equation}
We also use the following shorthand notation for the `trilinear form'
$g^{(3)}\colon(A\otimes A)\otimes A\to\1$ which is defined by:
\begin{equation}
\label{eq_trilinear}
  g^{(3)}=\xy
    (0,4)*\bbtriple{}
  \endxy := \xy
    (-2,8)*\bbmu{};
    (2,8)*\bbid{};
    (0,0)*\bbpair{}
  \endxy \;,
\end{equation}
and which has the following cyclic symmetry:
\begin{equation}
  \xy
    (0,8)*\bbbraidr{};
    (0,0)*\bbtriple{}
  \endxy = \xy
    (-4,8)*\bbid{};
    (0,8)*\bbid{};
    (4,8)*\bbid{};
    (0,0)*\bbtriple{}
  \endxy = \xy
    (0,8)*\bbbraidl{};
    (0,0)*\bbtriple{}
  \endxy.
\end{equation}

\begin{lem}
\label{lem_pairing}
Let $(\cal{C},\otimes,\1,\alpha,\lambda,\rho,\tau)$ be a symmetric
monoidal category. Every symmetric Frobenius algebra object
$(A,\mu,\eta,\Delta,\epsilon)$ in $\cal{C}$ gives rise to a
non-degenerate symmetric invariant pairing $g:=\epsilon\circ\mu$ on
$A$ with inverse $g^\ast:=\Delta\circ\eta$. Conversely, given an
algebra object $(A,\mu,\eta)$ in $\cal{C}$ and a non-degenerate
symmetric invariant pairing $g$ on $A$, there is a symmetric Frobenius
algebra object $(A,\mu,\eta,\Delta,\eta)$ with
$\Delta:=(\mu\otimes\id_A)\circ\alpha_{A,A,A}^{-1}\circ(\id_A\otimes
g^\ast)\circ\rho_A^{-1}$ and
$\epsilon:=g\circ(\id_A\otimes\eta)\circ\rho_A^{-1}$.
\end{lem}

The defining equations used in this lemma can be read diagrammatically
as:
\begin{equation}
\label{eq_semi_frob}
  \xy
    (0,4)*\bbpair{}
  \endxy := \xy
    (0,8)*\bbmu{};
    (0,0)*\bbepsilon{}
  \endxy,\qquad\xy
    (0,4)*\bbdual{}
  \endxy := \xy
    (0,8)*\bbeta{};
    (0,0)*\bbdelta{}
  \endxy\; ,\qquad\mbox{and}\qquad\xy
    (0,4)*\bbdelta{}
  \endxy := \xy
    (-4,8)*\bbid{};
    (2,8)*\bbdual{};
    (-2,0)*\bbmu{};
    (4,0)*\bbid{}
  \endxy\; ,\qquad\xy
    (0,4)*\bbepsilon{}
  \endxy := \xy
    (-2,8)*\bbid{};
    (2,8)*\bbeta{};
    (0,0)*\bbpair{}
  \endxy.
\end{equation}
Notice that every symmetric Frobenius algebra object
$(A,\mu,\eta,\Delta,\epsilon)$ in $\cal{C}$ is a rigid object of
$\cal{C}$ with left-dual\footnote{It is right-dual at the same time,
but we do not refer to this property in the following.}
$(A,\epsilon\circ\mu,\Delta\circ\eta)$.

\subsection{Strongly separable algebras}

Every rigid algebra object in a symmetric monoidal category is
equipped with a canonical pairing. Recall first the special case
of $\cal{C}=\cat{Vect}_k$ for an arbitrary field $k$. Let $A$ be a
finite-dimensional algebra over $k$ and denote the left-regular
representation by $L\colon A\to\End_k(A),a\mapsto L_a$ with
$L_a\colon A\to A,b\mapsto ab$. By associativity, $L_{ab}=L_a\circ
L_b$ for all $a,b\in A$. The trace of the matrices of the
left-regular representation equips $A$ with a canonical bilinear
form,
\begin{equation}
  g_\mathrm{can}\colon A\otimes A\to k,\quad
  a\otimes b\mapsto \tr_A(L_{ab}),
\end{equation}
which can be shown to be symmetric and invariant. We are interested in
those algebras for which this canonical bilinear form is
non-degenerate. These are the strongly separable algebras. Let us
first recall the definition.

\begin{defn}
Let $A$ be an algebra over a commutative ring $r$. We denote by
$A^\mathrm{op}$ the opposite algebra of $A$, by $A^e=A\otimes
A^\mathrm{op}$ its enveloping algebra and by $\mu\colon A^e\to A$,
$a\otimes b\mapsto ab$ the augmentation mapping. A is called
\emph{separable} if there is an element $e\in A^e$ (called a
\emph{separability idempotent}) such that,
\begin{enumerate}
\item
  $(a\otimes 1)e = (1\otimes a)e$ holds in $A^e$ for all $a\in A$.
\item
  $\mu(e)=1$.
\end{enumerate}
$A$ is called \emph{strongly separable} if $A$ is separable with a
separability idempotent that satisfies $\tau_{A,A}(e)=e$.
\end{defn}

\begin{thm}[see, for example~\cite{Aguiar}]
\label{thm_stronglysep}
Let $A$ be an algebra over any field $k$. Then the following are
equivalent:
\begin{enumerate}
\item
  $A$ is finite-dimensional over $k$, and the canonical bilinear form
  is non-degenerate.
\item
  $A$ is strongly separable.
\end{enumerate}
\end{thm}

Every strongly separable algebra therefore carries a canonical
symmetric Frobenius algebra structure by Lemma~\ref{lem_pairing}. The
following definition of a canonical pairing for generic $\cal{C}$
reduces to the canonical bilinear form in the case of
$\cal{C}=\cat{Vect}_k$.

\begin{prop}
Let $(\cal{C},\otimes,\1,\alpha,\lambda,\rho,\tau)$ be a symmetric
monoidal category and $(A,\mu,\eta)$ be an algebra object in $\cal{C}$
such that the object $A$ is rigid with left-dual
$(A^\ast,\ev_A,\coev_A)$. Then there is a symmetric invariant pairing
on $A$ given by,
\begin{equation}
\label{eq_pairing}
  g_\mathrm{can}:=
    \ev_A\circ\tau_{A,A^\ast}\circ(\mu\otimes\id_{A^\ast})\circ\alpha_{A,A,A^\ast}^{-1}
      \circ(\mu\otimes\coev_A)\circ\rho_{A\otimes A}^{-1} =
  \xy
    (-2,16)*\bbmu{};
    (4,16)*\bbcap{};
    (0,8)*\bbmu{};
    (4,8)*\bbdlu{};
    (2,0)*\bbbraiddu{};
    (2,-8)*\bbcup{}
  \endxy.
\end{equation}
which we call the \emph{canonical pairing}.
\end{prop}

\begin{defn}
A rigid algebra object in a symmetric monoidal category is called
\emph{strongly separable} if its canonical pairing is non-degenerate.
\end{defn}

By Theorem~\ref{thm_stronglysep}, this notion of a strongly separable
algebra object in some symmetric monoidal category agrees with the
usual definition in the case $\cal{C}=\cat{Vect}_k$. We are not aware
of any such result for the more general case of modules over a
commutative ring. In order to illustrate how strong the condition of
strong separability is, we include the following results and examples
from~\cite{Aguiar,Pierce}.

\begin{thm}
Let $A$ be an algebra over some field $k$.
\begin{enumerate}
\item
  If $A$ is strongly separable, then $A$ is finite-dimensional,
  separable, and semisimple.
\item
  If $A$ is separable and commutative, then $A$ is strongly separable.
\item
  If $A$ is finite-dimensional and semisimple and $\chr k=0$, then
  $A$ is strongly separable.
\item
  If $A$ is finite-dimensional and semisimple and $k$ is a perfect
  field, then $A$ is separable.
\end{enumerate}
\end{thm}

\begin{myexample}
\label{ex_groupalg}
Let $k$ be a field and $G$ be a finite group.
\begin{enumerate}
\item
  If $\chr k$ does not divide the order of $G$, then the group
  algebra $k[G]$ is strongly separable.
\item
  If $\chr k$ divides the order of $G$, then $k[G]$ is neither
  semisimple nor separable.
\end{enumerate}
\end{myexample}

\begin{myexample}
\label{ex_matrixalg}
Let $k$ be a field and $M_n(k)$ be the algebra of $(n\times
n)$-matrices over $k$.
\begin{enumerate}
\item
  If $\chr k$ does not divide $n$, then $M_n(k)$ is
  strongly separable.
\item
  If $\chr k$ divides $n$, then $M_n(k)$ is semisimple and separable,
  but not strongly separable.
\end{enumerate}
\end{myexample}

In both examples, the non-degeneracy of the canonical bilinear form is
a convenient criterion for strong separability. We explain below why
in the Examples~\ref{ex_groupalg}(2) and~\ref{ex_matrixalg}(2), the
state sum construction fails. In particular, for a finite field of
non-zero characteristic $p$, the original Fukuma--Hosono--Kawai state
sum~\cite{FHK} cannot be applied to the $(p\times p)$-matrix algebra
$A:=M_p(k)$ although $k$ is perfect and $M_p(k)$ is
finite-dimensional, separable, and semisimple.

\subsection{Strongly separable symmetric Frobenius algebras}

In this section, we compare the pairing $\epsilon\circ\mu$ of a
generic symmetric Frobenius algebra with the canonical pairing. They
differ by multiplication with a central element which we call the
\emph{window element}\footnote{This terminology is inspired by the
open-closed cobordism~\cite{LP05} that is associated with this
element.}.

In a generic locally small symmetric monoidal category
$(\cal{C},\otimes,\1,\alpha,\lambda,\rho,\tau)$, we use the
terminology \emph{element of} $A$ for a morphism $a\colon\1\to A$. The
set $\Hom(\1,A)$ of elements of an algebra object $(A,\mu,\eta)$ in
$\cal{C}$ forms a monoid with respect to convolution $a\cdot
b:=\mu\circ(a\otimes b)\circ\lambda_\1^{-1}$ for elements
$a,b\in\Hom(\1,A)$ and with unit $\eta$. An element $a\in\Hom(\1,A)$
is called \emph{central} if it is contained in the commutative
submonoid,
\begin{equation}
\label{eq_centremor}
  \cal{Z}(A):=\{\,a\in\Hom(\1,A)\colon\quad
    \mu\circ(a\otimes\id_A)\circ\lambda_A^{-1}
     = \mu\circ(\id_A\otimes a)\circ\rho_A^{-1}\,\}.
\end{equation}
The set of invertible elements of $A$ forms a group
${\Hom(\1,A)}^\times\subseteq\Hom(\1,A)$, and the set of
invertible central elements
${\cal{Z}(A)}^\times:=\cal{Z}(A)\cap{\Hom(\1,A)}^\times\leq{\Hom(\1,A)}^\times$
a subgroup. This means in particular that the inverse of every
central element is central, too. $\cal{Z}(A)$ acts on $\Hom(A,A)$
by
\begin{equation}
  \cal{Z}(A)\times\Hom(A,A)\to\Hom(A,A),\quad (a,f)\mapsto a\cdot
  f:=\mu\circ(a\otimes f)\circ\lambda_A^{-1}.
\end{equation}
We also have $(a\cdot\id_A)\circ\eta=a$ and
$(a\cdot\id_A)\circ(b\cdot\id_A)=(a\cdot b)\cdot\id_A$ for all
$a,b\in\cal{Z}(A)$ as well as
\begin{equation}
  \mu\circ((a\cdot\id_A)\otimes\id_A) = (a\cdot\id_A)\circ\mu =
  \mu\circ(\id_A\otimes(a\cdot\id_A)),
\end{equation}
and for a Frobenius algebra object $(A,\mu,\eta,\Delta,\epsilon)$
also,
\begin{equation}
  ((a\cdot\id_A)\otimes\id_A)\circ\Delta
  = \Delta\circ(a\cdot\id_A)
  = (\id_A\otimes(a\cdot\id_A))\circ\Delta.
\end{equation}
The following diagrams show an element $a\in\Hom(\1,A)$, the morphism
$a\cdot\id_A$, and the centrality condition:
\begin{equation}
  a = \xy
    (0,4)*\bbelement{a}
  \endxy,\qquad
  a\cdot\id_A = \xy
    (-2,8)*\bbelement{a};
    (2,8)*\bbid{};
    (0,0)*\bbmu{}
  \endxy,\qquad \xy
    (-2,8)*\bbelement{a};
    (2,8)*\bbid{};
    (0,0)*\bbmu{}
  \endxy = \xy
    (-2,8)*\bbid{};
    (2,8)*\bbelement{a};
    (0,0)*\bbmu{}
  \endxy
\end{equation}

\begin{defn}
Let $(A,\mu,\eta,\Delta,\epsilon)$ be a symmetric Frobenius algebra
object in a locally small symmetric monoidal category $\cal{C}$. The
\emph{window element} of $A$ is defined by,
\begin{equation}
\label{eq_window}
  a := \mu\circ\Delta\circ\eta = \xy
    (0,12)*\bbeta{};
    (0,4)*\bbdelta{};
    (0,-4)*\bbmu{}
  \endxy
\end{equation}
\end{defn}

The window element is a central element. The comparison between the
pairing $\epsilon\circ\mu$ of a generic symmetric Frobenius algebra
and the canonical pairing can be done as follows.

\begin{prop}
\label{prop_canwindow}
Let $\cal{C}$ be a locally small symmetric monoidal category and
$(A,\mu,\eta,\Delta,\epsilon)$ be a symmetric Frobenius algebra object
in $\cal{C}$. Then the canonical pairing of $A$ is given by,
\begin{equation}
\label{eq_paircomp}
  g_\mathrm{can} = \epsilon\circ(a\cdot\id_A)\circ\mu,
\end{equation}
where $a$ denotes the window element.
\end{prop}

\begin{proof}
Notice that $A$ is a rigid object of $\cal{C}$, and so it makes sense
to study the canonical pairing~\eqref{eq_pairing}. We use the diagrams
of~\eqref{eq_pairdiag} for $g:=\epsilon\circ\mu$ and
$g^\ast:=\Delta\circ\eta$:
\begin{equation}
  g_\mathrm{can} =
  \xy
    (-2,16)*\bbmu{};
    (4,16)*\bbcap{};
    (0,8)*\bbmu{};
    (4,8)*\bbdlu{};
    (2,0)*\bbbraiddu{};
    (2,-8)*\bbcup{}
  \endxy = \xy
    (-2,16)*\bbmu{};
    (4,16)*\bbdual{};
    (0,8)*\bbmu{};
    (4,8)*\bbdl{};
    (2,0)*\bbbraid{};
    (2,-8)*\bbpair{}
  \endxy = \xy
    (-2,16)*\bbmu{};
    (4,16)*\bbdual{};
    (0,8)*\bbmu{};
    (4,8)*\bbdl{};
    (2,0)*\bbpair{}
  \endxy = \xy
    (-2,24)*\bbeta{};
    (2,16)*\bbmu{};
    (-2,16)*\bbdelta{};
    (-2,8)*\bbmu{};
    (2,8)*\bbid{};
    (0,0)*\bbmu{};
    (0,-8)*\bbepsilon{};
  \endxy = \xy
    (2,16)*\bbmu{};
    (-2,8)*\bbelement{a};
    (2,8)*\bbid{};
    (0,0)*\bbmu{};
    (0,-8)*\bbepsilon{};
  \endxy
\end{equation}
The first equality is the definition; for the second one, we have
exploited the fact that $A$ satisfies the zig-zag
identities~\eqref{eq_dualzigzag} both with $(A^\ast,\ev_A,\coev_A)$
and with $(A,g,g^\ast)$; the third equality is symmetry; the fourth
one follows from the axioms of a Frobenius algebra; and the fifth is
the definition of the window element.
\end{proof}

In the above proposition, the pairing $g:=\epsilon\circ\mu$ is always
non-degenerate whereas the canonical pairing is non-degenerate if and
only if the algebra is strongly separable.

\begin{thm}
Let $(\cal{C},\otimes,\1,\alpha,\lambda,\rho,\tau)$ be a locally
small symmetric monoidal category and
$(A,\mu,\eta,\Delta,\epsilon)$ be a symmetric Frobenius algebra
object in $\cal{C}$ with window element $a$. Then the following
are equivalent:
\begin{enumerate}
\item
  The algebra object $(A,\mu,\eta)$ is strongly separable.
\item
  The window element is invertible.
\end{enumerate}
\end{thm}

\begin{proof}
Let $(A,\mu,\eta)$ be strongly separable. Then the canonical pairing
is non-degenerate and therefore has an inverse
$g_\mathrm{can}^\ast\colon\1\to A\otimes A$. We denote the pairing
$g:=\epsilon\circ\mu$ and its inverse $g^\ast:=\Delta\circ\eta$ by the
diagrams of~\eqref{eq_pairdiag} and the canonical pairing and its
inverse by,
\begin{equation}
  g_\mathrm{can} = \xy
    (0,4)*\bbpairp{}
  \endxy,\qquad
  g_\mathrm{can}^\ast = \xy
    (0,4)*\bbdualp{}
  \endxy.
\end{equation}
Define $\tilde a:=\lambda_A\circ(\epsilon\otimes\id_A)\circ
g_\mathrm{can}^\ast$. Then $\tilde a$ is the inverse of $a$
because,
\begin{equation}
  \xy
    (-4,8)*\bbelement{a};
    (4,8)*\bbelement{\tilde a};
    (0,0)*\bbmediummu{};
  \endxy = \xy
    (0,16)*\bbdualp{};
    (-2,8)*\bbepsilon{};
    (2, 8)*\bbid{};
    (-2,0)*\bbelement{a};
    (2, 0)*\bbid{};
    (0,-8)*\bbmu{};
  \endxy = \xy
    (0,16)*\bbdualp{};
    (-2,8)*\bbepsilon{};
    (2, 8)*\bbid{};
    (-2,0)*\bbelement{a};
    (2, 0)*\bbid{};
    (0,-8)*\bbmu{};
    (6,-8)*\bbdual{};
    (2,-16)*\bbpair{};
    (8,-16)*\bbid{}
  \endxy = \xy
    (0,16)*\bbdualp{};
    (8,16)*\bbdual{};
    (-2,8)*\bbepsilon{};
    (4, 8)*\bbmu{};
    (10,8)*\bbid{};
    (0, 0)*\bbelement{a};
    (4, 0)*\bbid{};
    (10,0)*\bbid{};
    (2,-8)*\bbmu{};
    (10,-8)*\bbid{};
    (2,-16)*\bbepsilon{};
    (10,-16)*\bbid{}
  \endxy = \xy
    (0, 8)*\bbdualp{};
    (8, 8)*\bbdual{};
    (-2,0)*\bbepsilon{};
    (4, 0)*\bbpairp{};
    (10, 0)*\bbid{}
  \endxy = \xy
    (0, 8)*\bbdual{};
    (-2,0)*\bbepsilon{};
    (2, 0)*\bbid{}
  \endxy = \xy
    (0,4)*\bbeta{};
  \endxy.
\end{equation}
The first equality is the definition of $\tilde a$; the second one a
zig-zag identity for $(A,g,g^\ast)$; the third one is a consequence of
the axioms of a Frobenius algebra; the fourth one
is~\eqref{eq_paircomp}; the fifth one is a zig-zag identity for
$(A,g_\mathrm{can},g_\mathrm{can}^\ast)$; and the last equation holds
by the axioms of a Frobenius algebra.

Conversely, let the window element $a$ have an inverse $a^{-1}$. Then
a similar computation shows that the canonical pairing satisfies the
zig-zag identities with the inverse $g_\mathrm{can}^\ast = \Delta\circ
a^{-1}$ and is therefore non-degenerate.
\end{proof}

Combined with Proposition~\ref{prop_canwindow}, the preceding
theorem implies that the symmetric Frobenius algebra structures of
a given strongly separable algebra are characterized by the
invertible central elements. The proposition below describes the
extent to which the notion of a strongly separable symmetric
Frobenius algebra generalizes the notion of a special symmetric
Frobenius algebra.

\begin{prop}
Let $\cal{C}$ be a locally small symmetric monoidal category and
$(A,\mu,\eta,\Delta,\epsilon)$ be a symmetric Frobenius algebra object
in $\cal{C}$ such that $\dim A$ is invertible in $\Hom(\1,\1)$. Then
the following are equivalent:
\begin{enumerate}
\item
  The Frobenius algebra object $(A,\mu,\eta,\Delta,\epsilon)$ is
  special with $\epsilon\circ\eta=\xi_\1\cdot\id_\1$ and
  $\mu\circ\Delta=\xi_A\cdot\id_A$ for some invertible
  $\xi_\1,\xi_A\in\Hom(\1,\1)$.
\item
  The algebra object $(A,\mu,\eta)$ is strongly separable, and the
  window element is of the form $a=\zeta\cdot\eta$ for some invertible
  $\zeta\in\Hom(\1,\1)$.
\end{enumerate}
In this case, $\xi_A=\zeta$ and $\xi_\1=\xi_A^{-1}\dim A$.
\end{prop}

\begin{proof}
If $A$ is special, the window element is
$a=\mu\circ\Delta\circ\eta=(\xi_A\cdot\id_A)\circ\eta =
\xi_A\cdot\eta$. It is invertible with $a^{-1}=\xi_A^{-1}\cdot\eta$,
and so $(A,\mu,\eta)$ is strongly separable.

Conversely, if $(A,\mu,\eta)$ is strongly separable with window
element $a=\zeta\cdot\eta$ for some invertible $\zeta\in\Hom(\1,\1)$,
then the second condition of~\eqref{eq_special} holds with invertible
$\xi_A=\zeta$. For a symmetric Frobenius algebra object in a symmetric
monoidal category, the second condition of~\eqref{eq_special} implies
the first one with $\xi_\1=\xi_A^{-1}\dim A$:
\begin{equation}
\dim A =
  \xy
    (0,12)*\bbcap{};
    (0,4)*\bbbraiddu{};
    (0,-4)*\bbcup{}
  \endxy = \xy
    (0,12)*\bbdual{};
    (0,4)*\bbbraid{};
    (0,-4)*\bbpair{}
  \endxy = \xy
    (0,8)*\bbdual{};
    (0,0)*\bbpair{}
  \endxy = \xy
    (0,16)*\bbeta{};
    (0,8)*\bbdelta{};
    (0,0)*\bbmu{};
    (0,-8)*\bbepsilon{};
  \endxy = \xi_A \;\xy
    (0,8)*\bbeta{};
    (0,0)*\bbepsilon{}
  \endxy.
\end{equation}
Since $\dim A$ is invertible by assumption, so is $\xi_\1$.
\end{proof}

\begin{rem}
Given any strongly separable symmetric Frobenius algebra object
$(A,\mu,\eta,\Delta,\epsilon)$ with window element $a$ in a locally
small symmetric monoidal category $\cal{C}$, the identity
\begin{equation}
\label{eq_bubble}
  (a^{-1}\cdot\id_A)\circ\mu\circ\Delta =
  \xy
    (4,12)*\bbdelta{};
    (4,4)*\bbmu{};
    (2,-4)*\bbainv{1};
  \endxy = \xy
    (0,8)*\bbid{};
    (0,0)*\bbid{}
  \endxy = \id_A
\end{equation}
generalizes the `bubble move' of Fukuma--Hosono--Kawai from the
canonical symmetric Frobenius algebra structure to the case of a
generic symmetric Frobenius algebra structure. In
Section~\ref{sect_statesum}, we explain why this generalization is
needed in order to obtain a sharp invariant of open-closed
cobordisms from the state sum.

For the algebras of Example~\ref{ex_groupalg}(2) and
Example~\ref{ex_matrixalg}(2) which are not strongly separable,
the morphism $\mu\circ\Delta$ is zero, and so there is no way of
obtaining an analogue of the `bubble move'.
\end{rem}

\subsection{Knowledgeable Frobenius algebras}
\label{sect_kfrob}

We have shown in~\cite{LP05} that open-closed TQFTs, \ie\ symmetric
monoidal functors $Z\colon\cat{2Cob}^\mathrm{ext}\to\cal{C}$ where
$\cal{C}$ is a symmetric monoidal category, are characterized by
\emph{knowledgeable Frobenius algebras} in $\cal{C}$. Here we just recall
the definition. For more details, we refer the reader to~\cite{LP05}.

\begin{defn}
Let $(\cal{C},\otimes,\1,\alpha,\lambda,\rho,\tau)$ be a symmetric
monoidal category. A \emph{homomorphism of algebras} $f\colon A\to
A^\prime$ between two algebra objects $(A,\mu,\eta)$ and
$(A^\prime,\mu^\prime,\eta^\prime)$ in $\cal{C}$ is a morphism $f$ of
$\cal{C}$ such that:
\begin{equation}
f \circ \mu = \mu'\circ (f \otimes f) \qquad\mbox{and}\qquad f
\circ \eta = \eta'.
\end{equation}
\end{defn}

\begin{defn}
Let $(\cal{C},\otimes,\1,\alpha,\lambda,\rho,\tau)$ be a symmetric
monoidal category. A \emph{knowledgeable Frobenius algebra}
$(A,C,\imath,\imath^\ast)$ in $\cal{C}$ consists of,
\begin{itemize}
\item
  a symmetric Frobenius algebra $A=(A,\mu_A,\eta_A,\Delta_A,\epsilon_A)$,
\item
  a commutative Frobenius algebra $C=(C,\mu_C,\eta_C,\Delta_C,\epsilon_C)$,
\item
  morphisms $\imath\colon C\to A$ and $\imath^\ast\colon A\to C$ of $\cal{C}$,
\end{itemize}
such that $\imath\colon C\to A$ is a homomorphism of algebra objects
in $\cal{C}$ and,
\begin{alignat}{2}
\label{eq_kfrob1}
  \mu_A\circ(\imath\otimes\id_A)
    &= \mu_A\circ\tau_{A,A}\circ(\imath\otimes\id_A)
    &\qquad& \mbox{(knowledge),}\\
\label{eq_kfrob2}
  \epsilon_C\circ\mu_C\circ(\id_C\otimes\imath^\ast)
    &= \epsilon_A\circ\mu_A\circ(\imath\otimes\id_A)
    &\qquad& \mbox{(duality),}\\
\label{eq_kfrob3}
  \mu_A\circ\tau_{A,A}\circ\Delta_A
    &= \imath\circ\imath^{\ast}
    &\qquad& \mbox{(Cardy condition).}
\end{alignat}
\end{defn}

Sometimes the folk theorem on the characterization of open-closed
TQFTs is stated in such a way that it includes the condition
$C=Z(A)$. The following example shows that there exist knowledgeable
Frobenius algebras and thereby open-closed TQFTs in which this
condition does not hold.

\begin{myexample}
\label{ex_notcentre}
Let $k$ be a field, $\chr k\neq 2$ and $n\in\N$ such that $\chr k$
does not divide $n$. Assume that there exists some $\alpha\in k$ such
that $\alpha^2=-1/2$ (for example $k=\C$).

Let $A=M_n(k)$ be the $n\times n$-matrix algebra over $k$. Choose
a $k$-basis ${(e_{ij})}_{1\leq i,j\leq n}$ of $A$ such that the
multiplication is given by $\mu_A(e_{ij}\otimes
e_{k\ell})=\delta_{jk}\,e_{i\ell}$ and the unit by
$\eta_A(1)=\sum_{i=1}^ne_{ii}$. The algebra $A$ forms a symmetric
Frobenius algebra with
$\Delta_A(e_{ij})=\alpha^{-1}\sum_{k=1}^ne_{ik}\otimes e_{kj}$ and
$\epsilon_A(e_{ij})=\alpha\,\delta_{ij}$. We compute
$\mu_A\circ\Delta_A = n\alpha^{-1}\cdot\id_A$ and the window
element $a_A=n\alpha^{-1}\cdot\eta_A$. It is invertible with
$a_A^{-1}=n^{-1}\alpha\cdot\eta_A$, and so $A$ is strongly
separable. In fact, $A$ is special with $\xi_A=n\alpha^{-1}$ and
$\xi_\1=n\alpha$. Obviously, $Z(A)\cong k$.

Let $C=k[X]/(X^2-1)$. A $k$-basis is given by $(1,X)$. $C$ becomes
a commutative Frobenius algebra with $\Delta_C(1)=1\otimes
X+X\otimes 1$, $\Delta_C(X)=1\otimes1+X\otimes X$,
$\epsilon_C(1)=0$, and $\epsilon_C(X)=1$. We compute
$(\mu_C\circ\Delta_C)(c)=2X\,c$ for all $c\in C$, and the window
element is $a_C=2X$. It is invertible with $a_C^{-1}=X/2$, and so
$C$ is strongly separable, too, but it is not special.

If we define $\imath\colon C\to A$ by $\imath(1)=\eta_A(1)$ and
$\imath(X)=-\eta_A(1)$, and $\imath^\ast\colon A\to C$ by
$\imath^\ast(e_{ij})=\delta_{ij}\,\alpha(X-1)$, then
$(A,C,\imath,\imath^\ast)$ forms a knowledgeable Frobenius
algebra. Observe that $Z(A)$ is $1$-dimensional over $k$, but $C$
is $2$-dimensional, and so $Z(A)\not\cong C$.
\end{myexample}

\subsection{Idempotents}

In this section, we show that every strongly separable symmetric
Frobenius algebra $A$ in an Abelian symmetric monoidal category
$\cal{C}$ gives rise to a knowledgeable Frobenius algebra
$(A,C,\imath,\imath^\ast)$ in $\cal{C}$. In $\cat{Vect}_k$, $C$ is
isomorphic to the centre of $A$. In general, it arises as the
image of the following canonical idempotent.

\begin{prop}
\label{prop_idempotent}
Let $(A,\mu_A,\eta_A,\Delta_A,\epsilon_A)$ be a strongly separable
symmetric Frobenius algebra object in a locally small symmetric
monoidal category $\cal{C}$ and let $a^{-1}$ denote the inverse of the
window element of $A$. Then the morphism
\begin{equation}
\label{eq_idempotent}
  p = \xy
    (0,0)*\bbP{};
  \endxy := (a^{-1}\cdot\id_A)\circ\mu_A\circ\tau_{A,A}\circ\Delta_A =
  \;
  \xy
    (0,-4)*\bbmediummu{};
    (-4,-1)*\xycircle(2.9,2.9){-}="t1";
    (-4,-1)*{\scriptstyle a^{\scriptscriptstyle -1}};
    (4,4)*\bbcardy{};
  \endxy
\end{equation}
has the following properties,
\begin{enumerate}
\item
  $p^2=p$,
\item
  $p\circ\eta_A = \eta_A$,
\item
  $\epsilon_A\circ p = \epsilon_A$,
\item
  $p\circ\mu_A\circ(p\otimes p) = \mu_A\circ(p\otimes p) = p\circ\mu_A\circ(p\otimes\id_A)
   =p\circ\mu_A\circ(\id_A\otimes p)$,
\item
  $(p\otimes p)\circ\Delta_A\circ p = (p\otimes p)\circ\Delta_A = (p\otimes\id_A)\circ\Delta_A\circ p
   =(\id_A\otimes p)\circ\Delta_A\circ p$,
\item
  $c=p\circ c$ for all $c\in\cal{Z}(A)$,
\item
  $(c\cdot\id_A)\circ p = p\circ(c\cdot\id_A)$ for all $c\in\cal{Z}(A)$,
\item
  $\mu_A\circ(p\otimes\id_A)=\mu_A\circ\tau_{A,A}\circ(p\otimes\id_A)$.
\end{enumerate}
\end{prop}

In $\cat{Vect}_k$, condition~(1) states that $p$ is a projector;
condition~(8) says that its image is contained in the centre
$Z(A)$, and condition~(6) says that the centre $Z(A)$ is contained
in the image of $p$, and so $p$ projects onto the centre $Z(A)$.
Whereas this $Z(A)$ arises as a subspace $Z(A)=\im p\subseteq A$,
the \emph{centre} $\cal{Z}(A)$ according to~\eqref{eq_centremor}
consists of morphisms $\1\to A$. In $\cat{Vect}_k$, one can
evaluate any such morphism $a\in\cal{Z}(A)$ at the unit $1\in k$
of the field and finds that $a(1)\in Z(A)\subseteq A$.

Note that the idempotent~\eqref{eq_idempotent} is precisely
$p=\mu_A\circ\tau_{A,A}\circ\Delta_A^{(\mathrm{can})}$ where
$\Delta_A^{(\mathrm{can})}$ refers to the canonical symmetric
Frobenius algebra structure on $A$.

In the state sum, the idempotent~\eqref{eq_idempotent} appears
whenever a unit interval is closed to a circle, \ie\ it is closely
related with the generators $\imath$ and $\imath^\ast$
of~\eqref{eq_generators}. The image of an idempotent can be
defined in any Abelian category as follows.

\begin{prop}[see, for example~\cite{Borceux}]
\label{prop_splitting} Let $\cal{C}$ be an Abelian category and
$p\colon A\to A$ be an idempotent. The image factorization of $p$
yields an object $p(A)$, called the \emph{image} of $p$, which is
unique up to isomorphism, together with morphisms $\coim p\colon
A\to p(A)$ (called the \emph{co-image}) and $\im p\colon p(A)\to
A$ (called the \emph{image}) such that the following diagram
commutes:
\begin{equation} \label{eq_factorization}
\begin{aligned}
\xymatrix@!C=2.3pc{
  A\ar[r]^-{\coim p}\ar[dr]_{p}&p(A)\ar[d]^{\im p}\\
  &A
}
\end{aligned}
\end{equation}
Since $\cal{C}$ is Abelian, the idempotent $p$ is split. The
splitting is given precisely by the two morphisms of the image
factorization, and so we have $\id_{p(A)}=\coim p\circ\im p$.
Therefore, the short exact sequence
\begin{equation}
\begin{aligned}
\xymatrix@!C=2.2pc{
  0\ar[r]&N_p\ar[r]^-{\ker p}&A\ar@<.5ex>[r]^-{\coim p}
   &p(A)\ar[r]\ar@<.5ex>[l]^-{\im p}&0
},
\end{aligned}
\end{equation}
is split as indicated. Here $N_p$ denotes the kernel of $p$. This
determines the structure of $A\cong N_p\oplus p(A)$ in terms of
the following biproduct:
\begin{equation}
\xymatrix@!C=4pc{
  N_p\ar@<.5ex>[r]^(0.4){\ker p}&
    N_p\oplus p(A)\ar@<.5ex>[l]^(0.6){\coker p}\ar@<.5ex>[r]^(0.6){\coim p}&
    p(A)\ar@<.5ex>[l]^(0.4){\im p}
}.
\end{equation}
The sequence from right to left is split exact, too.
\end{prop}

\begin{thm}
\label{thm_kfrob}
Let $\cal{C}$ be an Abelian symmetric monoidal category and
$(A,\mu,\eta,\Delta,\epsilon)$ be a strongly separable symmetric
Frobenius algebra object in $\cal{C}$ with window element $a$. Then
there exists a knowledgeable Frobenius algebra
$(A,C,\imath,\imath^\ast)$ where $C=p(A)$ is the image of the
idempotent~\eqref{eq_idempotent}, $\imath=\im p$, and
$\imath^\ast=\coim p\circ(a\cdot\id_A)$. The commutative Frobenius
algebra structure of $C$ is given by,
\begin{eqnarray}
  \mu_C      &=& \coim p\circ\mu_A\circ(\im p\otimes\im p),\\
  \eta_C     &=& \coim p\circ\eta_A,\\
  \Delta_C   &=& (\coim p\otimes\coim p)\circ\Delta_A\circ(a\cdot\id_A)\circ\im p,\\
  \epsilon_C &=& \epsilon_A\circ(a^{-1}\cdot\id_A)\circ\im p.
\end{eqnarray}
\end{thm}

\begin{proof}
The proof uses Proposition~\ref{prop_idempotent} and
Proposition~\ref{prop_splitting}.
\end{proof}

We show below in Section~\ref{sect_statesum} that this
knowledgeable Frobenius algebra is precisely the one that is
obtained from our generalized state sum for the strongly separable
algebra $A$. The following proposition introduces two families of
morphisms that are needed in order to show that the morphisms
associated with triangulated open-closed cobordisms do not depend
on the triangulation of the boundary.

Let $(A,\mu,\eta,\Delta,\epsilon)$ be a Frobenius algebra object in a
locally small symmetric monoidal category $\cal{C}$. For $k\in\N$, we
denote by
\begin{equation}
  \mu^{(k+1)}:=\mu\circ(\mu^{(k)}\otimes\id_A),\qquad
  \mu^{(2)}:=\mu,\qquad
  \mu^{(1)}:=\id_A
\end{equation}
and by
\begin{equation}
  \Delta^{(k+1)}:=(\Delta^{(k)}\otimes\id_A)\circ\Delta,\qquad
  \Delta^{(2)}:=\Delta,\qquad
  \Delta^{(1)}:=\id_A
\end{equation}
the iterated multiplication and comultiplication. We also write
$A^{\otimes(k+1)}:=A^{\otimes k}\otimes A$, $A^{\otimes 1}:=A$ and
$A^{\otimes 0}:=\1$, and for $a\in\cal{Z}(A)$,
$a^{k+1}\cdot\id_A:=(a^k\cdot\id_A)\circ(a\cdot\id_A)$ and
$a^0\cdot\id_A:=\id_A$.

\begin{prop} \label{prop_PQ}
Let $\cal{C}$ be a locally small symmetric monoidal category and
$(A,\mu,\eta,\Delta,\epsilon)$ be a strongly separable symmetric
Frobenius algebra object in $\cal{C}$ with window element $a$.
Then for $k,\ell\in\N$, the morphisms
\begin{alignat}{2}
  P_{k\ell}&:=\Delta^{(k)}\circ(a^{-(k-1)}\cdot\id_A)\circ\mu^{(\ell)}&\colon
    A^{\otimes\ell}\to A^{\otimes k},\\
  Q_{k\ell}&:=\Delta^{(k)}\circ(a^{-(k-1)}\cdot\id_A)\circ p\circ\mu^{(\ell)}&\colon
    A^{\otimes\ell}\to A^{\otimes k},
\end{alignat}
satisfy
\begin{equation} \label{eq_projproperty}
  P_{k\ell}\circ P_{\ell m} = P_{km}\qquad\mbox{and}\qquad
  Q_{k\ell}\circ Q_{\ell m} = Q_{km}
\end{equation}
for all $k,\ell,m\in\N$. Here $p$ denotes the idempotent
of~\eqref{eq_idempotent}. In particular, $P_{kk}$ and $Q_{kk}$ are
idempotents, and we have $P_{11}=\id_A$ and $Q_{11}=p$.
\end{prop}

\begin{proof}
In any symmetric Frobenius algebra, we have
\begin{equation}
  \mu^{(k)}\circ\Delta^{(k)}=a^{(k-1)}\cdot\id_A,
\end{equation}
which implies both claims.
\end{proof}

\begin{cor}
\label{cor_blackinv}
Let $\cal{C}$ be an Abelian symmetric monoidal category and
$(A,\mu,\eta,\Delta,\epsilon)$ be a strongly separable symmetric
Frobenius algebra object in $\cal{C}$. Then there are isomorphisms
\begin{equation}
  P_{kk}(A^{\otimes k})\cong A\qquad\mbox{and}\qquad
  Q_{kk}(A^{\otimes k})\cong p(A)
\end{equation}
for all $k\in\N$.
\end{cor}

\begin{proof}
The isomorphisms with their inverses are given by
\begin{alignat}{1}
  \Phi_k=\coim P_{kk}\circ P_{k1}&\colon A\to P_{kk}(A^{\otimes k}),\\
  \Phi_k^{-1}=P_{1k}\circ\im P_{kk}&\colon P_{kk}(A^{\otimes k})\to A
\end{alignat}
as well as
\begin{alignat}{1}
  \Psi_k=\coim Q_{kk}\circ Q_{k1}\circ\im p&\colon p(A)\to Q_{kk}(A^{\otimes k}),\\
  \Psi_k^{-1}=\coim p\circ Q_{1k}\circ \im Q_{kk}&\colon Q_{kk}(A^{\otimes k})\to p(A).
\end{alignat}
\ 
\end{proof}

%
\section{Open-closed cobordisms}
%
\label{sect_topology}

\subsection{Smooth open-closed cobordisms}

In this section, we briefly review the definition of the category
$\cat{2Cob}^{\mathrm{ext}}$ of open-closed cobordisms. These are
smooth $2$-manifolds with corners that have a particular global
structure as follows.

Recall that a \emph{smooth $k$-manifold with corners} $M$ is a
topological $k$-manifold such that every point has a neighbourhood
homeomorphic to an open subset of $\R_+^k:={[0,\infty)}^k$. The
transition functions are required to be the restrictions to $\R_+^k$
of diffeomorphisms between open subsets of $\R^k$.

For each $p\in M$, we define $c(p)\in\N_0$ to be the number of zero
coefficients in local coordinates $\phi(p)\in\R_+^k$. The result is
independent of the chosen coordinate system. A \emph{connected face
of} $M$ is the closure of a component of $\{\,p\in
M\colon\,c(p)=1\,\}$. A \emph{face} is a free union of pairwise
disjoint connected faces. A $k$-dimensional \emph{manifold with faces}
is a smooth $k$-manifold with corners such that every $p\in M$ is
contained in $c(p)$ different connected faces.

A $k$-dimensional $\left<2\right>$-manifold $M$ is a $k$-dimensional
manifold with faces with a specified pair $(\del_0 M,\del_1 M)$ of
faces of $M$ such that $\del_0 M\cup\del_1 M=\del M$ (the boundary of
$M$ as a topological manifold) and such that $\del_0 M\cap\del_1 M$ is
a face of both $\del_0 M$ and $\del_1 M$. A diffeomorphism $f\colon
M\to N$ between $\left<2\right>$-manifolds $M$ and $N$ is a
diffeomorphism of the underlying manifolds with corners that satisfies
$f(\del_0 M)=\del_0 N$ and $f(\del_1 M)=\del_1 N$.

In the following, we are interested in $2$-dimensional
$\left<2\right>$-manifolds. The following is a typical example:
\begin{equation}
\begin{aligned}
\psset{xunit=.4cm,yunit=.4cm}
\begin{pspicture}[0.5](4,5.4)
  \rput(2, 0){\multl}
  \rput(1, 2.5){\medidentl}
  \rput(2.92, 2.5){\ctl}
  \rput(2, 5.25){$M$}
\end{pspicture}
\qquad
\begin{pspicture}[0.5](4,5.4)
  \rput(2, 5.255){$\partial_0 M$}
  \psline[linewidth=0.5pt](1.5,4.5)(.5,4.5)
  \psline[linewidth=0.5pt](2.5,0)(1.5,0)
  \psellipse[linewidth=0.5pt](3,4.5)(.5,0.2)
\end{pspicture}
\qquad
\begin{pspicture}[0.5](4,5.4)
  \rput(2, 5.25){$\partial_1 M$}
  \rput(2,0){
    \psline[linewidth=0.5pt](-1.5,4.5)(-1.5,2.5)
    \psbezier[linewidth=0.5pt](-1.5,2.5)(-1.5,1.5)(-.4,1.3)(-.5,0)
    \psline[linewidth=0.5pt](-.5,4.5)(-.5,2.5)
  }
  \rput(3,0){
    \psbezier[linewidth=0.5pt](-.5,0)(-.6,1.3)(.5,1.5)(.5,2.5)
  }
  \rput(3.08,2.5){
    \psbezier[linewidth=0.5pt](-.58,0)(-.48,.5)(-.48,.7)(-.08,1)
    \psbezier[linewidth=0.5pt](-.08,1)(.32,.7)(.32,.5)(.42,0)
  }
  \rput(2,.5){
    \psbezier[linewidth=0.5pt](-.5,2)(-.6,1)(0.6,1)(.5,2)
  }
\end{pspicture}
\qquad
\begin{pspicture}[0.5](4,5.4)
  \rput(2, 5.25){$\partial_0 M \cup\partial_1 M$}
  \rput(2,0){
     \psline[linewidth=0.5pt](-1.5,4.5)(-1.5,2.5)
     \psbezier[linewidth=0.5pt](-1.5,2.5)(-1.5,1.5)(-.4,1.3)(-.5,0)
     \psline[linewidth=0.5pt](-.5,4.5)(-.5,2.5)
  }
  \rput(3,0){
    \psbezier[linewidth=0.5pt](-.5,0)(-.6,1.3)(.5,1.5)(.5,2.5)
  }
  \rput(3.08,2.5){
    \psbezier[linewidth=0.5pt](-.58,0)(-.48,.5)(-.48,.7)(-.08,1)
    \psbezier[linewidth=0.5pt](-.08,1)(.32,.7)(.32,.5)(.42,0)
  }
  \rput(2,.5){
    \psbezier[linewidth=0.5pt](-.5,2)(-.6,1)(0.6,1)(.5,2)
  }
  \psline[linewidth=0.5pt](1.5,4.5)(.5,4.5)
  \psline[linewidth=0.5pt](2.5,0)(1.5,0)
  \psellipse[linewidth=0.5pt](3,4.5)(.5,0.2)
\end{pspicture}
\qquad
\begin{pspicture}[0.5](4,5.4)
  \rput(2, 5.25){$\partial_0 M \cap\partial_1 M$}
  \psdots(.5,4.5)(1.5,4.5)(1.5,0)(2.5,0)
\end{pspicture}
\end{aligned}
\end{equation}

An \emph{open-closed cobordism} is a compact oriented $2$-dimensional
$\left<2\right>$-manifold $M$ whose distinguished faces we denote by
$(\del_0 M,\del_1 M)$. We call $\del_0 M$ the \emph{black} boundary
and $\del_1 M$ the \emph{coloured} boundary. Two open-closed
cobordisms are considered equivalent if there is an orientation
preserving diffeomorphism of $\left<2\right>$-manifolds that restricts
to the identity on the black boundary.

The black boundary $\del_0 M$ of an open-closed cobordism is
diffeomorphic to a free union of circles $S^1$ and unit intervals
$I=[0,1]$. One can glue open-closed cobordisms along the components of
their black boundary just as one intuitively expects and as it is
indicated by the pictures above. For the details, we refer to
Section~3 of~\cite{LP05}. There we have defined the symmetric monoidal
category $\cat{2Cob}^{\mathrm{ext}}$ as the skeleton of this extended
cobordism category. Its objects are finite sequences
$\vec{n}=(n_1,\cdots,n_k)$ with $n_i\in\{0,1\}$ for all $i$ which
define smooth oriented compact $2$-manifolds $C_{\vec{n}}$ that form
the disjoint union of a unit interval for all $n_i=1$ and of a circle
for all $n_i=0$. The morphisms of $\twocob$ are equivalence classes of
open-closed cobordisms between these $C_{\vec{n}}$. The tensor product
is the free union of manifolds, \ie\ juxtaposition of the
corresponding diagrams, and the composition of morphisms is the gluing
of the open-closed cobordisms along their black boundaries. The
identity morphisms are cylinders over the compact oriented
$1$-manifolds that represent the objects.

For the purpose of the present work, it is sufficient to keep in mind
that any morphism of the category $\cat{2Cob}^{\mathrm{ext}}$ can be
obtained from a finite number of copies of the
generators~\eqref{eq_generators} by taking tensor products and by
taking compositions.

We have also shown in~\cite{LP05} that any two equivalent
open-closed cobordisms are related by a finite sequence of moves.
These moves are precisely the defining equations of a
knowledgeable Frobenius algebra $(A,C,\imath,\imath^\ast)$ when
the operations $\mu_A,\Delta_A,\ldots,\imath^\ast$ are replaced by
the morphisms depicted in~\eqref{eq_generators}. For example, to
the associative law in the symmetric Frobenius algebra
$(A,\mu_A,\eta_A,\Delta_A,\epsilon_A)$, there corresponds the
following move:
\begin{equation}
\begin{aligned}
\psset{xunit=.3cm,yunit=.3cm}
\begin{pspicture}[0.5](4,5.5)
  \rput(2,0){\multl}
  \rput(3,2.5){\multl}
  \rput(1,2.5){\curveleftl}
\end{pspicture}
\quad\longleftrightarrow\quad
\begin{pspicture}[0.5](4,5.5)
  \rput(2,0){\multl}
  \rput(1,2.5){\multl}
  \rput(3,2.5){\curverightl}
\end{pspicture}
\end{aligned}
\end{equation}
We can summarize the results of~\cite{LP05} on the structure of the
category $\cat{2Cob}^{\mathrm{ext}}$ as follows.

\begin{thm}
The category $\cat{2Cob}^{\mathrm{ext}}$ of open-closed cobordisms is
the strict symmetric monoidal category freely generated by a
knowledgeable Frobenius algebra object
$(A,C,\imath,\imath^\ast)$. This generating knowledgeable Frobenius
algebra object consists of the diffeomorphism type $\vec n=(1)$ of the
unit interval $A=I$ which forms a symmetric Frobenius algebra
$(A,\mu_A,\eta_A,\Delta_A,\epsilon_A)$; the diffeomorphism type $\vec
n=(0)$ of the circle $C=S^1$ which forms a commutative Frobenius
algebra $(C,\mu_C,\eta_C,\Delta_C,\epsilon_C)$; together with the
morphisms $\imath$ and $\imath^\ast$. The morphisms
$\mu_A,\eta_A,\ldots,\imath^\ast$ are precisely the equivalence
classes of the open-closed cobordisms depicted
in~\eqref{eq_generators}.
\end{thm}

\subsection{Combinatorial open-closed cobordisms}

Open-closed cobordisms can be triangulated as follows. We use the
terminology of~\cite{RS}.

Given an open-closed cobordism $M$, the underlying topological
manifold is a compact oriented $2$-manifold with boundary. We
therefore have a finite simplicial complex $K$ whose underlying
polyhedron we denote by $|K|\subseteq\R^p$ for some $p$, and a
homeomorphism $T_M\colon |K|\to M$ which we call a
\emph{triangulation}. The simplicial complex $K$ satisfies the
conditions that guarantee that $|K|$ forms an oriented topological
$2$-manifold, \ie\ the link of each $d$-simplex is a
$(1-d)$-sphere iff the simplex is in the interior of $|K|$, and it
is a $(1-d)$-ball iff the simplex is in the boundary of $|K|$.
Furthermore, for each $2$-simplex $\sigma$, it is specified
whether $\sigma$ or its opposite oriented simplex $\sigma^\ast$ is
contained in $|K|$, and each $1$-simplex in the interior of $|K|$
appears as a face of precisely two $2$-simplices with opposite
induced orientations.

If $M$ and $N$ are equivalent open-closed cobordisms, their
underlying topological manifolds are homeomorphic. If we have
triangulations $T_M\colon|K|\to M$ and ${\tilde T}_M\colon|L|\to
N$ with simplicial complexes $K$ and $L$, Pachner's
theorem~\cite{Pachner} says that the simplicial complexes $K$ and
$L$ are related by a finite sequence of moves. These moves are the
\emph{bistellar moves} (called the 1-3 and 2-2 move),
\begin{equation}
\label{eq_bistellar}
\begin{aligned}
\psset{xunit=.4cm,yunit=.4cm}
\begin{pspicture}(4,3)
  \pspolygon[fillstyle=solid,fillcolor=lightgray](0,0)(4,0)(2,2.8)(0,0)
  \psdots(0,0)(4,0)(2,2.8)
\end{pspicture}
\end{aligned}
\quad\longleftrightarrow\quad
\begin{aligned}
\psset{xunit=.4cm,yunit=.4cm}
\begin{pspicture}(5,3)
  \pspolygon[fillstyle=solid,fillcolor=lightgray](0,0)(4,0)(2,2.8)(0,0)
  \psdots(0,0)(4,0)(2,2.8)
  \psline(2,1.1)(0,0)
  \psline(2,1.1)(4,0)
  \psline(2,1.1)(2,2.8)
  \psdots(2,1.1)
\end{pspicture}
\end{aligned}
\quad\mbox{and}\quad
\begin{aligned}
\psset{xunit=.4cm,yunit=.4cm}
\begin{pspicture}(3,3.5)
  \pspolygon[fillstyle=solid,fillcolor=lightgray](0,0)(3,0)(3,3)(0,3)(0,0)
  \psline(0,0)(3,3)
  \psdots(0,0)(3,0)(3,3)(0,3)(0,0)
\end{pspicture}
\end{aligned}
\quad\longleftrightarrow\quad
\begin{aligned}
\psset{xunit=.4cm,yunit=.4cm}
\begin{pspicture}(3.5,3.5)
 \pspolygon[fillstyle=solid,fillcolor=lightgray](0,0)(3,0)(3,3)(0,3)(0,0)
 \psline(3,0)(0,3)
 \psdots(0,0)(3,0)(3,3)(0,3)(0,0)
\end{pspicture}
\end{aligned},
\end{equation}
applicable to all $2$-simplices, and the \emph{elementary shellings}
\begin{equation}
\label{eq_shelling}
\begin{aligned}
\begin{pspicture}(1,3)
  \pspolygon[fillstyle=solid,fillcolor=lightgray, linecolor=lightgray](0,0)(1,0)(1,3)(0,3)(0,0)
  \psline(1,0)(1,3)
  \psdots(1,0)(1,1)(1,2)(1,3)
\end{pspicture}
\end{aligned}
\quad\longleftrightarrow\quad
\begin{aligned}
\begin{pspicture}(2,3)
  \pspolygon[fillstyle=solid,fillcolor=lightgray, linecolor=lightgray](0,0)(1,0)(1,3)(0,3)(0,0)
  \pspolygon[fillstyle=solid,fillcolor=lightgray](1,1)(1,2)(1.8,1.5)(1,1)
  \psline(1,0)(1,3)
  \psdots(1,0)(1,1)(1,2)(1,3)(1.8,1.5)
\end{pspicture}
\end{aligned}
\quad\mbox{and}\quad
\begin{aligned}
\begin{pspicture}(1,3)
  \pspolygon[fillstyle=solid,fillcolor=lightgray, linecolor=lightgray](0,0)(1,0)(1,3)(0,3)(0,0)
  \psline(1,0)(1,3)
  \psline(1,1)(.3,1.5)
  \psline(1,2)(.3,1.5)
  \psdots(1,0)(1,1)(1,2)(1,3)(.3,1.5)
\end{pspicture}
\end{aligned}
\quad\longleftrightarrow\quad
\begin{aligned}
\begin{pspicture}(1.5,3)
  \pspolygon[fillstyle=solid,fillcolor=lightgray,linecolor=lightgray](1,0)(0,0)(0,3)(1,3)(1,0)
  \pspolygon[fillstyle=solid,fillcolor=white,linecolor=white](1,1)(1,2)(.3,1.5)(1,1)
  \psline(1,0)(1,1)
  \psline(1,1)(.3,1.5)
  \psline(1,2)(.3,1.5)
  \psline(1,2)(1,3)
  \psdots(1,0)(1,1)(1,2)(1,3)(.3,1.5)
\end{pspicture}
\end{aligned},
\end{equation}
applicable to certain $2$-simplices some of whose faces coincide
with the boundary. The interior of the manifold is indicated by
the shading in our pictures. Recall that for finite simplicial
complexes which represent compact manifolds with non-empty
boundary, each bistellar move can be obtained from a finite
sequence of elementary shellings.

The set of corners $\del_0 M\cap\del_1 M$ of every open-closed
cobordism $M$ is a finite set. Given some triangulation
$T_M\colon|K|\to M$, we can apply a finite sequence of elementary
shellings in order to subdivide the $1$-simplices in the boundary
in such a way that to every corner of $M$, there corresponds a
$0$-simplex in $K$, \ie\ that $\del_0 M\cap\del_1 M\subseteq
T_M(|K_0|)$ where $K_0\subseteq K$ denotes the $0$-skeleton of
$K$. From now on we assume, without loss of generality, that every
triangulation has this property. Given a $1$-simplex $\sigma\in K$
in the boundary, we therefore have either
$T_M(|\sigma|)\subseteq\del_0 M$ or $T_M(|\sigma|)\subseteq\del_1
M$, \ie\ the $1$-simplices in the boundary are either \emph{black}
or \emph{coloured}.

Both elementary shellings of~\eqref{eq_shelling} replace two
boundary $1$-simplices (edges) by a single edge or vice versa. For
triangulations with the special property, each of the elementary
shellings~\eqref{eq_shelling} belongs to one of the following four
types:
\begin{enumerate}
\item
  two black edges $\longleftrightarrow$ one black edge,
\item
  two coloured edges $\longleftrightarrow$ one coloured edge,
\item
  one black and one coloured edge $\longleftrightarrow$ one black edge,
\item
  one black and one coloured edge $\longleftrightarrow$ one coloured edge.
\end{enumerate}
It is not difficult to see that the elementary shellings of type
(3.) and (4.) can be obtained from a finite sequence of bistellar
moves and elementary shellings of type (1.) and (2.).

When we construct open-closed TQFTs in Section~\ref{sect_statesum}
below, we consider triangulations of the open-closed cobordisms
and then show that the linear map associated with every given
cobordism is invariant under the bistellar
moves~\eqref{eq_bistellar} and under elementary shellings of type
(1.) and (2.). Then this linear map is independent of the choice
of the triangulation.

\subsection{Smoothing theory}

When one studies smooth manifolds by combinatorial techniques, the
relation between combinatorial and smooth manifolds is described by
two types of theorems:
\begin{itemize}
\item
  \emph{Triangulation:} Every compact smooth manifold (with boundary)
  admits a Whitehead triangulation. If two such manifolds are
  diffeomorphic, then their triangulations are related by a finite
  sequence of the appropriate Pachner moves.
\item
  \emph{Smoothing:} Given a finite simplicial complex $K$ that
  satisfies the conditions which ensure that its underlying polyhedron
  $|K|$ forms a topological manifold (with boundary), one needs to
  know (a) under which conditions there exists a smooth manifold that
  has $|K|$ as its triangulation and (b) whether the resulting smooth
  manifold is unique up to diffeomorphism.
\end{itemize}
Such theorems are available in order to compare smooth manifolds
with boundary and combinatorial manifolds with boundary, but we
are not aware of any systematic treatment for manifolds with
corners, manifolds with faces, or $\left<2\right>$-manifolds.

In the preceding section, we have solved the triangulation problem
for open-closed cobordisms by resorting to the underlying
topological manifold which is just a topological $2$-manifold with
boundary. It admits a triangulation, and this triangulation is
unique up to combinatorial equivalence, \ie\ Pachner moves, by the
validity of the Combinatorial Triangulation Conjecture and the
\emph{Hauptvermutung} for $2$-dimensional manifolds, see, for
example~\cite{Ranicki}. We have then dealt with the corner points
`by hand'.

The other direction, a solution to the smoothing problem, is not
needed if one is just interested in a combinatorial construction of
open-closed TQFTs. For completeness, we nevertheless sketch how one
can obtain the corresponding smoothing theorem: Let $K$ be a finite
simplicial complex that triangulates an open-closed cobordism. Then
every $1$-simplex in the boundary is either black or coloured as we
have explained above. The underlying polyhedron $|K|$ together with
this partitioning of the boundary is already sufficient to read off
the topological invariants defined in Section~3.2.4 of~\cite{LP05}. By
the normal form of open-closed cobordisms of Definition 3.18
of~\cite{LP05}, there exists an open-closed cobordism with the given
invariants, and by Corollary~3.24 of~\cite{LP05}, it is unique up to
equivalence.

%
\section{State Sum Construction}
%
\label{sect_statesum}

We begin this section with an overview of the state sum
construction in informal language.

Given a strongly separable symmetric Frobenius algebra 
$(A,\mu_A,\eta_A,\Delta_A,\epsilon_A)$ in an Abelian symmetric
monoidal category $\cal{C}$ and a connected open-closed cobordism $M$
with triangulation $T_M\colon|K|\to M$, we construct a morphism $Z(M)$
in $\cal{C}$.

\piccaption{
  $\vec{n}=(1,0)$, $\vec{n}'=(1)$, $h_1=2$, $h_2=5$, $h_3=4$,
  $m_1=7$, and $m_2=4$.
}
\parpic[r]{$\includegraphics{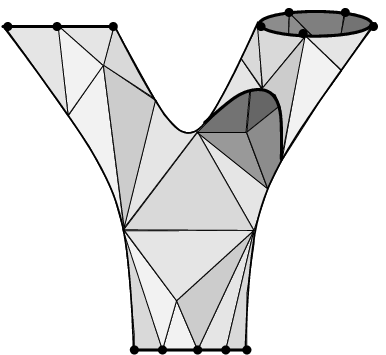}$
} 
For the duration of this section let $M$ be a connected open-closed
cobordism with source $\partial_0 M^{{\rm in}}:=\vec{n}=(n_1,
\cdots , n_k)$ and target $\partial_0 M^{{\rm
out}}:=\vec{n}'=(n'_1, \cdots , n'_{k'})$. Let $j$ enumerate the
black boundary components of $M$ so that $h_j$ denotes the number
of 1-simplices in the triangulation of the component $n_j$ for $1
\leq j \leq k$ or the component $n_j'$ for $k+1 \leq j \leq k+k'$.
The number of 1-simplices of $\partial_0 M^{{\rm in}}$ is given by
the sum $m_1:=\sum_{j=1}^{k}h_j$, and the number of 1-simplices of
$\partial_0 M^{{\rm out}}$ by the sum
$m_2:=\sum_{j=k+1}^{k+k'}h_j$.

As a first step to constructing the morphism $Z(M)$, we construct
a morphism $Z_{T_M}(M)\colon A^{\otimes m_1}\to A^{\otimes m_2}$.
These morphisms depend on the triangulation of the black boundary,
but they are already invariant under bistellar moves and under
elementary shellings of type~(2.), \ie\ those in which all the
involved boundary edges are coloured.

Define the symbol $A^{(n_j)}$ corresponding to the boundary
component $n_j$ to be $A$ if $n_j=1$ and $p(A)$ if $n_j=0$ and
define $A^{\otimes \vec{n}}$ to be the ordered tensor product $
\bigotimes_{j=1}^k A^{(n_j)} $.  Likewise, we set $A^{\otimes
\vec{n}'}$ equal to the ordered tensor product
$\bigotimes_{j=k+1}^{k+k'} A^{(n'_j)}$.

In Section~\ref{sec_indep-black}, we show that the isomorphisms
$P_{kk}(A^{\otimes k}) \cong A$ and $Q_{kk}(A^{\otimes k}) \cong
p(A)$ of Corollary~\ref{cor_blackinv} correspond to triangulated
cylinders over $I$ or $S^1$. We construct a map $Z(M)\colon
A^{\otimes \vec{n}}\to A^{\otimes \vec{n'}}$ using these
isomorphisms and the morphism $Z_{T_M}(M)$.  Since the claim of
Corollary~\ref{cor_blackinv} is independent of $k$, and since the
isomorphisms used in that corollary correspond to triangulated
cylinders over $I$ or $S^1$, the invariance under bistellar moves
and elementary shellings of type (2.) can be used to show
independence of the boundary triangulation. The morphism $Z(M)$ is
then also invariant under elementary shellings of type~(1.), \ie\
those involving the black boundary. $Z(M)$ is therefore
independent of the triangulation and thus well-defined for the
open-closed cobordism $M$.

One can verify explicitly that composition and disjoint union work
as required, and so the state sum defines an open-closed TQFT
$Z\colon\cat{2Cob}^{\mathrm{ext}}\to\cal{C}$. The objects of
$\cal{C}$ associated with the interval and the circle are $A$ and
$p(A)$, respectively, by construction. What is the knowledgeable
Frobenius algebra that characterizes this TQFT?

In order to answer this question, we compute the morphisms of
$\cal{C}$ associated with the generating open-closed
cobordisms~\eqref{eq_generators} and show that the open-closed
TQFT is characterized by the knowledgeable Frobenius algebra of
Theorem~\ref{thm_kfrob}.

\subsection{Defining the state sum}

\begin{figure}
$ \xy
 (-48.4,34.1)*{\scriptscriptstyle a^{ ^{\txt\tiny{-1} }} };
 (-29.6,34)*{\scriptscriptstyle a^{ ^{\txt\tiny{-1} }}};
 (-37.7,-47)*{\scriptscriptstyle a^{ ^{\txt\tiny{-1} }}};
 (0,0)*{\includegraphics{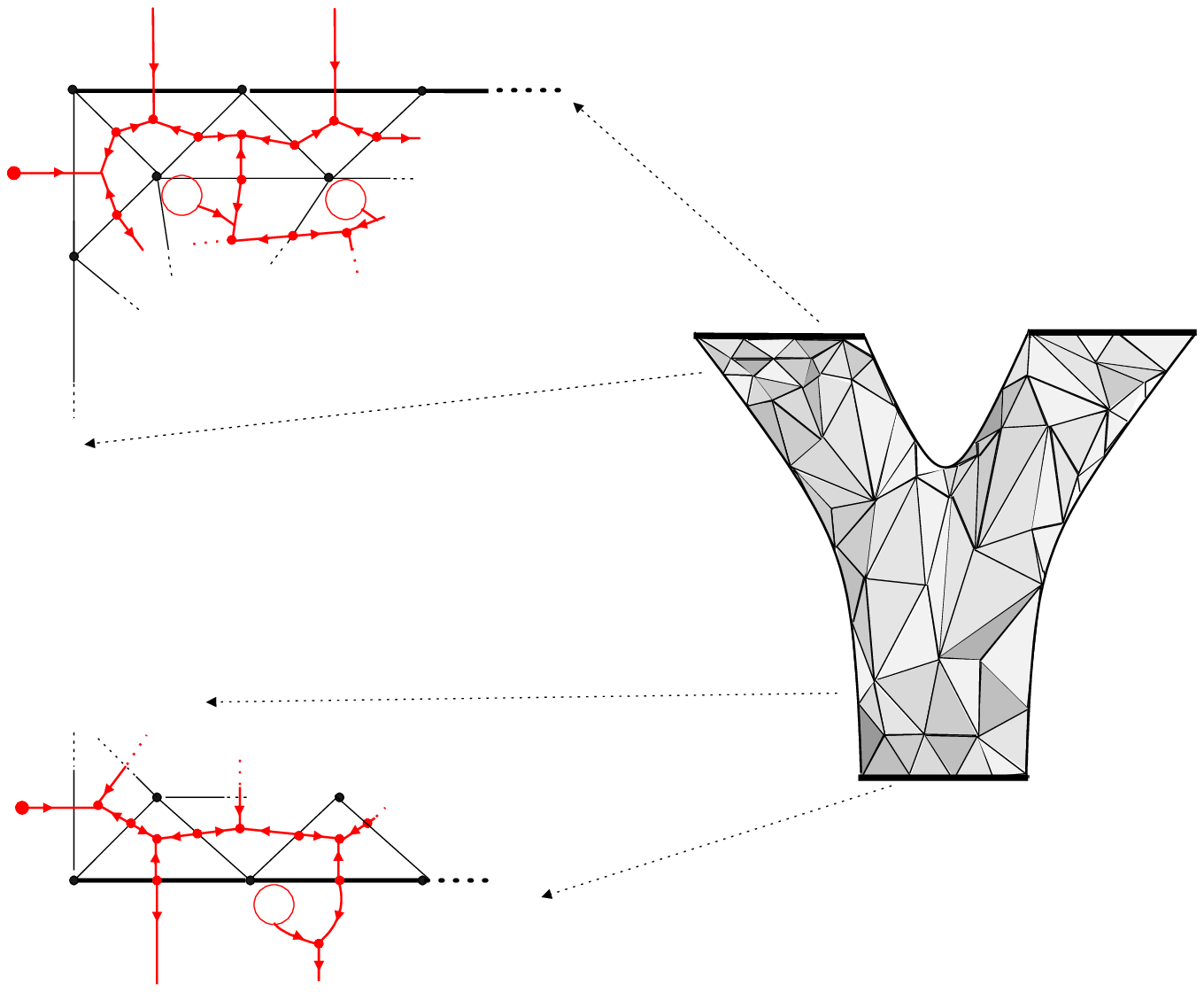}}; (-33,55)*{A}; (-54,55)*{A};
 (-35,-55)*{A}; (-56,-55)*{A}; (-72,35)*{\eta_A};
 (-72,-37)*{\eta_A}; \endxy $
\caption{
  This figure illustrates the state sum for an open-closed cobordisms $M$.
}\label{corner75}
\end{figure}

We first describe how to construct the morphism $Z_{T_M}(M)\colon
A^{\otimes m_1}\to A^{\otimes m_2}$. It is defined by a string
diagram in $\cal{C}$ obtained from the graph Poincar\'{e} dual to
the triangulation, see Figure~\ref{corner75}. By the coherence
theorem for symmetric monoidal categories, it does not matter how
one projects the Poincar{\'e} dual graph onto the drawing plane.

For every $2$-simplex (triangle), we put a `trilinear form'
$g^{(3)}$ (\cf~\eqref{eq_trilinear}), and for every edge in the
interior, we have an inverse bilinear form
$g^\ast=\Delta_A\circ\eta_A$. Note that $g^{(3)}$ has a symmetry
under the cyclic group $C_3$, but not in general under the
symmetric group $S_3$, and so this assignment depends on the
orientation.

For every edge on the coloured boundary $\del_1 M$, we put a unit
$\eta_A$. For every interior $0$-simplex (vertex), we multiply the
resulting morphism by the inverse $a^{-1}$ of the window element.
Since $a^{-1}$ is central and the cobordism connected, it does not
matter where in the diagram we do this.

At this stage, we have a morphism $A^{\otimes (m_1+m_2)}\to\1$ of
$\cal{C}$. Finally, for every edge in the black out-boundary
$\del_0 M^{\mathrm{out}}$, we put a $g^\ast$, too, in order to
turn this into a morphism $A^{\otimes m_1}\to A^{\otimes m_2}$.
Then, for every vertex in the black out-boundary that is not a
corner, we multiply by $a^{-1}$.

The terminology \emph{sum} in `state sum' is justified by the
following point of view: If $\cal{C}=\cat{Vect}_k$ and if one
chooses a basis of $A$ and expands all linear maps in this basis,
the state sum contains a sum over the basis vectors for each edge
in the interior of $M$. This is the \emph{sum} involved in the
state sum.

The morphisms specified by the string diagram have two important
properties.
\begin{itemize}
\item
  Gluing triangulated open-closed cobordisms along a common black
  boundary that is triangulated with the same number of edges,
  corresponds to the composition of morphisms.
\item
  The disjoint union of open-closed cobordisms gives the tensor
  product of morphisms.
\end{itemize}

\noindent
The definition reads in detail as follows.

\begin{defn}
\label{def_defstatesum}
Let $(A,\mu_A,\eta_A,\Delta_A,\epsilon_A)$ be a strongly separable
symmetric Frobenius algebra in an Abelian symmetric monoidal category
$\cal{C}$. Let $M$ be an open-closed cobordism with triangulation
$T_M\colon|K|\to M$. Let $K^{(j)}\subseteq K$ denote the set of
$j$-simplices, $j\in\{0,1,2\}$.

We characterize the edges, \ie\ the elements $\sigma_{\{i,j\}}\in
K^{(1)}$, by two-element sets $\{i,j\}\subseteq K^{(0)}$, $i\neq j$,
of vertices. The oriented triangles $\sigma_{(i,j,k)}\in K^{(2)}$ are
characterized by triples $(i,j,k)\in K^{(0)}\times K^{(0)}\times
K^{(0)}$ of vertices, modulo a permutation by a $3$-cycle.

We define the morphism $Z_{T_M}(M)\colon A^{\otimes m_1}\to A^{\otimes
m_2}$ as a composition
\begin{equation}
\label{eq_threeparts}
  Z_{T_M}(M) := Z_{T_M}^{(2)}\circ (a^{-k}\cdot\id_{A^{\otimes N}})
    \circ\tau\circ Z_{T_M}^{(1)}.
\end{equation}
where $N=m_2+|\{\,\sigma\in K^{(1)}\colon\,\sigma\subseteq\del
M\,\}| +2|\{\,\sigma\in K^{(1)}\colon\,\sigma\subseteq
M\backslash\del M\,\}| =m_2+3|K^{(2)}|$. The power of the inverse
window element in~\eqref{eq_threeparts} is $k=|\{\,\sigma\in
K^{(0)}\colon\,\sigma\subseteq M\backslash\del M\,\}|
+|\{\,\sigma\in K^{(0)}\colon\,\sigma\subseteq\del_0
M^\mathrm{out}\backslash(\del_0 M\cap\del_1 M)\,\}|$ --- the
number of interior vertices plus the number of vertices on the
outgoing edge that are not corners. We exploit the coherence
theorem for monoidal categories and suppress the associativity and
unit constraints of $\cal{C}$ and define
\begin{equation}
\label{eq_firstpart}
  Z_{T_M}^{(1)} := \biggl(\bigotimes_{j=1}^{m_1}\id_A\biggr)
    \otimes\biggl(\bigotimes_{j=1}^{m_2} g^\ast\biggr)
    \otimes\biggl(\bigotimes_{\ontop{\sigma\in K^{(1)}\colon}{\sigma\subseteq M\backslash\del M}}g^\ast\biggr)
    \otimes\biggl(\bigotimes_{\ontop{\sigma\in K^{(1)}\colon}{\sigma\subseteq\del_1 M}}\eta_A\biggr)\colon
    A^{\otimes m_1}\to A^{\otimes N}.
\end{equation}
and
\begin{equation}
\label{eq_secondpart}
  Z_{T_M}^{(2)} := \biggl(\bigotimes_{j=1}^{m_2}\id_A\biggr)
    \otimes\biggl(\bigotimes_{\sigma\in K^{(2)}}g^{(3)}\biggr)\colon
    A^{\otimes N}\to A^{\otimes m_2},
\end{equation}
The morphism $\tau\colon A^{\otimes N}\to A^{\otimes N}$ permutes
the tensor factors. In order to specify this permutation, we
associate the factors of the target of~\eqref{eq_firstpart} and
those of the domain of~\eqref{eq_secondpart} with the edges
$\sigma_{\{i,j\}}\in K^{(1)}$. This is denoted by superscripts
such as $A^{\{i,j\}}$. The permutation $\tau$ is specified by
requiring that it maps each factor $A^{\{i,j\}}$ to one whose
superscript is the same edge.

The superscripts for the $A$'s in the target of~\eqref{eq_firstpart}
are as follows. We go through the factors of~\eqref{eq_firstpart} from
left to right.
\begin{itemize}
\item
  For every edge $\sigma_{\{i,j\}}$ in the black in-boundary $\del_0
  M^{\mathrm{in}}$, we have $\id_A\colon A\to A^{\{i,j\}}$. There are
  $m_1$ edges of this sort.
\item
  For every edge $\sigma_{\{i,j\}}$ in the black out-boundary $\del_0
  M^{\mathrm{out}}$, we have $g^\ast\colon\1\to A^{\{i,j\}}\otimes
  A^{\{i,j\}}$. This edge therefore appears twice as a superscript,
  but due to the symmetry of $g^\ast$, we need not distinguish the
  two. There are $m_2$ edges of this sort.
\item
  For every edge $\sigma_{\{i,j\}}\subseteq M\backslash\del M$ in the
  interior, we have $g^\ast\colon\1\to A^{\{i,j\}}\otimes
  A^{\{i,j\}}$. Again the superscript occurs twice, and we do not
  distinguish.
\item
  For every edge $\sigma_{\{i,j\}}\subseteq\del_1 M$ in the coloured
  boundary, we have $\eta_A\colon\1\to A^{\{i,j\}}$.
\end{itemize}
The superscripts for the $A$'s in the domain of~\eqref{eq_secondpart}
are as follows.
\begin{itemize}
\item
  For every edge $\sigma_{\{i,j\}}$ in the black out-boundary $\del_0
  M^{\mathrm{out}}$, we have $\id_A\colon A^{\{i,j\}}\to A$.
\item
  For every oriented triangle $\sigma_{(i,j,k)}\in K^{(2)}$, we have
  $g^{(3)}\colon A^{\{i,j\}}\otimes A^{\{j,k\}}\otimes
  A^{\{k,i\}}\to\1$. Due to the cyclic symmetry of the `trilinear
  form' $g^{(3)}$, this morphism is invariant under permutations of
  the triple $(i,j,k)$ by a $3$-cycle.
\end{itemize}
Notice that the edges that appear as superscripts in the target
of~\eqref{eq_firstpart} and those in the domain
of~\eqref{eq_secondpart} agree including their multiplicities, and
that the permutation $\tau$ is well defined.
\end{defn}

\noindent
See~\eqref{eq_exampleSS} for an example of the diagram produced by the
state sum.

\subsection{Invariance under Pachner moves}

\begin{prop}
\label{prop_interiorcoloured} 
For a connected open-closed cobordism $M$ with triangulation $T_M$,
the state sum $Z_{T_M}(M)$ is invariant under the 1-3 and 2-2 Pachner
moves and under the elementary shellings of type (2.).
\end{prop}

\begin{proof}
The 2-2 Pachner move follows from the cyclic symmetry of the
`trilinear form' $g^{(3)}$.
\begin{equation}
 \psset{xunit=.5cm,yunit=.5cm}
 \begin{pspicture}[.5](3.2,3.5)
  \pspolygon[fillstyle=solid,fillcolor=lightgray](1.5,0)(3,1.5)(1.5,3)(0,1.5)(1.5,0)
  \psline(1.5,0)(1.5,3)
  \psdots(1.5,0)(3,1.5)(1.5,3)(0,1.5)
  \psline[linecolor=red,linewidth=.8pt](.75,2.25)(1,1.5)
  \psline[linecolor=red,linewidth=.8pt](.75,.75)(1,1.5)
  \psline[linecolor=red,linewidth=.8pt](1,1.5)(2,1.5)
  \psline[linecolor=red,linewidth=.8pt](2.25,2.25)(2,1.5)
  \psline[linecolor=red,linewidth=.8pt](2.25,.75)(2,1.5)
  \psdots[linecolor=red](1,1.5)(1.5,1.5)(2,1.5)
  \rput(.7,2.6){$\scriptstyle i$}
  \rput(.7,.4){$\scriptstyle j$}
  \rput(2.3,2.6){$\scriptstyle l$}
  \rput(2.3,.4){$\scriptstyle k$}
    \end{pspicture}
\quad = \quad
 \xy
 (0,-8)*\bbtriple{};
 (-4,0)*\bbainv{4};
 (0,0)*\bbdl{};
 (6,0)*\bbdual{};
 (12,-8)*\bbtriple{};
 (12,0)*\bbdr{};
 (16,0)*\bbdr{};
 (14,8)*\bbid{i};
 (10,8)*\bbid{};
 (2,8)*\bbid{};
 (-2,8)*\bbid{};
 (-3,7)*{\scriptstyle k};
 (1,7)*{\scriptstyle l};
 (9,7)*{\scriptstyle i};
 (13,7)*{\scriptstyle j};
 \endxy
 \quad = \quad
 \xy
 (0,-8)*\bbtriple{};
 (-4,0)*\bbainv{4};
 (0,0)*\bbdl{};
 (6,0)*\bbdual{};
 (12,-8)*\bbtriple{};
 (12,0)*\bbdr{};
 (16,0)*\bbdr{};
 (14,8)*\bbid{i};
 (10,8)*\bbid{};
 (2,8)*\bbid{};
 (-2,8)*\bbid{};
 (-3,7)*{\scriptstyle l};
 (1,7)*{\scriptstyle i};
 (9,7)*{\scriptstyle j};
 (13,7)*{\scriptstyle k};
 \endxy
 \quad = \quad
 \psset{xunit=.5cm,yunit=.5cm}
 \begin{pspicture}[.5](3.2,3.5)
  \pspolygon[fillstyle=solid,fillcolor=lightgray](1.5,0)(3,1.5)(1.5,3)(0,1.5)(1.5,0)
  \psline(0,1.5)(3,1.5)
  \psdots(1.5,0)(3,1.5)(1.5,3)(0,1.5)
   \psline[linecolor=red,linewidth=.8pt](.75,2.25)(1.5,2)
  \psline[linecolor=red,linewidth=.8pt](.75,.75)(1.5,1)
  \psline[linecolor=red,linewidth=.8pt](1.5,2)(1.5,1)
  \psline[linecolor=red,linewidth=.8pt](2.25,2.25)(1.5,2)
  \psline[linecolor=red,linewidth=.8pt](2.25,.75)(1.5,1)
  \psdots[linecolor=red](1.5,2)(1.5,1.5)(1.5,1)
    \rput(.7,2.6){$\scriptstyle i$}
  \rput(.7,.4){$\scriptstyle j$}
  \rput(2.3,2.6){$\scriptstyle l$}
  \rput(2.3,.4){$\scriptstyle k$}
    \end{pspicture}
\end{equation}
The 1-3 Pachner move is slightly more difficult because it
involves subdividing a triangle which inserts an additional
internal vertex. It makes use of the bubble
move~\eqref{eq_bubble}:
\begin{equation}
\psset{xunit=.4cm,yunit=.4cm}
\begin{pspicture}[.5](5,3)
  \pspolygon[fillstyle=solid,fillcolor=lightgray](0,0)(4,0)(2,2.8)(0,0)
  \psdots(0,0)(4,0)(2,2.8)
         \psline(2,1.1)(0,0)
         \psline(2,1.1)(4,0)
         \psline(2,1.1)(2,2.8)
        \psdots(2,1.1)
\end{pspicture}
\quad = \quad
\xy
 (0,-8)*\bbtriple{};
 (-4,0)*\bbainv{4};
 (0,0)*\bbdl{};
 (6,0)*\bbdual{};
 (12,-8)*\bbtriple{};
 (12,0)*\bbdr{};
 (18,0)*\bbdual{};
 (24,-8)*\bbtriple{};
 (24,0)*\bbdr{};
 (28,0)*\bbdr{};
 (10,8)*\bbid{};
 (2,8)*\bbid{};
 (22,8)*\bbid{};
 (10,16)*\bbid{};
 (2,16)*\bbid{};
 (22,16)*\bbid{};
 (2,8)*\bblrlong{};
 (22,8)*\bbrllong{};
 (12,16)*\bbhugedual{};
\endxy
\quad = \quad
\vcenter{\xy
 (0,0)*\bbtriple{};
 (0,0)*\bbid{};
 (-2,-8)*\bbainv{4};
 (-2,-16)*\bbdelta{};
 (-2,-24)*\bbpair{};
\endxy}\nn
\end{equation}
\begin{equation}
 = \quad
\vcenter{\xy
 (0,-8)*\bbtriple{};
 (-4,0)*\bbainv{3};
 (-2,8)*\bbid{};
  (0,0)*\bbdl{};
 (2,8)*\bbid{};
  (4,0)*\bbdl{};
 (6,8)*\bbid{};
\endxy}
\quad = \quad
\psset{xunit=.4cm,yunit=.4cm}
\begin{pspicture}[.5](4,3)
  \pspolygon[fillstyle=solid,fillcolor=lightgray](0,0)(4,0)(2,2.8)(0,0)
  \psdots(0,0)(4,0)(2,2.8)
\end{pspicture}
\end{equation}

\parpic[r]{
$ \psset{xunit=.5cm,yunit=.5cm}
\begin{aligned}
 \begin{pspicture}(2,3)
  \pspolygon[fillstyle=solid,fillcolor=lightgray, linecolor=lightgray](0,0)(2,0)(2,3)(0,3)(0,0)
  \psline(2,0)(2,3)
  \psdots(2,.75)(2,2.25)(1,1.5)
    \psline(2,.75)(1,1.5)
    \psline(2,2.25)(1,1.5)
    \psline(1,1.5)(1,2.25)
    \psline[linestyle=dotted](1,1.5)(1,2.65)
    \psline(1,1.5)(.4,1.8)
    \psline[linestyle=dotted](1,1.5)(.1,1.95)
    \psline(1,1.5)(.7,.9)
    \psline[linestyle=dotted](1,1.5)(.5,.5)
    \end{pspicture}
\end{aligned}
\qquad \leftrightsquigarrow \qquad
\begin{aligned}
  \begin{pspicture}(2.1,3)
  \pspolygon[fillstyle=solid,fillcolor=lightgray, linecolor=lightgray](0,0)(2,0)(2,3)(0,3)(0,0)
  \pspolygon[fillstyle=solid,fillcolor=white,linecolor=white](2,.75)(2,2.25)(1,1.5)(2,.75)
  \psline(2,.75)(2,0)
   \psline(2,2.25)(2,3)
  \psdots(2,.75)(2,2.25)(1,1.5)
    \psline(2,.75)(1,1.5)
    \psline(2,2.25)(1,1.5)
    \psline(1,1.5)(1,2.25)
    \psline[linestyle=dotted](1,1.5)(1,2.65)
    \psline(1,1.5)(.4,1.8)
    \psline[linestyle=dotted](1,1.5)(.1,1.95)
    \psline(1,1.5)(.7,.9)
    \psline[linestyle=dotted](1,1.5)(.5,.5)
    \end{pspicture}
\end{aligned}
$} There are two elementary shellings~\eqref{eq_shelling} of type
(2.). Recall that the state sum assigns to each edge of the coloured
boundary the algebra unit $\eta_A\maps\1\to A$. The first move
of~\eqref{eq_shelling} follows directly from the unit axioms. The
second move turns an interior vertex into an exterior vertex (featured
to the right). This move follows from the bubble
move~\eqref{eq_bubble}:
\begin{equation}
 \xy
 (-2,-8)*\bbainv{1};
 (0,0)*\bbmu{};
 (0,8)*\bbdelta{};
 (0,16)*\bbeta{};
 \endxy
 \quad = \quad
  \xy
 (0,0)*\bbmu{};
 (2,8)*\bbeta{};
 (-2,8)*\bbeta{};
 \endxy
\end{equation}
\end{proof}

Note that the bubble move~\eqref{eq_bubble} is required to prove the
above proposition. This is the reason why we cannot define the state
sum for the non strongly separable algebras of
Example~\ref{ex_groupalg}(2) and~\ref{ex_matrixalg}(2).

For convenience, we sometimes use degenerate triangulations in which
the two vertices in the boundary of an edge agree. In this case it is
always understood that we apply bistellar moves and elementary
shellings in order to turn them into proper simplicial complexes.

An example showing the diagram produced by the state sum on the
torus $T^2$ is depicted below:
\begin{equation} \label{eq_exampleSS}
\psset{xunit=.4cm,yunit=.4cm}
        \begin{pspicture}[.4](4,4)
        \rput(2.5,0){\zigc}
        \rput(2.65,1.99){\zagc}
        \end{pspicture}
       \qquad  \rightsquigarrow \qquad
\psset{xunit=.8cm,yunit=.8cm}
\begin{pspicture}[.4](3,4)
  \pspolygon[linecolor=lightgray,fillstyle=solid,fillcolor=lightgray](0,0)(3,0)(3,3)(0,0)
  \pspolygon[linecolor=lightgray,fillstyle=solid,fillcolor=lightgray](0,0)(3,3)(0,3)(0,0)
  \psline[linewidth=1pt](0,0)(3,3)
  \psline[linestyle=dashed,linewidth=1pt](0,0)(3,0)
  \psline[linestyle=dashed,linewidth=1pt](0,3)(3,3)
  \psline[linestyle=dotted,linewidth=1pt](0,0)(0,3)
  \psline[linestyle=dotted,linewidth=1pt](3,0)(3,3)
  \psdots(0,0)(3,0)(3,3)(0,3)
  \psline[linewidth=1pt,linecolor=red]{->}(1.5,1.5)(1.2,1.8)
    \psline[linewidth=1pt,linecolor=red](1,2)(1.5,1.5)
  \psline[linewidth=1pt,linecolor=red]{->}(1.5,1.5)(1.8,1.2)
    \psline[linewidth=1pt,linecolor=red](2,1)(1.5,1.5)
  \psline[linewidth=1pt,linecolor=red]{->}(3,1.5)(2.4,1.2)
    \psline[linewidth=1pt,linecolor=red](2,1)(3,1.5)
  \psline[linewidth=1pt,linecolor=red]{->}(1.5,0)(1.8,.6)
    \psline[linewidth=1pt,linecolor=red](2,1)(1.5,0)
  \psline[linewidth=1pt,linecolor=red]{->}(1.5,3)(1.2,2.4)
    \psline[linewidth=1pt,linecolor=red](1,2)(1.5,3)
  \psline[linewidth=1pt,linecolor=red]{->}(0,1.5)(.6,1.8)
    \psline[linewidth=1pt,linecolor=red](1,2)(0,1.5)
  \psdots[linecolor=red](1,2)(2,1)(1.5,3)(1.5,0)(0,1.5)(3,1.5)(1.5,1.5)
  \psdots[dotstyle=pentagon](1.5,0)(1.5,3)
  \psdots[dotstyle=triangle](3,1.5)(0,1.5)
\end{pspicture}
\qquad  \rightsquigarrow \qquad
    \vcenter{ \xy
    (-2,8)*\bbrl{};
    (6,16)*\bbrl{};
    (0,16)*\bbid{};
    (12,16)*\bbid{};
    (2,8)*\bbrl{};
    (-2,-8)*\bbtriple{};
    (-6,0)*\bbainv{1};
    (-2,0)*\bbdl{};
    (2,0)*\bbdl{};
    (6,8)*\bbdual{};
    (2,24)*\bbdual{};
    (10,24)*\bbdual{};
    (10,8)*\bblr{};
    (14,8)*\bblr{};
    (6,16)*\bblr{};
    (8,0)*\bbid{};
    (12,0)*\bbid{};
    (16,0)*\bbid{};
    (12,-8)*\bbtriple{};
  \endxy}
\end{equation}
Here we have used the triangulation of the torus as a rectangle
where the dotted and dashed lines are identified in the usual way.
After the identifications this triangulation has a single interior
vertex and hence the single factor of $a^{-1}$ that appears in the
string diagram on the right.
\subsection{Independence of the triangulation of black boundaries}
\label{sec_indep-black}

We now define a morphism $Z(M)$ from the morphism $Z_{T_M}(M)$ which
does not depend on the choice of triangulation of the black
boundary. Observe that to each black boundary component $n_j$
triangulated with $h_j$ edges, we have associated the vector space
$A^{\otimes h_j}$.

\begin{prop}
\label{prop_blackboundary} For the triangulations $T_{I\times
I}^{k\ell}$ and $T_{S^1\times I}^{k\ell}$ of the flat strip
$I\times I$ and the cylinder $S^1\times I$ with $\ell$ incoming
edges and $k$ outgoing edges in their black boundaries, the state
sum of Definition~\ref{def_defstatesum} yields the morphisms
$P_{k\ell}\maps A^{\otimes\ell}\to A^{\otimes k}$ and
$Q_{k\ell}\maps A^{\otimes\ell}\to A^{\otimes k}$ of
Proposition~\ref{prop_PQ}. That is,
\begin{eqnarray}
  Z_{T_{I\times I}^{k\ell}}(I\times I) &=& P_{k\ell},\\
  Z_{T_{S^1\times I}^{k\ell}}(S^1\times I) &=& Q_{k\ell}.
\end{eqnarray}
\end{prop}

\begin{proof}
Write down the string diagram defining the state sum, \cf\
Figure~\ref{corner75}, and use the bubble move and the axioms of a
symmetric Frobenius algebra.

We here include the simplest triangulations of $S^1\times I$ and
$I\times I$ and the associated morphisms for $k=\ell=1$:
\begin{equation}
  Z_{T_{S^1\times I}^{11}} (
\psset{xunit=.5cm,yunit=.5cm}
\begin{pspicture}[.4](2,2.5)
  \rput(1,0){\identc}
\end{pspicture}
) \quad = \qquad
\psset{xunit=.5cm,yunit=.5cm}
\begin{pspicture}[.3](4,4)
  \pspolygon[linecolor=lightgray,fillstyle=solid,fillcolor=lightgray](0,0)(3,0)(3,3)(0,0)
  \pspolygon[linecolor=lightgray,fillstyle=solid,fillcolor=lightgray](0,0)(3,3)(0,3)(0,0)
  \psline[linewidth=.5pt](0,0)(3,3)
  \psline[linewidth=1.5pt](0,0)(3,0)
  \psline[linewidth=1.5pt](0,3)(3,3)
  \psline[linewidth=.5pt,linestyle=dashed](0,0)(0,3)
  \psline[linewidth=.5pt,linestyle=dashed](3,0)(3,3)
  \psdots[linewidth=.5pt](0,0)(3,0)(3,3)(0,3)
  \psline[linewidth=1pt,linecolor=red]{->}(1.5,1.5)(1.1,1.9)
    \psline[linewidth=1pt,linecolor=red](1,2)(1.5,1.5)
  \psline[linewidth=1pt,linecolor=red]{->}(1.5,1.5)(1.9,1.1)
    \psline[linewidth=1pt,linecolor=red](2,1)(1.5,1.5)
  \psline[linewidth=1pt,linecolor=red]{->}(3,1.5)(2.5,1.25)
    \psline[linewidth=1pt,linecolor=red](2,1)(3,1.5)
  \psline[linewidth=1pt,linecolor=red]{->}(2,1)(1.65,.3)
    \psline[linewidth=1pt,linecolor=red](2,1)(1.25,-.5)
  \psline[linewidth=1pt,linecolor=red]{->}(1.25,3.5)(1.05,2.4)
    \psline[linewidth=1pt,linecolor=red](1,2)(1.25,3.5)
  \psline[linewidth=1pt,linecolor=red]{->}(0,1.5)(.5,1.75)
    \psline[linewidth=1pt,linecolor=red](1,2)(0,1.5)
 \psdots[linecolor=red,dotsize=.15](3,1.5)(0,1.5)
  \psdots[linecolor=red](1,2)(2,1)(1.5,1.5)
  \psline[arrowsize=2.5pt 10, arrowinset=.8, linestyle=dashed]{->}(0,3)(0,.7)
  \psline[arrowsize=2.5pt 10, arrowinset=.8,linestyle=dashed]{->}(3,3)(3,.7)
\end{pspicture}
 \quad = \quad \vcenter{ \xy
    (0,0)*\bbcardy{};
    (-2,-8)*\bbainv{1};
\endxy}
\end{equation}
\begin{equation}
  Z_{T_{I\times I}^{11}} (
\psset{xunit=.5cm,yunit=.5cm}
\begin{pspicture}[.4](2,2.5)
  \rput(1,0){\identl}
\end{pspicture}
)\quad = \qquad
 \psset{xunit=.5cm,yunit=.5cm}
\begin{pspicture}[.3](4,4)
  \pspolygon[linecolor=lightgray,fillstyle=solid,fillcolor=lightgray](0,0)(3,0)(3,3)(0,0)
  \pspolygon[linecolor=lightgray,fillstyle=solid,fillcolor=lightgray](0,0)(3,3)(0,3)(0,0)
  \psline[linewidth=.5pt](0,0)(3,3)
  \psline[linewidth=1.5pt](0,0)(3,0)
  \psline[linewidth=1.5pt](0,3)(3,3)
  \psline[linewidth=.5pt](0,0)(0,3)
  \psline[linewidth=.5pt](3,0)(3,3)
  \psdots[linewidth=.5pt](0,0)(3,0)(3,3)(0,3)
  \psline[linewidth=1pt,linecolor=red]{->}(1.5,1.5)(1.1,1.9)
    \psline[linewidth=1pt,linecolor=red](1,2)(1.5,1.5)
  \psline[linewidth=1pt,linecolor=red]{->}(1.5,1.5)(1.9,1.1)
    \psline[linewidth=1pt,linecolor=red](2,1)(1.5,1.5)
  \psline[linewidth=1pt,linecolor=red]{->}(3.5,1.5)(2.7,1.25)
    \psline[linewidth=1pt,linecolor=red](2,1)(3.5,1.5)
  \psline[linewidth=1pt,linecolor=red]{->}(2,1)(1.65,.3)
    \psline[linewidth=1pt,linecolor=red](2,1)(1.25,-.5)
  \psline[linewidth=1pt,linecolor=red]{->}(1.25,3.5)(1.05,2.4)
    \psline[linewidth=1pt,linecolor=red](1,2)(1.25,3.5)
  \psline[linewidth=1pt,linecolor=red]{->}(-.5,1.5)(.3,1.75)
    \psline[linewidth=1pt,linecolor=red](1,2)(-.5,1.5)
 \psdots[linecolor=red,dotsize=.15](3.5,1.5)(-.5,1.5)
  \psdots[linecolor=red](1,2)(2,1)(1.5,1.5)
\end{pspicture}
 \quad = \quad \quad \ \; \vcenter{ \xy
    (2,0)*\bbid{};
    \endxy}
\end{equation}
\end{proof}

\parpic[l]{
$\xy
 (2,6)*{P_{44}};
 (2,-4)*{M};
 (-2,-4)*{};
\endxy$
$
\begin{aligned}
\includegraphics{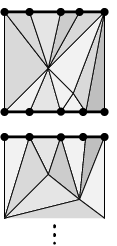}
\end{aligned}
$ } Given any triangulated open-closed cobordism $M$ with a black
boundary component homeomorphic to $I$ and triangulated with
$\ell$ edges, one can now glue a suitably triangulated cylinder
$I\times I$ to that boundary. By
Proposition~\ref{prop_interiorcoloured}, this yields the same
morphism $Z_{T_M}(M)$. Similarly, for every black boundary
component homeomorphic to $S^1$ and triangulated with $\ell$
edges, one can glue a suitably triangulated $S^1\times I$ to that
boundary, again leaving $Z_{T_M}(M)$ unchanged. It is therefore
sufficient to consider the restriction of $Z_{T_M}(M)$ to the
appropriate images of the idempotents $P_{\ell\ell}$ and
$Q_{\ell\ell}$, respectively. We therefore define:

\begin{defn}\label{def_tildeZ}
For every open-closed cobordism $M$ with triangulation $T_M$, we
define the state sum ${\tilde Z}_{T_M}(M)$ by subsequently pre-
and post-composing $Z_{T_M}(M)$ with the following morphisms: For
each $n_j \in \vec{n}=\del_0 M^{\mathrm{in}}$ triangulated with
$h_j$ edges, pre-composition with $\im P_{h_jh_j}$ if $n_j=1$ and
pre-composition with $\im Q_{h_jh_j}$ if $n_j=0$; For each $n'_j
\in \vec{n}'=\del_0 M^{\mathrm{out}}$ triangulated with $h_j$
edges, post-composition with $\coim P_{h_jh_j}$ if $n_j=1$ and
post-composition with $\coim Q_{h_jh_j}$ if $n_j=0$.

If we write $R^{(0)}_{k\ell}:=P_{k\ell}$ and
$R^{(1)}_{k\ell}:=Q_{k\ell}$, then the above composite is the
morphism
\begin{eqnarray}
\label{eq_tildez}
  {\tilde Z}_{T_M}(M) &=&
  \bigl(\bigotimes_{j=k+1}^{k+k^\prime}\coim R^{(n_j)}_{h_jh_j}\bigr)
  \circ Z_{T_M}(M)\circ
  \bigl(\bigotimes_{j=1}^k\im R^{(n_j)}_{h_jh_j}\bigr)\colon\nn\\
 &&
 \bigotimes_{j=1}^{k}R^{(n_j)}_{h_jh_j}(A^{\otimes h_j}) \to
    \bigotimes_{j=k+1}^{k+k'}R^{(n_j)}_{h_jh_j}(A^{\otimes h_j}).
\end{eqnarray}
\end{defn}

One can now use the isomorphisms of Corollary~\ref{cor_blackinv}
in order to relate the ${\tilde Z}_{T_M}(M)$ for different
triangulations of the black boundary as follows. The morphism
${\tilde Z}_{T_M}(M)$ is completely determined by the
triangulation of the boundary $\partial_0 M$ by
Proposition~\ref{prop_interiorcoloured}. Hence, the morphism
${\tilde Z}_{T_M}(M)$ associated to a triangulation $T_M$ is
related to the morphism ${\tilde Z}_{T'_M}(M)$ obtained from a
different triangulation ${T'_M}$ by gluing on cylinders whose
boundaries are appropriately triangulated. These cylinders yield
precisely the morphisms $P_{k\ell}$ and $Q_{k\ell}$.

\begin{defn}
\label{def_Z} For every open-closed cobordism $M$, we choose a
triangulation $T_M$. We define the state sum $Z(M)$ by
subsequently pre- and post-composing ${\tilde Z}_{T_M}(M)$ with
the following morphisms: For each $n_j \in \vec{n}=\del_0
M^{\mathrm{in}}$ triangulated with $h_j$ edges, pre-composition
with $\Phi_{h_j}$ if $n_j=1$ and pre-composition with $\Psi_{h_j}$
if $n_j=0$; For each $n'_j \in \vec{n}'=\del_0 M^{\mathrm{out}}$
triangulated with $h_j$ edges, post-composition with
$\Phi_{h_j}^{-1}$ if $n_j=1$ and post-composition with
$\Psi_{h_j}^{-1}$ if $n_j=0$. This yields the morphism
\begin{equation}
\label{eq_independent}
  Z(M) = \bigl(\bigotimes_{j=k+1}^{k+k^\prime}{(\Xi_{h_j}^{(n_j)})}^{-1}\bigr)
    \circ{\tilde Z}_{T_M}(M)\circ
    \bigl(\bigotimes_{j=1}^k\Xi_{h_j}^{(n_j)}\bigr)\colon
  A^{\otimes \vec{n}} \to A^{\otimes \vec{n'}},
\end{equation}
where we write $\Xi_{h_j}^{(0)}:=\Psi_{h_j}$ and
$\Xi_{h_j}^{(1)}:=\Phi_{h_j}$.
\end{defn}

The definition of $Z(M)$ is illustrated below:
\begin{equation}
 \vcenter{
\xy
 (0,0)*{\includegraphics{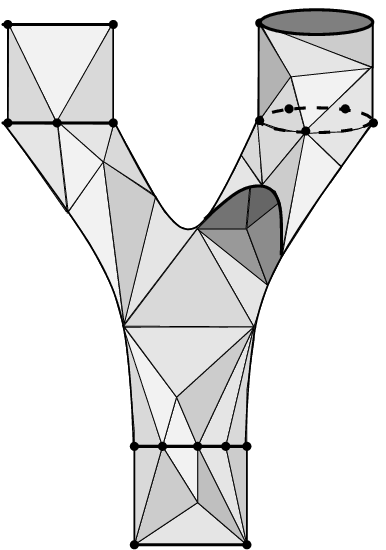}};
 \endxy}\quad
 \vcenter{
  \xy {\ar@{~>}^{Z} (-4,0);(4,0)}\endxy}\quad\vcenter{
 \xy
 (5,49)*+{A \otimes p(A)}="1";
 (5,37)*+{P_{22}\big(A^{\otimes 2}\big) \otimes Q_{55}\big( A^{\otimes 5}\big)}="2";
 (5,8)*+{P_{44}\big(A^{\otimes 4}\big)}="3";
 (5,-2)*+{A}="4";
  {\ar^{\Phi_2 \otimes \Psi_5} "1";"2"};
 {\ar^{{\tilde Z}_{T_M}(M)} "2";"3"};
 {\ar^{\Phi_4^{-1}}"3";"4"};
 \endxy}
\end{equation}

\begin{thm}
The morphism~\eqref{eq_independent} is well defined, \ie\ it does
not depend on the triangulation $T_M$ of $M$. In particular, it is
independent of the numbers $h_j$ of edges in
Definition~\ref{def_tildeZ} and Definition~\ref{def_Z}.
\end{thm}

\begin{proof}
Insert~\eqref{eq_tildez} into~\eqref{eq_independent} and draw the
cylinders over $I$ and over $S^1$ whose triangulations are given
by $\im P_{h_jh_j}\circ \Phi_{h_j}=P_{h_j1}$, \etc\ and glue them
to the triangulation used in the state sum $Z_{T_M}(M)$ of
Definition~\ref{def_defstatesum}. The invariance under bistellar
moves and elementary shellings of type (2.) of
Proposition~\ref{prop_interiorcoloured} then implies the theorem.
\end{proof}

\subsection{Open-closed Topological Quantum Field Theories}

From Definition~\ref{def_defstatesum}, it is obvious that the state
sum $Z(M)$ associates with the composition of open-closed cobordisms
the composition of morphisms of $\cal{C}$ and with the disjoint union
of open-closed cobordisms the tensor product of morphisms in
$\cal{C}$. It is not difficult to see that we get a symmetric monoidal
functor $Z\colon\cat{2Cob}^{\mathrm{ext}}\to\cal{C}$, \ie\ an
open-closed TQFT.

In this section, we show that this open-closed TQFT is the one
characterized by the knowledgeable Frobenius algebra of
Theorem~\ref{thm_kfrob}.

\subsubsection{Generators via the state sum construction}

Below we provide a choice of triangulation for some of the
generators in $\cat{2Cob}^{{\rm ext}}$.
\begin{gather}
\begin{aligned}
 \psset{xunit=.5cm,yunit=.5cm}
 \begin{pspicture}[.4](3,2.5)
    \rput(1,0){\multc}
 \end{pspicture}
 & \quad &
 \psset{xunit=.5cm,yunit=.5cm}
 \begin{pspicture}[.4](3,1.5)
    \rput(1,.4){\birthc}
 \end{pspicture}
  & \quad &
  \psset{xunit=.5cm,yunit=.5cm}
 \begin{pspicture}[.4](3,2.5)
        \rput(1,0){\multl}
 \end{pspicture}
 & \quad &
 \psset{xunit=.5cm,yunit=.5cm}
        \begin{pspicture}[.4](3,1.5)
        \rput(1,.4){\birthl}
        \end{pspicture}
  & \quad &
  \psset{xunit=.5cm,yunit=.5cm}
 \begin{pspicture}[.4](3,2.5)
        \rput(1,0){\ltc}
 \end{pspicture}
 \\
  \psset{xunit=.3cm,yunit=.3cm}
\begin{pspicture}[.4](7,6.7)
  \pspolygon[linecolor=lightgray,fillstyle=solid,fillcolor=lightgray]
        (2,0)(0,2)(.5,4.5)(3.5,5.3)(6.5,4.5)(7,2)(5,0)(2,0)
  \psline[linewidth=.5pt](2,0)(.5,4.5)
  \psline[linewidth=.5pt](2,0)(3.5,5.3)
  \psline[linewidth=.5pt](5,0)(3.5,5.3)
  \psline[linewidth=.5pt](5,0)(6.5,4.5)
  \psline[linewidth=1.5pt](2,0)(5,0)
  \psline[linewidth=1.5pt](0,2)(.5,4.5)
  \psline[linewidth=1.5pt](7,2)(6.5,4.5)
  \psline(2,0)(0,2)
  \psline(5,0)(7,2)
  \psline(.5,4.5)(3.5,5.3)
  \psline(6.5,4.5)(3.5,5.3)
  \psdots(2,0)(5,0)(3.5,5.3)(0,2)(7,2)(6.5,4.5)(.5,4.5)
  \psline[arrowsize=2.5pt 10, arrowinset=.8]{->}(5,0)(6,1)
  \psline[arrowsize=2.5pt 10, arrowinset=.8]{->}(5,0)(6.3,1.3)
  \psline[arrowsize=2.5pt 10, arrowinset=.8]{->}(3.5,5.3)(5,4.9)
  \psline[arrowsize=2.5pt 10, arrowinset=.8]{->}(3.5,5.3)(5.45,4.75)
  \psline[arrowsize=2.5pt 10, arrowinset=.8]{->}(0,2)(1,1)
  \psline[arrowsize=2.5pt 10, arrowinset=.8]{->}(.5,4.5)(2,4.9)
  \rput(4,6.2){$\scs a^{-1}$}
   \rput(6.5,0){$\scs a^{-1}$}
\end{pspicture}
& \quad &
 \psset{xunit=.35cm,yunit=.35cm}
 \begin{pspicture}[.4](5,4)
  \pspolygon[fillstyle=solid,fillcolor=lightgray](0,0)(4,0)(2,2.2)(0,0)
  \psdots(0,0)(4,0)(2,2.2)
  \psline[linewidth=1.5pt](0,0)(4,0)
  \psline[arrowsize=2.5pt 10, arrowinset=.8]{->}(2,2.2)(1,1.1)
  \psline[arrowsize=2.5pt 10, arrowinset=.8]{->}(2,2.2)(3,1.1)
  \rput(4.9,.1){$\scs  a^{-1}$}
  \rput(2.4,2.9){$\scs a^{-1}$}
    \end{pspicture}
  & \quad &
 \psset{xunit=.25cm,yunit=.25cm}
\begin{pspicture}[.4](7.5,6)
  \pspolygon[linecolor=lightgray,fillstyle=solid,fillcolor=lightgray]
        (2,0)(0,2.5)(2,5)(5,5)(7,2.5)(5,0)(2,0)
  \psline[linewidth=1.5pt](2,0)(5,0)
  \psline[linewidth=1.5pt](0,2.5)(2,5)
  \psline[linewidth=1.5pt](7,2.5)(5,5)
   \psline(2,5)(5,5)
   \psline(0,2.5)(2,0)
   \psline(5,0)(7,2.5)
   \psline[linewidth=.5pt](2,0)(2,5)
   \psline[linewidth=.5pt](5,0)(5,5)
   \psline[linewidth=.5pt](2,0)(5,5)
  \psdots(2,0)(0,2.5)(2,5)(5,5)(7,2.5)(5,0)
\end{pspicture}
  & \quad &
  \psset{xunit=.35cm,yunit=.35cm}
  \begin{pspicture}[.4](5,4)
  \pspolygon[fillstyle=solid,fillcolor=lightgray](0,0)(4,0)(2,2.2)(0,0)
  \psdots(0,0)(4,0)(2,2.2)
  \psline[linewidth=1.5pt](0,0)(4,0)
    \end{pspicture}
   & \quad &
   \psset{xunit=.25cm,yunit=.25cm}
    \begin{pspicture}[.4](7,6)
    \pspolygon[linecolor=lightgray,fillstyle=solid,fillcolor=lightgray]
        (3,0)(0,0)(0,4)(6,4)(3,0)
    \psline[linewidth=1.5pt](3,0)(0,0)
    \psline[linewidth=1.5pt](0,4)(3,4)
    \psline(0,0)(0,4)
    \psline(3,0)(6,4)
    \psline(3,4)(6,4)
    \psline[linewidth=.5pt](0,0)(3,4)
      \psline[linewidth=.5pt](0,0)(6,4)
    \psline[arrowsize=2.5pt 10, arrowinset=.8]{->}(0,0)(0,2.2)
    \psline[arrowsize=2.5pt 10, arrowinset=.8]{->}(3,0)(4.5,2)
    \rput(4.5,.1){$\scs a^{-1}$}
    \psdots(3,0)(0,0)(0,4)(6,4)(3,4)
    \end{pspicture}
\end{aligned}
\end{gather}
Those edges with matching arrow heads on the triangulations are to be
identified. The black boundaries are depicted slightly thicker than
the coloured boundaries. A choice of triangulation for the remaining
generators is immediate from those above. The factors of $a^{-1}$ are
meant to remind the reader which vertices in the triangulation
contribute factors of $a^{-1}$.

Using these triangulations we can compute the morphisms
$Z_{T_M}(M)$ associated to the open-closed cobordisms $M$
generating $\cat{2Cob}^{{\rm ext}}$. For completeness, we include
the triangulation of the cylinders $S^1\times I$ and $I\times I$
as well.

\begin{gather}
\label{eq_trigenstart}
\begin{aligned}
\psset{xunit=.3cm,yunit=.3cm}
Z_{T_M}\bigl(
\begin{pspicture}[.4](3.5,3)
  \rput(2,0){\multc}
\end{pspicture}
\bigr) = \;\;
\vcenter{
  \xy 0;/r.22pc/:
    (0,0)*\bbmedmu{};
    (-2,-8)*\bbainv{2};
    (-4,8)*\bbcardy{};
    (4,8)*\bbcardy{};
  \endxy}
  \;\; = \;\;
\vcenter{
  \xy 0;/r.22pc/:
    (0,0)*\bbmedmu{};
    (-4,6)*\bbP{};
    (4,6)*\bbP{};
  \endxy}
&\qquad&
\psset{xunit=.3cm,yunit=.3cm}
Z_{T_M}\bigl(
\begin{pspicture}[.4](1.8,1.5)
  \rput(1,1){\birthc}
\end{pspicture}
\bigr)= \;\;
\vcenter{
  \xy
    (0,0)*\bbeta{};
    (-2,-8)*\bbainv{1};
  \endxy}
\end{aligned}\\
\begin{aligned}
\psset{xunit=.3cm,yunit=.3cm}
Z_{T_M}\bigl(
\begin{pspicture}[.4](3.5,2.5)
  \rput(2,0){\comultc}
\end{pspicture}
\bigr) = \;\;
\vcenter{
  \xy 0;/r.22pc/:
    (0,16)*\bbmeddelta{};
    (-6,-8)*\bbainv{3};
    (4,-8)*\bbid{};
    (-4,0)*\bbcardy{};
    (4,0)*\bbcardy{};
  \endxy}
  \;\; = \;\;
\vcenter{
  \xy 0;/r.22pc/:
    (0,14)*\bbmeddelta{};
    (-6,-6)*\bbainv{1};
    (4,-6)*\bbid{};
    (-4,0)*\bbP{};
    (4,0)*\bbP{};
  \endxy}
&\qquad&
\psset{xunit=.3cm,yunit=.3cm}
Z_{T_M}\bigl(
\begin{pspicture}[.4](2,1.5)
  \rput(1,0){\deathc}
\end{pspicture}
\bigr) = \;\;
\vcenter{
  \xy
    (-2,0)*\bbepsilon{};
  \endxy}
\end{aligned}
\end{gather}

\begin{gather}
\begin{aligned}
\psset{xunit=.3cm,yunit=.3cm}
Z_{T_M}\bigl(
\begin{pspicture}[.4](4,3)
  \rput(2,0){\multl}
\end{pspicture}
\bigr) = \;\;
\vcenter{
  \xy
    (0,0)*\bbmu{};
  \endxy}
\quad
\psset{xunit=.3cm,yunit=.3cm}
Z_{T_M}\bigl(
\begin{pspicture}[.4](2,1.5)
  \rput(1,1){\birthl}
\end{pspicture}
\bigr) = \;\;
  \xy
    (0,6)*\bbeta{};
  \endxy
\end{aligned}\\
\begin{aligned}
\psset{xunit=.3cm,yunit=.3cm}
Z_{T_M}\bigl(
\begin{pspicture}[.4](4,3)
  \rput(2,0){\comultl}
\end{pspicture}
\bigr) = \;\;
\vcenter{
  \xy
    (0,0)*\bbdelta{};
  \endxy}
\quad
\psset{xunit=.3cm,yunit=.3cm}
Z_{T_M}\bigl(
\begin{pspicture}[.4](2,1.5)
  \rput(1,1){\deathl}
\end{pspicture}
\bigr) = \;\;
  \xy
    (0,6)*\bbepsilon{};
  \endxy
\end{aligned}
\end{gather}
\begin{gather}
\begin{aligned}
\psset{xunit=.3cm,yunit=.3cm}
Z_{T_M}\bigl(
\begin{pspicture}[.4](2,2.5)
  \rput(1,0){\ctl}
\end{pspicture}
\bigr) = \;\;
\vcenter{
  \xy
    (0,0)*\bbcardy{};
  \endxy}
  \;\; = \;\;
\vcenter{
  \xy
    (0,-2)*\bbP{};
    (-2,-8)*\bba{};
  \endxy}
\quad
\psset{xunit=.3cm,yunit=.3cm}
Z_{T_M}\bigl(
\begin{pspicture}[.4](2,2.5)
  \rput(1.2,0){\ltc}
\end{pspicture}
\bigr) = \;\;
\vcenter{
  \xy
    (0,0)*\bbcardy{};
    (-2,-8)*\bbainv{1};
  \endxy}
  \;\; = \;\;
  \vcenter{ \xy
    (0,-2)*\bbP{};
  \endxy}
\end{aligned}
\end{gather}
\begin{gather}
\begin{aligned}
\psset{xunit=.3cm,yunit=.3cm}
Z_{T_M}\bigl(
\begin{pspicture}[.4](2,2.5)
  \rput(1.2,0){\identc}
\end{pspicture}
\bigr)= \;\;
\vcenter{
  \xy
    (0,0)*\bbcardy{};
    (-2,-8)*\bbainv{1};
  \endxy}
  \;\; = \;\;
\vcenter{
  \xy
    (0,-2)*\bbP{};
  \endxy}
\qquad
\psset{xunit=.3cm,yunit=.3cm}
Z_{T_M}\bigl(
\begin{pspicture}[.4](2,2.5)
  \rput(1,0){\identl}
\end{pspicture}
\bigr) = \;\;
\vcenter{
  \xy
    (2,0)*\bbid{};
  \endxy}
\end{aligned}
\label{eq_trigenend}
\end{gather}

\begin{thm}
\label{thm_statesum} Let $\cal{C}$ be an Abelian symmetric
monoidal category and $A$ be a rigid and strongly separable
algebra object in $\cal{C}$ that is equipped with the structure of
a symmetric Frobenius algebra. Then the state
sum~\eqref{def_defstatesum} defines an open-closed TQFT
$Z\colon\cat{2Cob}^\mathrm{ext}\to\cal{C}$. It is characterized by
the knowledgeable Frobenius algebra constructed from $A$ in
Theorem~\ref{thm_kfrob}.
\end{thm}

\begin{proof}
Using the triangulations of the generators given in
\eqref{eq_trigenstart}-\eqref{eq_trigenend}, compute the morphisms
$Z_{T_M}(M)$ for each generator of $\twocob$. Pre and post
composing with the relevant maps specified in
Definitions~\ref{def_tildeZ} and \ref{def_Z} produces the
knowledgeable Frobenius algebra $(A,Z(A),\imath,\imath^{\ast})$
defined in Theorem~\ref{thm_kfrob}. For example, $Z_{T_M}\left(
\psset{xunit=.16cm,yunit=.16cm} \begin{pspicture}[.4](4,1.5)
\rput(2.5,0){\multc} \end{pspicture}\right)=\mu_A \circ (p \otimes
p)$ so that ${\tilde Z}_{T_M}\left(
\psset{xunit=.16cm,yunit=.16cm}
\begin{pspicture}[.4](4,1.5) \rput(2.5,0){\multc}
\end{pspicture}\right)=\coim Q_{11} \circ \mu_A \circ(p \otimes
p)\circ(\im Q_{11} \otimes \im Q_{11})$. Noting that $Q_{11}=p$
and using the image factorization of $p$ \eqref{eq_factorization}
together with the idempotent property $p^2=p$ it is easy to check
that
\begin{equation} Z\left( \psset{xunit=.16cm,yunit=.16cm}
\begin{pspicture}[.4](4,1.5) \rput(2.5,0){\multc}
\end{pspicture}\right)=\coim p \circ \mu_A \circ (\im p \otimes \im p)
\end{equation} as specified in Theorem~\ref{thm_kfrob}.

Since $\cat{2Cob}^{\mathrm{ext}}$ is the strict symmetric monoidal
category freely generated by a knowledgeable Frobenius algebra object,
this uniquely determines a symmetric monoidal functor
$Z\colon\cat{2Cob}^{\mathrm{ext}}\to\cal{C}$.
\end{proof}

Recall that given an open-closed TQFT, the algebra object $A:=Z(I)$
does not necessarily determine the object $C:=Z(S^1)$.  Consider, for
example, the knowledgeable Frobenius algebra
$(A,C,\imath,\imath^\ast)$ of Example~\ref{ex_notcentre} in which
$C\not\cong Z(A)$, and secondly the knowledgeable Frobenius algebra
$(A,Z(A),\imath^\prime,\imath^{\prime\ast})$ constructed in
Theorem~\ref{thm_kfrob} based on the same $A$. Both characterize an
open-closed TQFT, but only the latter one can be obtained from the
state sum.

Conversely, in an open-closed TQFT, the object $Z(S^1)$ does not
determine the object $Z(I)$. This can be easily seen from
Example~\ref{ex_matrixexample} below.

\subsection{Examples}

In \cite{LP05} it was shown that connected open-closed cobordisms are
determined up to orientation-preserving diffeomorphism preserving the
black boundary by a set of topological invariants defined in the work
of Baas, Cohen, and Ram\'{\i}rez \cite{BCR}. These topological
invariants are the \emph{genus} (defined as the genus of the
underlying topological 2-manifold), the \emph{window number}, defined
as the number of components of $\partial_1 M$ diffeomorphic to $S^1$,
and the \emph{boundary permutation}. For a surface $M$ ($\partial_0 M
= \emptyset$) only the genus and window number are relevant.  In this
context we will refer to the window number as the number of punctures
in $M$.

Let $(A,C,\imath,\imath^\ast)$ be a knowledgeable Frobenius
algebra in a symmetric monoidal category $\cal{C}$. We call
$\mu_C\circ\Delta_C\colon C\to C$ the \emph{genus-one operator}
and $\imath^\ast\circ\imath\colon C\to C$ the \emph{window
operator}. The invariant associated to the connected surface
$M^\ell_k$ of genus $\ell$ with $k$ punctures is determined by
evaluating the morphism
\begin{equation}
  Z(M^\ell_k)=\varepsilon_C\circ\big( \imath^\ast\circ\imath \big)^{k}\circ \big( \mu_C\circ\Delta_C
  \big)^{\ell}\circ\eta_C\colon\1\to\1
\end{equation}
in $\cal{C}$.

In this section, we provide several examples of strongly separable
symmetric Frobenius algebras and use the genus-one operator and
the window operator to compute the state sum invariant
$Z\big(M^{\ell}_{k}\big)$.

\begin{myexample}
\label{ex_matrixexample}
Let $k$ be a field, $n\in\N$, and $m_1,\ldots,m_n\in\N$, and consider
the direct product\footnote{We write $\oplus$ because this is actually
the biproduct in the Abelian category $\cat{Vect}_k$.}
\begin{equation}
  A:=\bigoplus_{j=1}^n M_{m_j}(k)
\end{equation}
of matrix algebras. We choose a basis ${\{e_{pq}^{(j)}\}}_{1\leq
p,q\leq m_j, 1\leq j\leq n}$ of $A$ such that the multiplication reads
$\mu_A(e_{pq}^{(j)}\otimes
e_{rs}^{(\ell)})=\delta_{j\ell}\delta_{rq}e_{ps}^{(j)}$ with unit
$\eta_A(1)=\sum_{j=1}^n\sum_{p=1}^{m_j}e_{pp}^{(j)}$. The $k$-algebra
$(A,\mu_A,\eta_A)$ is strongly separable if and only if for all $j$,
$\chr k$ does not divide $m_j$. From now on we assume that this
condition holds.

The centre $Z(A)$ of $A$ has a basis ${\{z_j\}}_{1\leq j\leq n}$
of orthogonal idempotents $z_j:=\sum_{p=1}^{m_j}e_{pp}^{(j)}$,
\ie\ $\mu_A(z_j\otimes z_\ell)=\delta_{j\ell}z_j$. The symmetric
Frobenius algebra structures
$(A,\mu_A,\eta_A,\Delta_A,\epsilon_A)$ are characterized by the
invertible central elements $a=\sum_{j=1}^na_jz_j$, \ie\ $a_j\in
k\backslash\{0\}$ for all $j$, as follows:
\begin{eqnarray}
  \Delta_A(e_{pq}^{(j)}) &=& a_j m_j^{-1}\sum_{r=1}^{m_j}e_{pr}^{(j)}\otimes e_{rq}^{(j)},\\
  \epsilon_A(e_{pq}^{(j)}) &=& \delta_{pq} m_j a_j^{-1},
\end{eqnarray}
and indeed one finds $(\mu_A\circ\Delta_A\circ\eta_A)(1)=a$ for
the window element. This illustrates further the distinction
between special Frobenius algebras and strongly separable
Frobenius algebras. $A$ is special if and only if $a_i=a_j$ for
all $i,j$. We compute the idempotent $p$ of~\eqref{eq_idempotent}
as follows:
\begin{equation}
  p(e_{pq}^{(j)}) = \delta_{pq}m_j^{-1}\sum_{r=1}^{m_j}e_{rr}^{(j)},
\end{equation}
and indeed the image is $p(A)\cong Z(A)$ with the splitting
\begin{alignat}{2}
  \im p   &\colon p(A)\to A,\quad& z_j&\mapsto \sum_{p=1}^{m_j}e_{pp}^{(j)},\\
  \coim p &\colon A\to p(A),\quad& e_{pq}^{(j)}&\mapsto \delta_{pq}m_j^{-1}z_j.
\end{alignat}
The knowledgeable Frobenius algebra $(A,C,\imath,\imath^\ast)$ of
Theorem~\ref{thm_kfrob} for this algebra $A$ is given by the
following commutative Frobenius algebra structure
$(C,\mu_C,\eta_C,\Delta_C,\epsilon_C)$ on $C:=Z(A)$:
\begin{eqnarray}
  \mu_C(z_j\otimes z_\ell) &=& \delta_{j\ell}z_j,\\
  \eta_C(1) &=& \sum_{j=1}^n z_j,\\
  \Delta_C(z_j) &=& a_j^2m_j^{-2} z_j\otimes z_j,\\
  \epsilon_C(z_j) &=& m_j^2a_j^{-2},
\end{eqnarray}
together with
\begin{alignat}{5}
  \imath&\colon C&\to& A,\quad& z_j\mapsto&\sum_{p=1}^{m_j}e_{pp}^{(j)},\\
  \imath^\ast&\colon A&\to& C,\quad& e_{pq}^{(j)}\mapsto& a_jm_j^{-1}\delta_{pq}z_j.
\end{alignat}

We finally compute the genus-one operator
$(\mu_C\circ\Delta_C)(z_j)=a_j^2m_j^{-2}z_j$ and the window
operator $(\imath^\ast\circ\imath)(z_j)=a_jz_j$, and so the
invariant~\eqref{eq_surface} associated with the genus
$\ell$-surface with $k$ punctures, $k,\ell\in\N_0$, is
\begin{equation}
\label{eq_surface}
  Z(M^\ell_k)(1)=(\epsilon_C\circ{(\imath^\ast\circ\imath)}^k\circ{(\mu_C\circ\Delta_C)}^\ell\circ\eta_C)(1)
  = \sum_{j=1}^na_j^{k+2(\ell-1)}m_j^{-2(\ell-1)}.
\end{equation}
\end{myexample}

Fukuma--Hosono--Kawai~\cite{FHK} choose the canonical Frobenius
algebra structure on $A$, \ie\ $a=\eta$ and therefore $a_j=1$ for
all $j$. In this case, the invariant is blind to the \emph{window
number} $k$. With a generic symmetric Frobenius algebra structure,
however, one can easily obtain an invariant that can distinguish
any two inequivalent connected surfaces.

\begin{myexample}
Let $G$ be a finite group, $k$ a field, and $A:=k[G]$ be the group
algebra. We choose the basis ${\{g\}}_{g\in G}$ for $A$ and have
$\mu_A(g\otimes h)=gh$ for $g,h\in G$ and $\eta_A(1)=e$. The
$k$-algebra $(A,\mu_A,\eta_A)$ is strongly separable if and only if
$\chr k$ does not divide the order $|G|$ of $G$. We now assume that
this condition holds.

We denote by $[g]:=\{\,hgh^{-1}\colon\, h\in G\,\}\subseteq G$ the
conjugacy class of $g\in G$ and by $G/\sim:=\{\,[g]\colon\, g\in G\}$
the set of classes. Then the centre $Z(A)$ has the basis
${\{z_{[g]}\}}_{[g]\in G/\sim}$ where $z_{[g]}:=\sum_{h\in [g]}h$
denotes the class sum. We have the unit $\eta_A(1)=\sum_{[g]\in
G/\sim}z_{[g]}$ and $\mu_A(z_{[g]}\otimes z_{[h]})=\sum_{[\ell]\in
G/\sim}\mu_{[g],[h]}^{[\ell]}z_{[\ell]}$ for all $g,h\in G$ with some
$\mu_{[g],[h]}^{[\ell]}\in k$.

The $z_{[g]}$ are in general not orthogonal idempotents. Working with
a generic invertible central element in the basis
${\{z_{[g]}\}}_{[g]\in G/\sim}$ is not very instructive. If $k$ is
algebraically closed, the irreducible characters $\chi_\rho\colon G\to
k$ provide us with a basis ${\{z_\rho\}}_\rho$ of orthogonal
idempotents $z_\rho:=d_\rho{|G|}^{-1}\sum_{g\in G}\chi_\rho(g)g$,
$d_\rho=\chi_\rho(e)$, for $Z(A)$. We then get the same results as for
a direct product of $d_\rho\times d_\rho$-matrix algebras.

In the following, we restrict ourselves to the symmetric Frobenius
algebra structure
\begin{eqnarray}
  \Delta_A(g) &=& \sum_{h\in G} h\otimes h^{-1}g,\\
  \epsilon_A(g) &=& \begin{cases}
    1,&\text{if $g=e$}\\
    0,&\text{else}
  \end{cases}
\end{eqnarray}
which is characterized by the window element
$(\mu_A\circ\Delta_A\circ\eta_A)(1)=|G|e=|G|\eta_A(1)$. The symmetric
Frobenius algebra $(A,\mu_A,\eta_A,\Delta_A,\epsilon_A)$ is therefore
special in the sense of~\eqref{eq_special}. In this case
\begin{eqnarray}
  g^\ast(1) = (\Delta_A\circ\mu_A)(1) &=& \sum_{h\in G}h\otimes h^{-1},\nn\\
  g^{(3)}((g\otimes h)\otimes\ell)
    = (\epsilon_A\circ\mu_A\circ(\mu_A\otimes\id_A))((g\otimes h)\otimes\ell)
    &=&\begin{cases}
      1,&\text{if $gh\ell=e$}\nn\\
      0,&\text{else}
    \end{cases}
\end{eqnarray}
The state sum $Z(M)$ then agrees with the partition function of a
topological gauge theory with gauge group $G$ or, in other words, with
the volume of the moduli space of flat $G$-bundles on $M$. In the
state sum of Definition~\ref{def_defstatesum}, the window element
$|G|$ is divided out for every vertex in the interior of $M$ (this
prefactor of $Z(M)$ is sometimes called the \emph{anomaly}). In the
closed TQFT, the meaning of this factor is somewhat mysterious --- the
factor is merely needed in order to make the 1-3 Pachner move work ---
but in our extension to the open-closed TQFT, the factor $|G|$ is
directly related to the symmetric Frobenius algebra structure of $A$
and thereby to topology.
\end{myexample}

\begin{rem}
Although our state sum of Definition~\ref{def_defstatesum} requires an
oriented $2$-manifold, the previous example with the group algebra
$A=k[G]$ makes sense even for unoriented manifolds (without
boundary). This is possible because $A$ also has the structure of an
involutory Hopf algebra
$(A,\mu_A,\eta_A,\Delta_A^{\mathrm{Hopf}},\epsilon_A^{\mathrm{Hopf}},S_A)$
with
\begin{eqnarray}
  \Delta_A^{\mathrm{Hopf}}(g) &=& g\otimes g,\\
  \epsilon_A^{\mathrm{Hopf}}(g) &=& 1,\\
  S_A(g) &=& g^{-1},
\end{eqnarray}
with a co-integral $\sum_{g\in G}g$ and an integral
$g\mapsto\delta_G(g)$ where $\delta_G(e)=1$ and $\delta_G(g)=0$ for
all $g\neq e$. For this involutory Hopf algebra, one can evaluate
Kuperberg's $3$-manifold invariant~\cite{Kuperberg} which does not
refer to the $3$-simplices and therefore makes sense for (unoriented)
$2$-manifolds, too. In the oriented case, it agrees with our state
sum. The unoriented case is treated in more generality in~\cite{AN}.
\end{rem}

%
\section{State sums with D-branes}
%

Our next example, the groupoid algebra of a finite groupoid, also
yields the state sum of an open-closed TQFT in a straightforward way,
but in addition it provides us with an example of an $S$-coloured
open-closed TQFT, \cf\ Section~5 of~\cite{LP05}.

A \emph{groupoid} $\cal{G}=(X,G,s,t,\imath,\circ,{}^{-1})$ consists of
sets $X$ (\emph{objects}) and $G$ (\emph{morphisms}) and maps $s\colon
G\to X$ (\emph{source}), $t\colon G\to X$ (\emph{target}),
$\imath\colon X\to G$ (\emph{identity}), $\circ\colon
G{}_t\!\!\times_s G:=\{\,(h_1,h_2)\in G\times
G\colon\,t(h_1)=s(h_2)\,\}\to G$ (\emph{composition}, written from
left to right) and ${}^{-1}\colon G\to G$ (\emph{inversion}) such that
the following conditions are satisfied,
\begin{enumerate}
\item
  $s(\imath(x))=x$ and $t(\imath(x))=x$ for all $x\in X$,
\item
  $s(h_1\circ h_2)=s(h_1)$ and $t(h_1\circ h_2)=t(h_2)$ for all
  $(h_1,h_2)\in X{}_t\!\!\times_s X$,
\item
  $(h_1\circ h_2)\circ h_3 = h_1\circ(h_2\circ h_3)$ for all
  $h_1,h_2,h_3\in G$ for which $t(h_1)=s(h_2)$ and $t(h_2)=s(h_3)$,
\item
  $\imath(s(h))\circ h=h=h\circ\imath(t(h))$ for all $h\in G$,
\item
  $s(h^{-1})=t(h)$ and $t(h^{-1})=s(h)$ for all $h\in G$,
\item
  $h^{-1}\circ h=\imath(t(h))$ and $h\circ h^{-1}=\imath(s(h))$ for
  all $h\in G$.
\end{enumerate}
The groupoid is called \emph{finite} if $G$ is a finite set. For
every $x\in X$, we denote its connected component by
$[x]:=\{\,t(h)\colon\,h\in G, s(h)=x\,\}$. The groupoid is called
\emph{connected} if $X=[x]$ for some $x\in X$. For $x\in X$, the
\emph{star of $\cal{G}$ at $x$} is the set,
\begin{eqnarray}
  \st_{\cal{G}}(x) = \{\,g\in G\colon\,s(g)=x\,\}.
\end{eqnarray}
We denote the order of the star of $\cal{G}$ at $x\in X$ by
$N_{[x]}:=|\st_{\cal{G}}(x)|$. It depends only on the connected
component $[x]$ of $x\in X$.

Given a finite groupoid $\cal{G}=(X,G,s,t,\imath,\circ,{}^{-1})$ and
a field $k$, the \emph{groupoid algebra} $(k[G],\mu,\eta)$ is the free
vector space $k[G]$ on the set of morphisms with the operations,
\begin{eqnarray}
  \mu(h_1\otimes h_2) &=& \begin{cases}
    h_1\circ h_2,&\text{if $t(h_1)=s(h_2)$}\\
    0,           &\text{else}
  \end{cases}\\
  \eta(1) &=& \sum_{x\in X}\imath(x),
\end{eqnarray}
where $h_1,h_2\in G$.

\begin{myexample}
Let $\cal(G)=(X,G,s,t,\imath,\circ,{}^{-1})$ be a finite groupoid and
consider the groupoid algebra $A:=k[G]$. The $k$-algebra $A$ is
strongly separable if and only of $\chr k$ does not divide $N_{[x]}$
for any $x\in X$. From now on, we assume that this is the case.

We denote by $G^{(0)}:=\{\,g\in G\colon\,s(g)=t(g)\,\}\subseteq G$ the
set of automorphisms, by $[g]:=\{\,h\circ g\circ h^{-1}\colon\,h\in G,
t(h)=t(g)\,\}$ the conjugacy class of the automorphism $g\in G^{(0)}$,
and by $G^{(0)}/\sim:=\{\,[g]\colon\,g\in G^{(0)}\,\}$ the set of
conjugacy classes. Choose the basis ${\{h\}}_{h\in G}$ of $A$. We find
the centre $Z(A)\cong k[G^{(0)}/\sim]$ with a basis
${\{z_{[g]}\}}_{g\in G^{(0)}/\sim}$ where $z_{[g]}:=\sum_{h\in[g]}h$
denotes the class sum.

The canonical symmetric Frobenius algebra structure
$(A,\mu_A,\eta_A,\Delta_A,\epsilon_A)$ is given by
\begin{eqnarray}
  \epsilon_A(g) &=& \begin{cases}
    N_{[s(g)]},&\text{if $g=\imath(s(g))$}\\
    0,&\text{else}
  \end{cases}\\
  \Delta_A(g) &=& \frac{1}{N_{[t(g)]}}\sum_{h\in G\colon\, s(h)=s(g)} h\otimes(h^{-1}\circ g),
\end{eqnarray}
from which we obtain the canonical idempotent~\eqref{eq_idempotent}
\begin{equation}
  p(g) = \begin{cases}
    z_{[g]}/N_{[t(g)]},&\text{if $t(g)=s(g)$}\\
    0,&\text{else}
  \end{cases}
\end{equation}
with the image decomposition
\begin{alignat}{5}
  \im p&\colon Z(A)&\to& A,\quad&z_{[g]}\mapsto&\sum_{h\in[g]}h,\\
  \coim p&\colon A&\to& Z(A),\quad&g\mapsto&\begin{cases}
    z_{[g]}/N_{[t(g)]},&\text{if $s(g)=t(g)$}\\
    0,&\text{else.}
  \end{cases}
\end{alignat}
From these data, one can compute the knowledgeable Frobenius algebra
$(A,Z(A),\imath,\imath^\ast)$ that appears in Theorem~\ref{thm_kfrob}
with $\imath=\im p$ and $\imath^\ast=\coim p$. The state sum
construction therefore yields the corresponding open-closed TQFT.
\end{myexample}

There is, however, another point of view according to which the
groupoid algebra gives rise to an $X$-coloured knowledgeable
Frobenius algebra (Section~5 of~\cite{LP05}). Although this
example is rather trivial, it nicely illustrates where the various
structures appear.

\begin{myexample}
Let $\cal{G}=(X,G,s,t,\imath,\circ,{}^{-1})$ be a finite groupoid and
$k$ be a field such that $\chr k$ does not divide $N_{[x]}$ for any
$x\in X$. Denote by $\Hom(x,y)=\{\,g\in G\colon s(g)=x, t(g)=y\,\}$
the morphisms from $x$ to $y$. Then there is a family of vector spaces
$A_{xy}:=k[\Hom(x,y)]$. By restricting the operations of the groupoid
algebra $A=k[G]$ to the $A_{xy}$, we obtain the following linear maps:
\begin{eqnarray}
  \mu_{xyz}\colon A_{xy}\otimes A_{yz}\to A_{xz},\quad &g_1\otimes g_2&\mapsto g_1\circ g_2,\\
  \eta_x(1)\colon k\to A_{xx},\quad &1&\mapsto\imath(x),\\
  \Delta_{xyz}\colon A_{xz}\to A_{xy}\otimes A_{yz},\quad&g&
    \mapsto\frac{1}{N_{[t(g)]}}\sum_{h\in G\colon\,s(h)=x} h\otimes h^{-1}\circ g,\\
  \epsilon_x\colon A_{xx}\to k,\quad &g&\mapsto
    \begin{cases}
      N_{[s(g)]},&\text{if $g=\imath(x)$}\\
      0,&\text{else}
    \end{cases}
\end{eqnarray}
for $x,y,z\in X$. Similarly by restricting $\imath$ and $\imath^\ast$,
we find for all $x\in X$:
\begin{alignat}{5}
  \imath_x&\colon Z(A)&\to& A_{xx}, z_{[g]}\mapsto&\sum_{h\in[g]\colon\,h\in\Hom(x,x)}h,\\
  \imath^\ast&\colon A_{xx}&\to& Z(A), g\mapsto&\frac{1}{N_{[x]}}z_{[g]}.
\end{alignat}
Then we have an $X$-coloured knowledgeable Frob\-enius algebra
\begin{equation}
 (\{A_{xy}\},\{\mu_{xyz}\},\{\eta_x\},\{\Delta_{xyz}\},\{\epsilon_x\},
Z(A),\{\imath_x\},\{\imath^\ast_x\}).
\end{equation}
The commutative Frobenius algebra structure of $Z(A)$ is as in
Theorem~\ref{thm_kfrob}. In particular, each $A_{xx}$, $x\in X$,
forms a symmetric Frobenius algebra, the $\imath_x\colon Z(A)\to
A_{xx}$ are algebra homomorphisms, and each $A_{xy}$ forms an
$(A_{xx},A_{yy})$-bimodule with dual $A_{yx}$. Observe that the
state sum can be evaluated directly for the full groupoid algebra
\begin{equation}
  A = \bigoplus_{x,y\in X} A_{xy},
\end{equation}
and so the vector space associated with the unit interval is
precisely this direct sum. If one restricts it to the subspaces
$A_{xy}$ corresponding to the boundary colours $x,y\in X$ of a
given interval, one obtains an $X$-coloured open-closed TQFT. The
full state sum with $A$, however, contains more than just these
homogeneous elements. It includes their linear combinations as
well.
\end{myexample}

This last example is especially relevant in the context where the
open-closed cobordisms are interpreted as open and closed string
worldsheets. In this case, the colours of an $X$-coloured
knowledgeable Frobenius algebra are interpreted as the set of
boundary conditions, or D-branes, for the open strings.  The
decomposition of the finite groupoid algebra then allows the state
sum to compute topological invariants of open and closed string
worldsheets equipped with D-brane labels from the set of objects
$X$ of the groupoid $\cal{G}$.

\section*{Acknowledgments}

HP would like to thank Emmanuel College for a Research Fellowship.  HP
has performed a part of this work at Emmanuel College and at the
Department of Applied Mathematics and Theoretical Physics (DAMTP),
University of Cambridge, UK.  Both authors would like to thank John
Baez and Martin Hyland for helpful discussions and Ingo Runkel for
correspondence. We thank Ezra Getzler for bringing the
article~\cite{AN} to our attention. Both authors are grateful to the
European Union Superstring Theory Network for support.

\end{document}